\documentclass[11pt,reqno]{amsart}
\usepackage{amsmath,hyperref,graphicx}
\usepackage[utf8x]{inputenc}
\usepackage[T1]{fontenc}
\usepackage[english]{babel}
\usepackage{amsfonts}
\usepackage[export]{adjustbox}
\usepackage{amsthm}
\usepackage[dvipsnames]{xcolor}
\usepackage[text={6.5in, 8in}, centering]{geometry}
\usepackage{caption}
\usepackage{subcaption}

\newtheorem{thm}{Theorem}
\newtheorem*{thm*}{Theorem}

\usepackage{appendix}

\theoremstyle{definition}

\newcommand{\A}{\mathbf A}
\newcommand{\Aa}{\mathbf a}
\newcommand{\bb}{\mathbf b}

\newcommand{\nn}{\mathbf n}
\newcommand{\T}{\mathbf T}
\newcommand{\Tss}{\mathbf{T}_{ss}}

\newcommand{\tpq}{t_{pq}}

\newcommand{\pp}{\mathbf{e}_1}
\newcommand{\qq}{\mathbf{e}_2}
\newcommand{\Ts}{\T_s}
\newcommand{\Xt}{\mathbf X_t}
\newcommand{\X}{\mathbf X}
\newcommand{\Xs}{\mathbf X_s}
\newcommand{\Xss}{\mathbf X_{ss}}
\newcommand{\Tt}{\mathbf T_t}
\newcommand{\RR}{\mathbb R}
\newcommand{\Hbb}{\mathbb H}
\newcommand{\Hbf}{\mathbf H}
\newcommand{\Talg}{\mathbf T_{alg}}
\newcommand{\Tnum}{\mathbf T_{num}}
\newcommand{\Xalg}{\mathbf X_{alg}}
\newcommand{\Xnum}{\mathbf X_{num}}
\newcommand{\WW}{\widehat{W}}
\DeclareMathOperator{\arccosh}{arccosh}

\DeclareMathOperator{\sech}{sech}
\DeclareMathOperator{\sign}{sgn}

\newcommand{\mean}{\operatorname{mean}}

\title{Vortex Filament Equation for a Regular Polygon in the Hyperbolic Plane}
\author[F. de la Hoz]{Francisco de la Hoz}
\address[F. de la Hoz]{Department of Applied Mathematics and Statistics and Operations Research, Faculty of Science and Technology, University of the Basque Country UPV/EHU, Barrio Sarriena S/N, 48940 Leioa, Spain} 
\email{francisco.delahoz@ehu.eus}
\author[S. Kumar]{Sandeep Kumar}
\address[S. Kumar]{BCAM - Basque Center for Applied Mathematics, Alameda de Mazarredo 14, 48009 Bilbao, Spain.}
\email{skumar@bcamath.org}
\author[L. Vega]{Luis Vega}
\address[L. Vega]{BCAM, Department of Mathematics, Faculty of Science and Technology, University of the Basque Country UPV/EHU, Barrio Sarriena S/N, 48940 Leioa, Spain}
\email{luis.vega@ehu.eus, lvega@bcamath.org}

\date{\today}
\begin{document}
\newenvironment{red}{\textcolor{red}}

\begin{abstract} 
The aim of this paper is twofold. First, we show the evolution of the vortex filament equation (VFE) for a regular planar polygon in the hyperbolic space. Unlike in the Euclidean space, the planar polygon is open and both of its ends grow up exponentially, which makes the problem more challenging from a numerical point of view. However, using a finite difference scheme in space combined with a fourth-order Runge-Kutta method in time and fixed boundary conditions, we show that the numerical solution is in complete agreement with the one obtained by means of algebraic techniques. Second, as in the Euclidean case, we claim that, at infinitesimal times, the evolution of VFE for a planar polygon as the initial datum can be described as a superposition of several one-corner initial data. As a consequence, not only can we compute the speed of the center of mass of the planar polygon, but the relationship also allows us to compare the time evolution of any of its corners with the evolution in the Euclidean case. 
\end{abstract}

\maketitle
%\tableofcontents
% --------------------------------------------------------------------------------
\section{Introduction}
% --------------------------------------------------------------------------------
Consider the binormal flow
\begin{equation}
\label{eq:BF}
\Xt = \kappa \bb, 
\end{equation} 
where $t$ is the time, $\kappa$ the curvature, and $\bb$ the binormal component of the Frenet-Serret formulas. This equation first appeared in the work of Da Rios in 1906, as an approximation of the dynamics of a vortex filament (represented by $\X$) under Euler equations, and was later rederived by Arms and Hama in 1965 \cite{darios,arms}. This model is commonly known as the vortex filament equation (VFE). The flow, also called the localized induction approximation (LIA), can be expressed as
\begin{equation}
\label{eq:VFE}
\Xt = \Xs \wedge_{+} \Xss, 
\end{equation}
where $s$ is the arc-length parameter, and $\wedge_{+}$ is the usual cross product. The tangent vector $\T = \Xs$ satisfies 
\begin{equation}
\label{eq:SMP}
\Tt = \T \wedge_{+} \Tss, 
\end{equation}
and, during the time evolution, it preserves its magnitude, so we can assume that it takes values in the unit sphere, i.e., $\T\in\mathbb{S}^2$. Equation \eqref{eq:SMP} is called the Schr\"{o}dinger map equation onto the sphere and can be expressed in a more geometric way as 
\begin{equation}
\label{eq:SMP-geo}
\Tt = \mathbf J \mathbf{D}_s \Ts,
\end{equation} 
where $\mathbf{D}_s$ is the covariant derivative, and $\mathbf J$ is the complex structure of the sphere. By writing it in this way, \eqref{eq:SMP} can be extended to more general definition domains and images \cite{Khesin}. For instance, when the target space is chosen as the hyperbolic plane $\Hbb^2=\{ (x_1,x_2,x_3):-x_1^2+x_2^2+x_3^2=-1, x_1>0  \}$, i.e., a unit sphere in the Minkowski 3-space $\RR^{1,2} = \{(x_1,x_2,x_3): ds^2 = -dx_1^2 + dx_2^2 + dx_3^2\}$, the equivalent of \eqref{eq:SMP} is (see \cite{Ding1998})
\begin{equation}
\label{eq:SMP-hyp}
\Tt = \T \wedge_{-} \Tss, 
\end{equation}
and that of \eqref{eq:VFE} is
\begin{equation}
\label{eq:VFE-hyp}
\Xt = \Xs \wedge_{-} \Xss,
\end{equation}
where $\X\in \RR^{1,2}$, $\T\in \Hbb^2$, and the Minkowski cross product $\wedge_{-}$ is defined by (see \cite{HypGeoText1})
$$
\Aa \wedge_- \bb = ( - (a_2b_3 - a_3b_2), a_3b_1 - a_1b_3, a_1b_2-a_2b_1), \ \Aa, \bb \in \RR^{1,2}. 
$$
In this paper, we use the term \textit{hyperbolic} to refer to the case when $\T\in\Hbb^2$; and \textit{Euclidean}, when $\T\in\mathbb{S}^2$.

The Minkowski pseudo-scalar product is given by 
$$
\Aa \circ_- \bb = -a_1b_1 + a_2 b_2 + a_3b_3,
$$
which defines 
\begin{equation}
\label{eq:Mink-norm}
| \Aa |_0^2 = \Aa \circ_{-} \Aa.
\end{equation}
Thus, depending on whether $|\cdot|_0$ is positive, zero, or positive imaginary, the corresponding vector  can be classified as \textit{space-like}, \textit{light-like}, or \textit{time-like}, respectively. Since $\T\in\Hbb^2$, the corresponding $\X$ is called a time-like curve \cite{HypGeoText1,Lopez}. Note that depending on the sign of the first component of a time-like vector, it can be further classified as \textit{positive} or \textit{negative time-like}; for instance, in the definition of $\Hbb^2$ given above, we have considered only the positive time-like vectors. Let us also define the \textit{hyperbolic angle} between two positive (respectively, negative) time-like vectors $\Aa$ and $\bb$ as the unique nonnegative real number $\sigma(\Aa,\bb)$, such that
\begin{equation}
\label{eq:hyp-angle}
\Aa \circ_- \bb = -| \Aa |_0 | \bb |_0 \cosh(\sigma(\Aa,\bb)).
\end{equation}
In this work, we deal with vectors that are positive time-like; for simplicity of notation, we refer to them as time-like. On the other hand, for a sufficiently smooth curve $\X$ with curvature $\kappa$ and torsion $\tau$, the equivalent of the Frenet-Serret formulas in the hyperbolic setting is given by
\begin{equation}
\label{mat:SF-chap-hyp}
\begin{pmatrix} \T \\ \nn \\ \bb \end{pmatrix}_s = 
\begin{pmatrix} 0 & \kappa & 0\\ \kappa & 0 & \tau \\  0& -\tau & 0 \end{pmatrix} 
{.}\begin{pmatrix} \T \\ \nn \\ \bb \end{pmatrix},
\end{equation}
where the normal vector $\nn$ and binormal vector $\bb$ are space-like and, along with $\T$, form an orthonormal system  \cite{Lopez}. In addition, the corresponding filament function,
\begin{equation}
\label{eq:fila-fun-hyp}
\psi(s,t) = \kappa(s,t) e^{i\int_{0}^{s}\tau(s^\prime,t)ds^\prime},
\end{equation}
transforms \eqref{eq:SMP-hyp}--\eqref{eq:VFE-hyp} into the defocusing nonlinear Schr\"{o}dinger (NLS) equation \cite{hasimoto}:
\begin{align}\label{eq:NLS-hyp}
\psi_t = i\psi_{ss} - \frac{i}{2} \psi (|\psi|^2+A(t)), \ A(t)\in\RR.
\end{align}
As is well known, \eqref{eq:NLS-hyp} is a completely integrable system with infinitely many conservation laws. The simplest of these conservation laws is the one associated with the space $\mathcal L^2$. The rest of them involve an increasing number of derivatives of the solution with a jump of $1/2$ derivative from one law to the next one, if the regularity is measured using the class of Sobolev spaces. For each of these conservation laws, explicit solutions can be constructed. Moreover, the inverse scattering method can be used to build the solution for generic regular data. In this paper, we are motivated by a geometric problem and the possibility of having an initial condition with corners. Thanks to Hasimoto transformation, this implies considering initial data given by a sequence of Dirac delta functions, so that it belongs to the Sobolev space $\mathcal H^s$, with $s<-1/2$, and none of the conservation laws mentioned above can be used.
 
Note that VFE is time-reversible, i.e., if $\X(s,t)$ is a solution, then so is $\X(-s,-t)$. Bearing this in mind, an important property of VFE and hence, of the Schr\"{o}dinger map, is that it has a one-parameter family of regular self-similar solutions that develop a corner-shaped singularity in finite time. In other words, at the time of the formation of the singularity, i.e., $t=0$, the curve $\X$ has a corner, its tangent vector is a Heaviside-type function, and $\psi$ is a Dirac delta located at $s=0$. This was shown in \cite{GutierrezRivasVega2003} for the Euclidean case, and the hyperbolic case was studied in \cite{delahoz2007} (from now on, it will be referred to as the one-corner problem). Moreover, the well-posedness of the problem in the elliptic case has been established through a series of papers by Banica and Vega \cite{BV2,BV3}. 

On the other hand, the numerical study of the self-similar solutions was first done in \cite{buttke87}, and later in \cite{DelahozGarciaCerveraVega09}, where both the Euclidean and hyperbolic cases were considered. In \cite{DelahozGarciaCerveraVega09}, not only the formation of the singularity was captured, but the authors also started with a corner-shaped initial datum and recovered the self-similar solutions numerically. In all the cases, the choice of boundary conditions was found to be crucial. 

Although the problem of a curve with one corner that is otherwise smooth is well-understood both theoretically and numerically, the case of a polygonal curve has gained attention only recently \cite{JSt2012,JSt2015}. In \cite{HozVega2014}, a regular planar polygon with $M$ sides (which, from now on, will be referred to as the planar $M$-polygon) was considered as an initial curve in the Euclidean case; and using algebraic and numerical techniques, it was shown that the evolution of $\T$, and that of $\X$ after removing the vertical height, are $2\pi/M^2$-periodic in time. Moreover, at intermediate times that are rational multiples of $2\pi/M^2$, i.e., $\tpq=(2\pi/M^2)(p/q)$, with $\gcd(p,q)=1$, the planar $M$-polygon evolves in such a way that it has $Mq$ sides if $q$ is odd, and $Mq/2$ sides if $q$ is even, a behavior that is reminiscent of the so-called Talbot effect in optics \cite{BK,ET,Olver}. Let us also mention that at a macroscopic level, effects similar to those mentioned above were also observed in the case of real fluids \cite{GG,GGP}.

Another interesting aspect of the evolution of the planar $M$-polygon is the trajectory of any of its corners, which seems to be a multifractal and resembles the graph of Riemann's non-differentiable function \cite{Du}:
\begin{equation}
\sum_{k=1}^{\infty} \frac{\sin (\pi k^2)}{\pi k^2}, \ t \in [0,2].
\end{equation}
Recall that, for a given $M$, apart from the formation of new sides, the planar $M$-polygon evolves in the vertical direction with a constant speed $c_M$. Hence, bearing in mind the symmetries of the problem, the curve $\X(0,t)$ is planar. In \cite{HozVega2014}, it was denoted by $z_M(t)$, after removing the vertical height from it and projecting the resulting curve onto the complex plane. Then, strong numerical evidence was given, showing that, as $M$ tends to infinity, $z_M(t)$ converges to the complex version of Riemann's non-differentiable function: 
\begin{equation}
\phi(t) = \sum_{k=1}^{\infty} \frac{e^{\pi i k^2 t}}{i \pi k^2}, \ t\in[0,2].
\end{equation}
Recently, considering an $M$-sided polygon with nonzero torsion as the initial datum, new variants of $\phi(t)$ have been discovered in the trajectory of $\X(0,t)$, whose structure depends on the torsion introduced in the problem \cite{HozKumarVega2019}. Thus, by showing the existence of $\phi$ and its variants, it has been proved numerically that, in the Euclidean case, the time evolution of the smooth solutions of VFE, i.e., the circle, the helix and the straight line, is not stable. In other words, a particle can be placed on a curve arbitrarily close to a circle, helix or straight line, but, in the right topology, its trajectory converges to the graph of $\phi$ (or its variants). Moreover, this topology is motivated by some recent works on the well-posedness of VFE, which shows the existence of a new conservation law precisely at the critical scale of the problem $s=-1/2$, and that can be used for the solutions of the NLS equation associated to the self-similar solutions of VFE \cite{BanicaVega2018,BanicaVega2019,BanicaVega2020}.

Hence, we see that the evolution of $M$-sided polygons reveals many fascinating properties of VFE. With this motivation, another interesting problem is to look at the equivalent of a planar $M$-polygon in the hyperbolic setting and compare the evolutions of the two. It turns out that, in the absence of torsion, the corresponding polygon is a time-like curve that is characterized by a parameter $l>0$ representing the angle between any of its two sides. We refer to the polygonal curve as a \textit{planar $l$-polygon} (alternatively, in \cite{FF}, it is called an \textit{elementary $l$-convex polygon}).

Note that, unlike in the Euclidean case, the planar $l$-polygon is open and both of its endpoints tend to infinity (see Figure \ref{fig:l-polygon}). Furthermore, the corresponding tangent vector $\T$ lies on a unit hyperbola, and $\psi(s,0)$ is the $l$-periodic sum of Dirac deltas with coefficients that depend on the initial configuration of the planar $l$-polygon. Let us mention that, due to the mix of lack of regularity and periodicity, the well-posedness is quite challenging for this kind of problems. Recently, taking an initial datum consisting of polygonal lines that are asymptotically close to two straight lines at infinity, it has been proved that the problem is well-posed \cite{BanicaVega2018}. Moreover, using the appropriate topology, it has been shown in \cite{BanicaVega2019} that the solution also satisfies a conservation law.

The aim of this paper is twofold. First, we observe the evolution of \eqref{eq:SMP-hyp}--\eqref{eq:VFE-hyp} for a planar $l$-polygon as the initial datum, which, from now on, will be referred to as the $l$-polygon problem. In this regard, as in the Euclidean case, the algebraic solution is obtained by working at the level of the NLS equation. However, solving the problem numerically appears to be more challenging. In particular, in our numerical simulation, as we truncate the infinitely long $l$-polygon, the role of boundary conditions becomes very important. Moreover, as observed in the one-corner problem, due to the exponential growth of the tangent vector, working with all the values of the parameter $l$ becomes very difficult numerically. Bearing this in mind, we propose a numerical scheme (which will be explained in the following lines), and show a good agreement between the results thus obtained and the ones from the theoretical arguments. Then, as in \cite{HozVega2018}, we answer up to what extent the $l$-polygon problem and the one-corner problem are related. Consequently, not only can we compute the speed of the center of mass of the planar $l$-polygon, but the relationship also helps in comparing the trajectory of any of the corners of a regular planar polygon in both the Euclidean and hyperbolic cases.  

The structure of this paper is as follows. In Section \ref{sec:sol-l-pol-hyp}, we define the problem by formulating the main theoretical arguments that justify our numerical experiments. In particular, in Section \ref{sec:pb-def}, we introduce the parametric form of the initial data, and the relevant properties, such as symmetries. In Section \ref{sec:pb-formulation}, we observe that, as in the Euclidean case, the Galilean invariance of the NLS equation helps in obtaining the solution up to a function that depends on time. However, the function is now determined using the conservation law established for polygonal lines in \cite{BanicaVega2019}, an approach that was also employed in \cite{HozKumarVega2019}. Let us not forget that, in the case of curves with vanishing curvature, it is desirable to work with the parallel frame where the normal plane is spanned by the vectors $\pp$, $\qq$, whose space derivatives depend only on $\T$ \cite{Bishop}. In the hyperbolic setting, the corresponding parallel frame is given by \eqref{mat:SF-chap-hyp-thm}, where $\mathbf{e}_1$, $\mathbf{e}_2$ are the unit space-like normal and binormal vectors, respectively. Thus, by integrating the generalized Frenet-Serret formulas at times that are rational multiples of $l^2/(2\pi)$, we obtain the evolution of the curve $\X$ and of the tangent vector $\T$, up to a rigid movement. This has been illustrated in Section \ref{sec:alg-sol}, where knowing the rotations in the Minkowski 3-space (from now on, referred to as \textit{hyperbolic rotations}) is found to be quite essential \cite{MOME}. Moreover, the rigid movement can be determined by using the symmetries of the regular planar $l$-polygon and, in this way, we recover $\T$ completely. However, $\X$ is computed only up to a movement in the YZ-plane, which is obtained numerically in Section \ref{sec:num-res}. 

In Section \ref{sec:num-sol}, we study the numerical evolution of \eqref{eq:SMP-hyp}--\eqref{eq:VFE-hyp} for different values of the parameter $l$. Bearing in mind that, unlike in the Euclidean case, a planar $l$-polygon is of infinite length, we consider a planar $l$-polygon with only $M$ sides in our numerical simulations, i.e., such that its length is $L=lM$. We have found that Dirichlet boundary conditions on the tangent vector, with a finite difference discretization in space, combined with a fourth-order Runge--Kutta method in time, give the best numerical results, both in terms of computational cost and accuracy. These ideas have been offered in Section \ref{sec:num-mthd}. In Section \ref{sec:num-res}, we begin by calculating the movement of the center of mass in the YZ-plane, which allows us to compare the numerical solution with its algebraic counterpart (obtained in Section \ref{sec:alg-sol}). On the other hand, the trajectory of a corner of the $l$-polygon initially located at $s=0$ (i.e., $\X(0,t)$), although resembling Riemann's non-differentiable function, is quite different from its equivalent in the Euclidean case. Moreover, it converges to the function, as the parameter $l$ tends to zero. In Section \ref{sec:X0t}, we provide strong numerical evidence to prove this claim. Section \ref{sec:T-irr} is about the behavior of the tangent vector $\T$ near irrational times and its comparison with the tangent vector in the Euclidean case. 

Section \ref{sec:rel-l-1-corner} is based on the relationship between the $l$-polygon problem and the one-corner problem. In this regard, let us first briefly recall the main ideas of the one-corner problem. In \cite{delahoz2007}, the existence of the solutions of \eqref{eq:VFE-hyp} for the following initial datum are proved:
\begin{equation}
\label{eq:IVP-1corner}
\X(s,0) = \A^- s \chi_{(-\infty,0]}(s) + \A^+ s \chi_{[0,\infty)}(s), \ \A^\pm \in \Hbb^2,
\end{equation}
where, due to the rotation invariance of VFE, the unit vector can be chosen such that $\A^{\pm}=(A_1,\pm A_2,\pm A_3)^T$. The self-similar solutions of \eqref{eq:VFE-hyp} satisfying $\X(s,t) = \sqrt{t} \X(s/\sqrt{t},1)$, $t>0$, solve (see \cite{delahoz2007,buttke87})
\begin{align}
\label{eq:G-eqn}
\frac 1 2 \X(s/\sqrt{t},1) &- \frac{s}{2\sqrt{t}} \X^\prime(s/\sqrt{t},1) = \mathbf \X^\prime(s/\sqrt{t},1)  \wedge_{-} \mathbf \X^{\prime\prime}(s/\sqrt{t},1).
\end{align}
Then, from \eqref{mat:SF-chap-hyp} and \eqref{eq:G-eqn}, $\kappa(s,t) = c_0 / \sqrt{t}$ and $\tau(s,t)=s/(2t)$ can be obtained, where the constant $c_0$ characterizes the one-parameter family of smooth curves $\X$ that can be described using \eqref{mat:SF-chap-hyp}, $\Xs=\T$, and the initial conditions
\begin{equation}
\label{eq:ini-cond-1-crnr}
\begin{aligned}
\X(0,t) &= 2c_0 \sqrt{t}(0,0,1)^T, \\ 
\T(0,t) = (1,0,0)^T, \
\nn(0,t) &= (0,1,0)^T, \
\bb(0,t) = (0,0,1)^T.
\end{aligned}
\end{equation}
The parameter $c_0$ is the curvature of $\X(s,1)$, which, in turn, solves the following ODE \cite{delahoz2007}:
\begin{equation}
\label{eq:X-ODE}
\X^{\prime\prime\prime}(s,1) + \left(-c_0^2+\frac{s^2}{4}\right) \X^\prime(s,1) - \frac{s}{4} \X(s,1) = 0.
\end{equation}
With some abuse of notation, if we define the Fourier transform of $\X(s,1)$ by
\begin{equation*}
\hat \X(\xi) = \frac{1}{\sqrt{2\pi}} \int_{-\infty}^{\infty} \X(s,1) e^{-is\xi} ds, 
\end{equation*}
then it satisfies
\begin{equation}
\label{eq:X-ODE-FT}
\xi \hat{\X}^{\prime\prime}(\xi) + 3 \hat{\X}^\prime(\xi) + 4 \xi^3 \hat{\X}(\xi) + 4c_0^2 \xi \hat{\X}(\xi) = 0.
\end{equation}
That being said, in Section \ref{sec:num-exp-l-1-rel}, following the approach in \cite{HozVega2018}, we provide very strong numerical evidence to establish the connection between the two problems. As a consequence, in Section \ref{sec:frthr-rmrks}, an explicit expression for the speed of the center of mass of the planar $l$-polygon is given, according to which it moves in the vertical direction. Moreover, we also make some remarks on the trajectory of $\X(0,t)$. 

In Section \ref{sec:conclusion}, we discuss the main conclusions. Finally, recall that in \cite{delahoz2007}, a precise expression for the first component of the tangent vector $\A^\pm$ was obtained:
\begin{equation}
\label{eq:A1-c0}
A_1 = e^{\pi c_0^2/2},
\end{equation} 
which also relates $c_0$ to the time-like angle $\theta$ between $\A^+$ and $\A^-$:
\begin{equation}
\label{eq:A1-c0-ang}
\cosh(\theta) = -1 + 2A_1^2= -1 + 2 e^{\pi c_0^2}. 
\end{equation}
\noindent Thus, to conclude this paper, in Appendix \ref{apdx}, we provide calculations to obtain a compact expression for $A_1$, $A_2$, and $A_3$.

\section{A solution of $\Xt = \Xs \wedge_{-} \Xss$ for a planar $l$-polygon}
\label{sec:sol-l-pol-hyp}
% --------------------------------------------------------------------------------
One of the main goals of this paper is to obtain the solutions of \eqref{eq:VFE-hyp} and explain their dynamics, when regular planar $l$-polygons are considered as initial data. In this regard, by assuming uniqueness as in the Euclidean case, we prove the following theorem.
\begin{thm}
\label{thm:evo-l-pol}
	Assume that there exists a unique solution of the initial value problem
	\begin{equation}
	\Xt = \Xs \wedge_{-} \Xss,
	\end{equation}
	with $\X(s,0)$ being a regular planar $l$-polygon. Then, at a time $\tpq$ which is a rational multiple of $l^2/2\pi$, i.e., $\tpq\equiv(l^2/2\pi)(p/q)$, with  $p\in\mathbb{Z}$, $q\in \mathbb{N}$, $\gcd(p,q)=1$, the solution is a skew $l_q$-polygon, such that, in $s\in[d, d+l)$, for any $d\in\mathbb R$, $\X(s,t_{pq})$ has $q$ times as many sides (if $q$ odd) or $q/2$ times as many sides (if $q$ even) as $\X(s, 0)$. All the new sides have the same length, and the time-like angle $l_q$ between any two adjacent sides is constant. Moreover, the polygon at a time $\tpq$ is the solution of the generalized Frenet-Serret formulas
	\begin{equation}
	\label{mat:SF-chap-hyp-thm}
	\begin{pmatrix} \T(s,\tpq) \\ \pp(s,\tpq) \\ \qq(s,\tpq) \end{pmatrix}_s = 
	\begin{pmatrix} 0 & \alpha(s,\tpq) & \beta(s,\tpq)\\  \alpha(s,\tpq) & 0 & 0 \\ \beta(s,\tpq) & 0 & 0 \end{pmatrix} 
	{.}\begin{pmatrix} \T \\ \pp \\ \qq \end{pmatrix},
	\end{equation}
	where $\alpha(s,\tpq)+i\beta(s,\tpq)=\Psi(s,\tpq)$, and $\Psi(s,\tpq)$ is the $l$-periodic function defined over the first period $s\in[0,l)$ as
	\begin{equation}
	\Psi(s,\tpq)=
	\begin{cases}
	\displaystyle
	\frac{l_q}{\sqrt{q}} \sum_{m=0}^{q-1} G(-p,m,q) \delta(s-\tfrac{lm}{q}),  \ &\text{if $q$ odd,} \\[1em]
	\displaystyle
	\frac{l_q}{\sqrt{2q}} \sum_{m=0}^{q-1} G(-p,m,q) \delta(s-\tfrac{lm}{q}),  \ &\text{if $q$ even,} \\    
	\end{cases}
	\end{equation}
	with 
	$$
	G(a,b,c) = \sum_{n=0}^{c-1} e^{2\pi i(an^2+bn)/c}, \ a,b \in \mathbb{Z}, c\in \mathbb{Z} \backslash \{0\}
	$$
	being a generalized quadratic Gau{\ss} sum. The mutual time-like angle $l_q$ between any two sides of the new polygon is given by 
	\begin{equation}\label{eq:l_q-thm}
	l_q =
	\begin{cases}
	2\arccosh(\cosh^{1/q}(l/2)), \quad & \mbox{if $q$  odd}, \\
	2\arccosh(\cosh^{2/q}(l/2)), \quad & \mbox{if $q$  even}.
	\end{cases}
	\end{equation}
\end{thm}

Remark that both $\X(s, t_{pq})$ and $\X(s, 0)$ have obviously a countable infinite number of sides. Therefore,  in this paper, whenever we say loosely speaking that $\X(s, t_{pq})$ has $q$ times as many sides as $\X(s, 0)$, etc., it must be understood that $s$ is being taken over any half-open interval of length $l$. On the other hand, let us mention that the determination of $l_q$ follows from the conservation law established for polygonal lines in \cite{BanicaVega2019}. 

% --------------------------------------------------------------------------------
\subsection{Problem definition}
\label{sec:pb-def}
\begin{figure}[htbp!] \centering
	\includegraphics[width=0.5\textwidth, clip=true]{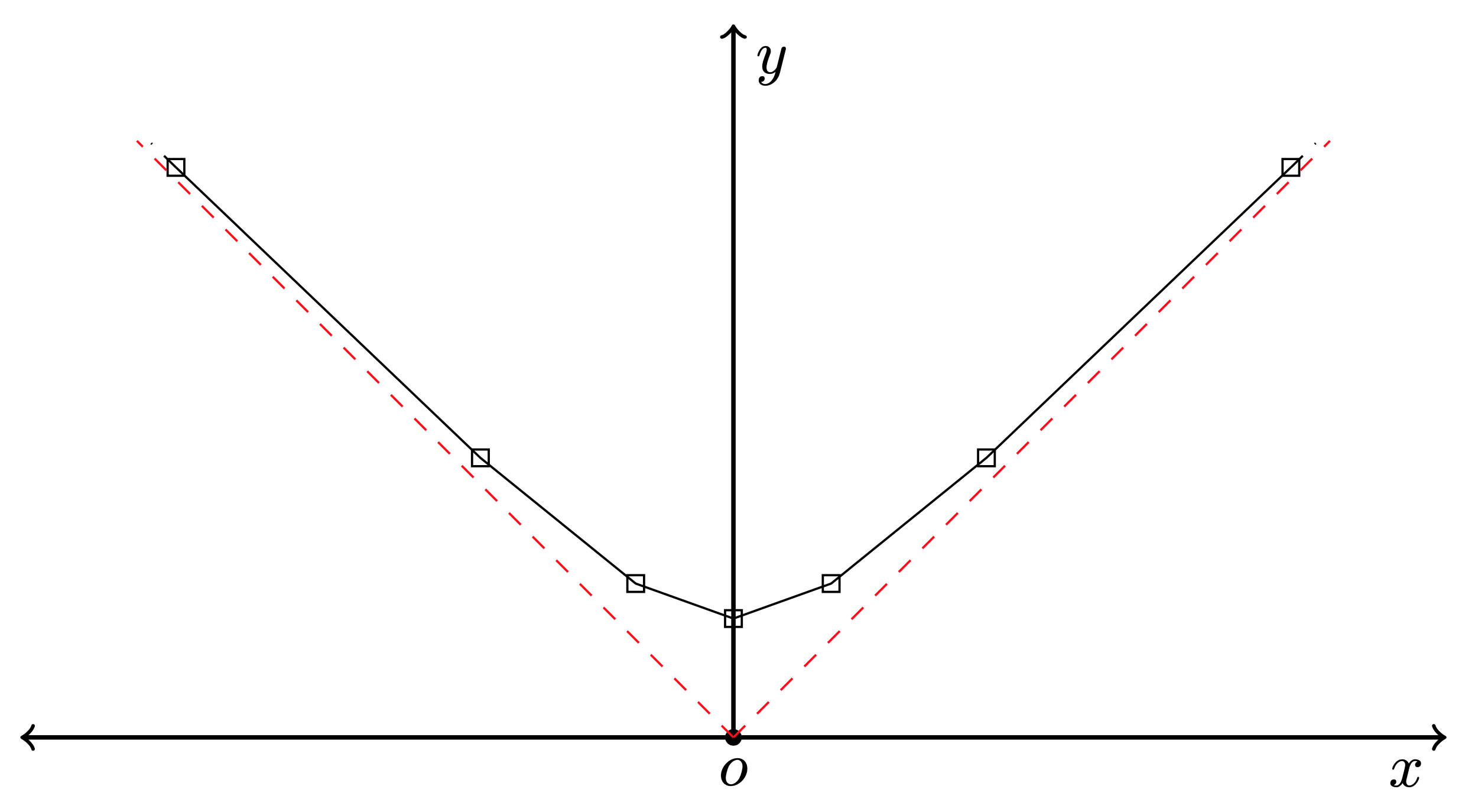}
	\caption{A planar $l$-polygon with vertices located at $s_n=nl$, $n\in \mathbb{Z}$, (black squares) and the asymptotes (dotted red lines).}
	\label{fig:l-polygon}  
\end{figure}

% --------------------------------------------------------------------------------
Given the parameter $l>0$, an arc-length parameterized planar $l$-polygon can be understood as a curve with curvature given by
\begin{equation} 
\label{eq:cur-ini-pol-hyp}
\kappa(s) = c_0 \sum_{k=-\infty}^{\infty} \delta(s-lk), \ s\in \RR.
\end{equation}
Here, the vanishing argument of the equally spaced Dirac deltas corresponds to the location of the corners, and the coefficient $c_0>0$ depends on the initial configuration of the curve. In particular, bearing in mind \eqref{eq:A1-c0-ang}, we choose 
\begin{equation}\label{eq:c0_exp_hyp}
c_0 = \left[ \frac{2}{\pi}\ln\left(\cosh\left(\frac{l}{2}\right)\right)\right]^{1/2}.
\end{equation}
Note that, in the absence of torsion, from \eqref{eq:fila-fun-hyp}, $\psi(s,0)$ is the curvature of the initial polygonal curve, i.e., $\psi(s,0)=\kappa(s)$, which is $l$-periodic. Moreover, since \eqref{eq:SMP-hyp}--\eqref{eq:VFE-hyp} are invariant under hyperbolic rotations, we can assume without loss of generality that the corresponding initial planar polygonal curve $\X(s,0)$ and its tangent vector $\T(s,0)$ lie on the XY-plane. Thus, by denoting the plane by OXY, for $ s_n=nl,\ n\in\mathbb{Z}$, we write the piecewise constant tangent vector $\T\equiv(T_1,T_2,T_3)^T$ as 
\begin{equation}\label{eq:T-ini-hyp}
\T(s,0) =  \left(\cosh \left(l/2+s_n \right), \sinh \left(l/2+s_n \right),0\right)^T, \ s\in(s_n,s_{n+1}).
\end{equation}
As a result, the vertices of the planar $l$-polygon $\X\equiv(X_1,X_2,X_3)^T$  can be expressed as
\begin{equation} \label{eq:X-ini-hyp}
\X(s_n,0) = \frac{(l/2)}{\sinh(l/2)} \left( \sinh\left(s_n \right), \cosh\left(s_n \right),0\right)^T,  
\end{equation} 
and, for $s\in(s_n,s_{n+1})$, the point $\X(s,0)$ lies in the segment that joins $\X(s_n,0)$ and $\X(s_{n+1},0)$ (see Figure \ref{fig:l-polygon}). Note that we have chosen $\T(s,0)$ as \eqref{eq:T-ini-hyp}, so that the vertex corresponding to $\X(0,0)$ in \eqref{eq:X-ini-hyp} lies on the $y$-axis. 

Thus, $\X(s,0)$ is a curve with infinite length and, from \eqref{eq:hyp-angle}, it follows that the hyperbolic angle between any two of its sides is constant and equal to $l$. Since one of the main concerns in this paper is to address the numerical evolution of a planar $l$-polygon, we work with a truncated curve with a finite number of sides. 

%--------------------------------------------------------------------------------
\subsubsection{Spatial symmetries of $\X$ and $\T$}
\label{sec:Spatial-sym}
% --------------------------------------------------------------------------------
The invariance of \eqref{eq:SMP-hyp}--\eqref{eq:VFE-hyp} under hyperbolic rotations follows from the invariance of the Minkowski cross product under them \cite{MOME}. Thus, given a hyperbolic rotation matrix $\mathbf{R}$, such that $\mathbf{R}\cdot \T(s,0)=\T(s,0)$ and $\mathbf{R}\cdot\X(s,0)=\X(s,0)$, if the solution is unique, then $\mathbf{R}\cdot \X(s,t)=\X(s,t)$, $\mathbf{R}\cdot \T(s,t)=\T(s,t)$, for all $t$. In particular, since $\X(s,0)$ and $\T(s,0)$, which are given respectively by \eqref{eq:X-ini-hyp} and \eqref{eq:T-ini-hyp}, are invariant under a rotation of time-like angle $nl$ about a space-like $z$-axis for all $n\in\mathbb{Z}$, it can be concluded that $\X(s,t)$ and $\T(s,t)$ are invariant under hyperbolic rotations, for all $t$.

One important consequence of these symmetries is that, for any time $t$, $\X(s+nl,t)$ always lies in the same orthogonal plane to the $z$-axis. Furthermore, as in the Euclidean case, \eqref{eq:SMP-hyp}--\eqref{eq:VFE-hyp} are mirror invariant, and, consequently, $\X(s,t)-\X(-s,t)$ is a positive multiple of $(1,0,0)^T$. This property plays an important role when constructing the algebraic solution.   
% --------------------------------------------------------------------------------
\subsection{Problem formulation and the behavior at rational multiples of the time period}
\label{sec:pb-formulation}
% --------------------------------------------------------------------------------
First, let us mention that, at the level of the NLS equation, the hyperbolic case is not much different from the Euclidean case; however, the obtention of $\X$ and $\T$ depends entirely on hyperbolic rotations \cite{MOME,HypGeoText1}. In this regard, following the approach in \cite{HozVega2014}, we observe that, by definition, $\psi(s,0)$ is $l$-periodic, and, since \eqref{eq:NLS-hyp} is invariant with respect to space translations, $\psi(s,t)$ is also $l$-periodic, for all $t \in \mathbb{R}$. On the other hand, $\psi(s,0)=e^{irks}\psi(s,0)$, $r=2\pi/l$, $l>0$; thus, from the Galilean invariances of \eqref{eq:NLS-hyp}, $\psi(s,t) = e^{irks - i (rk)^2 t} \psi(s-2rkt,t)$, for all $k$. Furthermore, since $\psi$ is periodic, using its Fourier coefficients, it can be expressed as
\begin{equation}
\label{eq:psi_st_hyp}
\psi(s,t) = \hat{\psi}(0,t) \sum_{k=-\infty}^{\infty}e^{i(rk)^2t+i(rk)s} ,
\end{equation}
where $\hat{\psi}(0,t)$ is a constant depending on time $t$. Due to the gauge invariance, we can take it to be real (see \cite{HozVega2014}), and
its value is computed explicitly by using a conservation law that will be explained in the following lines. Remark that $\psi(s,t)$ is periodic in time with period $2\pi/r^2$, or, $l^2/2\pi$, which we denote by $T_f$ in this paper.

Next, evaluating \eqref{eq:psi_st_hyp} at rational multiples of the time period $T_f$, i.e., at $t=t_{pq}= \frac{2\pi}{r^2}\frac{p}{q}$, $p\in\mathbb{Z}, q \in \mathbb{N}$, $\gcd(p,q)=1$, gives (See \cite[Section 3.3]{HozVega2014} for the intermediate steps)
\begin{equation}\label{eq:psi_tpqfinal}
\psi(s,t_{pq})	=
\begin{cases}
\displaystyle \frac{l}{\sqrt{q}} \hat{\psi}(0,t_{pq}) \sum_{m=0}^{q-1} e^{i\theta_m} \delta\left(s-\tfrac{ml}{q}\right), \  &\mbox{if $q$ odd},\\[0.5em]
\displaystyle \frac{l}{\sqrt{q/2}} \hat{\psi}(0,t_{pq})  \sum_{m=0}^{q/2-1} e^{i\theta_{2m+1}} \delta \left(s-\tfrac{({2m+1})l}{q}\right), \  &\mbox{if $q/2$ odd},\\[0.5em]
\displaystyle \frac{l}{\sqrt{q/2}} \hat{\psi}(0,t_{pq}) \sum_{m=0}^{q/2-1} e^{i\theta_{2m}} \delta\left(s-\tfrac{2ml}{q}\right), \  &\mbox{if $q/2$ even},
\end{cases}
\end{equation}
for $s\in(0,l)$. This implies that, at any rational time $\tpq$, a single side of the $l$-polygon at $t=0$ will evolve into $q$ sides, if $q$ is odd, and $q/2$ sides, if $q$ is even. Since it holds true for any $k\in\mathbb{Z}$, this would imply that the resulting polygon will have $q$ or $q/2$ times as many sides as the initial $l$-polygon. The new Dirac deltas thus formed are equally spaced and, as a result, all the sides of the new polygon are of equal length. Furthermore, the coefficients of Dirac deltas have equal modulus and are given by 
\begin{equation*}
c_q=
\begin{cases}
\frac{l}{\sqrt{q}}\hat{\psi}(0,\tpq), \ &\mbox{if $q$ is odd},\\[0.5em]
\frac{l}{\sqrt{q/2}}\hat{\psi}(0,\tpq), \ &\mbox{if $q$ is even}.
\end{cases}
\end{equation*}
Note that the conservation law established for the polygonal lines in \cite{BanicaVega2019} holds true for both the focusing and the defocusing NLS equation. Therefore, by following the approach in \cite{HozKumarVega2019}, we obtain $c_q = c_0/\sqrt{q}$, if $q$ is odd, and $c_q = c_0/\sqrt{2q}$, if $q$ is even, and 
\begin{equation}
\label{eq:psi-hat-0tpq}
\hat \psi(0,\tpq)=c_0/l. 
\end{equation}
On the other hand, \eqref{eq:c0_exp_hyp} holds true whenever a corner is created; for instance, in our case, at rational times $\tpq$. Then, from \eqref{eq:A1-c0-ang} and denoting the time-like angle between any two tangent vectors by $l_q$, 
\begin{equation}
\label{eq:cqq}
\cosh\left(\tfrac{l_q}{2}\right) = e^{\pi c_q^2/2}. 
\end{equation}
Moreover, since $c_q$ is independent from $k$, the angle $l_q$ is the same for all sides and, thus, using \eqref{eq:c0_exp_hyp}, \eqref{eq:psi-hat-0tpq}, \eqref{eq:cqq}  it can be expressed as
\begin{equation}\label{eq:l_q}
l_q =
\begin{cases}
2\arccosh(\cosh^{1/q}(l/2)), \quad & \mbox{if $q$  odd}, \\
2\arccosh(\cosh^{2/q}(l/2)), \quad & \mbox{if $q$  even}.
\end{cases}
\end{equation}  
% --------------------------------------------------------------------------------
\subsection{Algebraic solution}
\label{sec:alg-sol}
% --------------------------------------------------------------------------------
In order to construct the algebraic solution, as in \cite{HozVega2014}, we integrate the Frenet-Serret formulas \eqref{mat:SF-chap-hyp-thm}, taking
\begin{equation}
\label{eq:Psi-def}
\Psi(s,\tpq) = \frac{l_q}{c_q} \psi(s,\tpq) = \alpha(s,\tpq)+i\beta(s,\tpq),
\end{equation}
for $q$ odd, and similarly for $q$ even. 
Then, by expressing $\alpha  + i \beta  = l_q e^{i\theta}$, the integration yields
\begin{equation} \label{mat:Hm}
\mathbf{H}=\begin{pmatrix} \cosh(l_q) & \cos(\theta)  \sinh(l_q) & \sin(\theta)  \sinh(l_q) \\ \cos(\theta)  \sinh(l_q) & 1+\cos^2(\theta)(\cosh(l_q)-1) & \sin(\theta) \cos(\theta)(\cosh(l_q)-1)\\ \sin(\theta)  \sinh(l_q) & \sin(\theta) \cos(\theta) (\cosh(l_q) - 1) & 1+\sin^2(\theta)(\cosh(l_q)-1) \end{pmatrix},
\end{equation}
which is a hyperbolic rotation of angle $l_q$ about a space-like axis $(0,-\sin (\theta),\cos (\theta))^T$ \cite{MOME}. In other words, $\Hbf$ describes the transition from a vertex located at $s_k=-L/2+k(l/q)$, $k=0,1,\ldots,Mq-1$. By choosing the basis vectors $\tilde \T(s)$, $\tilde{\mathbf{e}}_1(s)$, $\tilde{\mathbf{e}}_2(s)$, such that they form an identity matrix at $s=s_0^-$, we obtain their values for the remaining $Mq$ sides by a subsequent action of $\Hbf$ corresponding to $\Psi(s,\tpq)$.  Additionally, $\tilde{\X}$, i.e., $\X$ up to a rigid movement, can be computed from $\tilde{\T}$ through
\begin{equation} \label{eq:X_alg_int}
\tilde{\X}(s_{k+1}) = 	\tilde{\X}(s_k) + \tfrac{l}{q} 	\tilde{\T}(s_k^+), \ k=0,1,\ldots,Mq,
\end{equation}
where $\tilde{\X}(s_0)$ can be assigned any value, for example, $\tilde{\X}(s_0)=(0,0,0)^T$.

Next, we determine the correct rotation by using the symmetries of the regular planar $l$-polygon. In order to align the polygon orthogonal to the $z$-axis, we use the fact that, at any time $t$, $\X(lk)$, lies in the XY-plane, for $k\in\mathbb{Z}$. Then, the resulting curve is rotated about the $z$-axis in such a way that $\X(l)-\X(-l)$ is a positive multiple of $(1,0,0)^T$. This can be done efficiently in the following way:
\begin{enumerate}
	\item Compute the unit time-like vectors $\mathbf{w}^+ = \tfrac{\tilde{\X}(l)-\tilde{\X}(0)}{|\tilde{\X}(l)-\tilde{\X}(0)|_0}$, $\mathbf{w}^- = \tfrac{\tilde{\X}(-l)-\tilde{\X}(0)}{|\tilde{\X}(-l)-\tilde{\X}(0)|_0}$.
	\item Compute the unit space-like vector $\hat{\mathbf{u}} = \tfrac{ \mathbf{w}^+ \wedge_{-} \mathbf{w}^-}{|\mathbf{w}^+ \wedge_{-} \mathbf{w}^-|_0}$.
	\item If the space-like vectors $\hat{\mathbf{u}}$ and $\hat{\mathbf{z}}=(0,0,1)^T$ are such that (see \cite{HypGeoText1})
	 \begin{enumerate}
	 	\item $(\hat{\mathbf{u}})_3 > | \hat{\mathbf{u}} |_0$, then the time-like angle $\nu_1 = \arccosh\left((\hat{\mathbf{u}})_3/|\mathbf{u}|_0\right)$, and $\hat{\mathbf{v}} = \hat{\mathbf{u}} \wedge_{-} \hat{\mathbf{z}}$ is a space-like vector,
	 	\item $(\hat{\mathbf{u}})_3 < | \hat{\mathbf{u}} |_0$, then the angle $\nu_1 = \arccos\left((\hat{\mathbf{u}})_3/|\mathbf{u}|_0\right)$, and $\hat{\mathbf{v}} = \hat{\mathbf{u}} \wedge_{-} \hat{\mathbf{z}}$ is a time-like vector,  
	 	\item  $(\hat{\mathbf{u}})_3 = | \hat{\mathbf{u}} |_0$, then $\nu_1=0$, and $\mathbf{L}_1$ is an identity matrix, then, $\mathbf{L}_1$ is a matrix performing a rotation of an angle $\nu_1$ about the axis $\hat{\mathbf{v}}/|\hat{\mathbf{v}} |_0$ \cite{MOME}.
	 	
	 \end{enumerate}
	\item Compute time-like vectors $\mathbf{w}_{rot}^+ = \mathbf{L}_1 \cdot \mathbf{w}^+$, $\mathbf{w}_{rot}^- = \mathbf{L}_1 \cdot \mathbf{w}^-$, and $\mathbf{w} =   \frac{\mathbf{w}_{rot}^+ - \mathbf{w}_{rot}^-}{| \mathbf{w}_{rot}^+ - \mathbf{w}_{rot}^-|_0}$. Then, $\nu_2 = \arccosh(\mathbf{w}\circ_{-}(1,0,0)^T)$ is the time-like angle, and $\mathbf{L}_2$ is the corresponding rotation about the axis given by $\tfrac{\mathbf{w} \wedge_{-} (1,0,0)}{|\mathbf{w} \wedge_{-} (1,0,0)|_{0}}$.

	\item Compute the desired rotation $\mathbf{L} = \mathbf{L}_2 \cdot \mathbf{L}_1,$ and $\T = \mathbf{L}\cdot \tilde{\T},$ $\X = \mathbf{L}\cdot \tilde{\X}$.  
\end{enumerate}
Thus, we obtain $\X$ and $\T$ correctly oriented. Although, the computation of $\T$ is complete, in order to fully determine $\X$, we need to compute the movement of its center of mass, which is done in the next section. Finally, from \eqref{eq:psi_tpqfinal}, \eqref{eq:psi-hat-0tpq} and \eqref{eq:Psi-def}, we conclude the proof of Theorem \ref{thm:evo-l-pol}.

% --------------------------------------------------------------------------------
\section{Numerical solution}
\label{sec:num-sol}
% --------------------------------------------------------------------------------
As mentioned previously, in order to simulate numerically the evolution, we consider an $l$-polygon of length $L$ that is now characterized by two parameters $l$ and $M$, such that $L=l \cdot M$. For our purposes, we take $M$ even, so that the initial curve $\X(s,0)$, $s\in[-L/2,L/2]$, has a vertex located at $s=0$ and the symmetries described in Section \ref{sec:Spatial-sym} apply. This also allows us to capture the time evolution of a corner initially located at $s=0$, i.e., $\X(0,t)$. Remark that $M$ is finite, but we are approximating an infinitely long polygon; so, in principle, more accurate results would be expected with a larger value of $M$. However, both $M$ and $l$ cannot be large, since due to the exponential growth of the Euclidean norm of the tangent vector $\T$, for a fixed $M$, a large value of $l$ causes the solution to blow up in a short time, making the numerical scheme unstable. This was also observed in the one-corner problem, where large values of $c_0$ lead to similar effects \cite{DelahozGarciaCerveraVega09}. On the other hand, a large $M$ value forces us to consider only small values of $l$. Let us not forget that, as we work with a truncated $l$-polygon, the role of the boundary conditions also becomes very important.

Our goal is to solve \eqref{eq:SMP-hyp}--\eqref{eq:VFE-hyp} numerically for the initial data given by \eqref{eq:T-ini-hyp}--\eqref{eq:X-ini-hyp}, for $s\in[-L/2,L/2]$. There have been several papers dedicated to the numerical treatment of \eqref{eq:SMP-hyp}--\eqref{eq:VFE-hyp} \cite{buttke87,DelahozGarciaCerveraVega09,HozVega2014}. For instance, for the Euclidean regular $M$-polygons, the coupled system is solved with a pseudo-spectral method in space and a fourth-order Runge--Kutta method in time \cite{HozVega2014,HozKumarVega2019}. In our case, a Chebyshev spectral discretization with an explicit scheme in time poses a severe restriction $|\Delta t|=\mathcal{O}(1/N^4)$, where $N$ is the number of nodes. On the other hand, due to its low order of accuracy, a second-order semi-implicit backward difference formula applied on the stereographic projection of \eqref{eq:SMP-hyp} does not serve our purpose, as we are interested in the evolution for all rational times, unlike in \cite{DelahozGarciaCerveraVega09}. Thus, after trying several numerical methods, we have found that both in terms of efficiency and computational cost, a fourth-order finite difference discretization in space with a fourth-order Runge--Kutta method in time and with fixed boundary conditions on $\T$ yield the best results.        

% --------------------------------------------------------------------------------
\subsection{Numerical method}
\label{sec:num-mthd}
% --------------------------------------------------------------------------------
We divide the interval $[-L/2,L/2]$ into $N+1$ equally spaced nodes $s_j = -L/2+ j L / N, \ j=0,1,\ldots,N$, with a step size $\Delta s=L/N$. The time interval $[0,T_{f}]$ has been discretized into $N_t+1$ equally spaced time steps $t_n=n \Delta t$, $n=0,1,\ldots,N_t$, with $\Delta t=T_f/N_t$. We denote $\X_j^{(n)} \equiv \X^{(n)}(s_j)\equiv \X(s_j,t_n)$, where $\X_j^{(0)}$ can be computed from \eqref{eq:X-ini-hyp} by using linear interpolation, and $\T_j^{(n)} \equiv \T^{(n)}(s_j)\equiv \T(s_j,t_n)$, where $\T(s_j,\cdot)=\T(s,\cdot)$, for $s_j\leq s<s_{j+1}$, if $s<0$, and $s_j< s\leq s_{j+1}$, if $s>0$. Thus, we obtain $N$ values of the piecewise constant tangent vector, each corresponding to $N$ segments, respectively.

In order to approximate the first and second derivatives, we use a fourth-order central difference scheme for the inner points, and in order to keep the same order of accuracy over the whole discretized domain, we employ a fourth-order forward/backward difference scheme for the boundary and its neighboring points; this results in banded differentiation matrices of size $(N+1)\times (N+1)$. Let us mention that, in order to maintain the dimensions of the vectors $\X_j^{(n)}$ and $\T_j^{(n)}$ consistent, we obtain the $N+1$ values of the piecewise continuous tangent vector $\T_j^{(0)}$ in the following way:
\begin{align*}
\tilde \T_0^{(0)} &= \T_0^{(0)}, \
\tilde \T_{j+1}^{(0)} = ( \T_j^{(0)} +  \T_{j+1}^{(0)})/2, \ j=0,1,\ldots,N-2, \\
\tilde \T_{N}^{(0)} &= \T_{N-1}^{(0)}, \
\T_j^{(0)} =  \tilde \T_j / | \tilde \T_j  |_0.
\end{align*}  
Hence, by fixing the boundary conditions for the tangent vector $\T$, which can be introduced explicitly, we solve the following initial-boundary value problem:
\begin{equation}\label{eq:fixedBCT}
\begin{cases}
\Tt(s,t) = \T(s,t) \wedge_{-} \Tss(s,t), \\
\Xt(s,t) = \Xs(s,t) \wedge_{-} \Xss(s,t) = \T(s,t) \wedge_{-} \Ts(s,t),  \\
\T(-L/2,t) = \left(\cosh \left(l/2-L/2 \right), \sinh \left(l/2-L/2 \right),0\right)^T, \\ \T(+L/2,t) = \left(\cosh \left(l/2+L/2 \right), \sinh \left(l/2+L/2 \right),0\right)^T, \ t \in [0,T_f], \\
\end{cases}
\end{equation}
with initial conditions $\X(s,0)$, $\T(s,0)$ given by \eqref{eq:T-ini-hyp}, \eqref{eq:X-ini-hyp}, respectively.
By using the space discretization mentioned above, we integrate \eqref{eq:fixedBCT} numerically by means of a fourth-order Runge--Kutta method in time. Moreover, in the numerical implementation, at the end of each time step $t_n$, we renormalize the tangent vector, so that $\T^{(n)}\in\Hbb^2$.  

% -------------------------------------------------------------------------------
%\subsubsection{Stability and computation cost}
% --------------------------------------------------------------------------------
To determine the stability constraints of the numerical scheme, we compute the maximum value of the time step $\Delta t$ for which the solution does not blow up. Thus, after giving different values to the parameters $N$, $M$, $l$, we obtain $\Delta t/\Delta s^2 = 0.5302\ldots$, i.e., $\Delta t = \mathcal{O}(\Delta s^2)$.

Let us remark that, in the case of regular polygons in the Euclidean space, the space derivatives are approximated at $N$ nodes by using the \texttt{fft} algorithm in MATLAB \cite{HozVega2014,HozKumarVega2019}. Due to the symmetries of the tangent vector, it was possible to do this by using only one side of the $M$-sided polygon, needing a computation cost of $\mathcal{O}((N/M)\log(N/M))$. However, in the current scenario with fixed boundary conditions, we work with all the sides of the truncated $l$-polygon, and the space derivatives are approximated with finite difference matrices of size $(N+1)\times(N+1)$, hence, making the problem challenging from a computational point of view as well. 

% --------------------------------------------------------------------------------
\subsection{Numerical results}
\label{sec:num-res}
% --------------------------------------------------------------------------------
Recall that, given any rational time, the computation of the algebraic solutions $\Xalg$ and $\Talg$ is entirely based on the assumption of uniqueness. In the following lines, we will see that, up to some numerical errors, the numerical solutions, denoted by $\Xnum$ and $\Tnum$, match very well the ones obtained from the theoretical arguments. However, remark that, in order to compare the two solutions, we need to specify the movement of $\Xalg$ at any rational times. This is done by computing the center of mass, which is given by the mean of $\X$, i.e.,
$$
\X^{mean}(t) = \frac{1}{L} \int_{-L/2}^{L/2}\X(s,t) ds.
$$
Thus, with the discretization mentioned above, we approximate the integral numerically by using the trapezoidal rule. Being our aim to analyze $\X^{mean}$ componentwise, we note, from the symmetries mentioned in Section \ref{sec:Spatial-sym}, that, for any given time $t$, the first component is equal to zero, while the second and third components, i.e., $X^{mean}_{2,0}$, and $X^{mean}_{3,0}$, are calculated as the mean of $N$ values of $X_2$ and $X_3$, respectively.

Here, $X_{3,0}^{mean}$ describes the position of the center of mass along the $z$-axis, i.e., the vertical height of the polygonal curve $\X$. After carrying out numerical simulations for different values of $M$ and $l$, it has been observed that $X_{3,0}^{mean}(t)$ can be very well approximated by means of a constant multiplied by $t$. More precisely,
\begin{equation}
\label{eq:ht-def}
X_{3,0}^{mean}(t) \approx \frac{X_{3,0}^{mean}(T_f)}{T_f} t = c_l^{num} \ t,
\end{equation}
where $c_l^{num}$ is the mean speed computed numerically, and its exact value is obtained as a consequence of Theorem \ref{thm:c_l_expression}.

\begin{figure}[htbp!] 
	\centering    
	\begin{subfigure}[b]{0.476\textwidth}
		\centering
		\includegraphics[width=\textwidth]{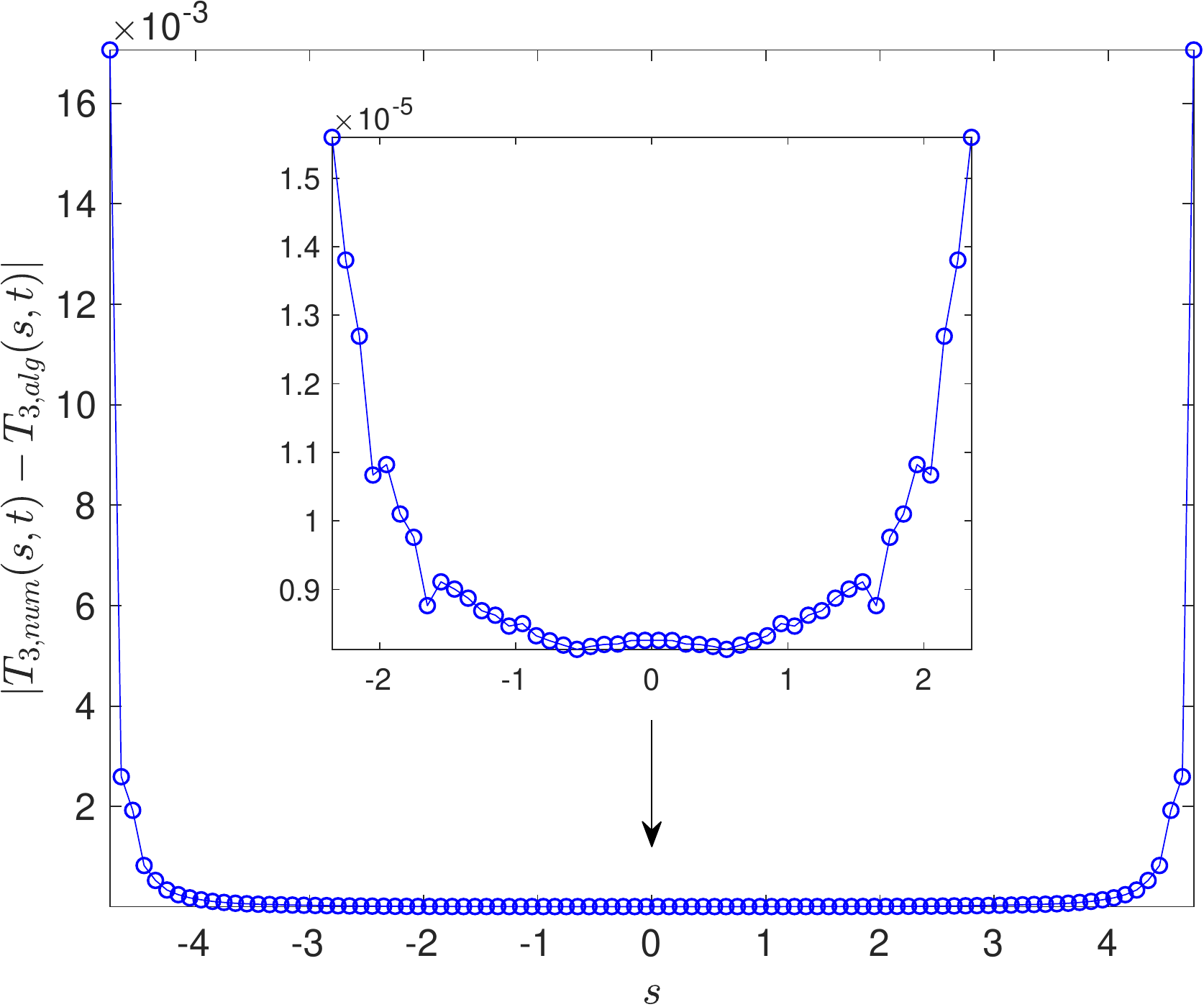}
		\caption{The error $|T_{3,num}(s,t)-T_{3,alg}(s,t)|$, at $t=T_f$, $s\in[-L/2,L/2]$, where the values represented by blue circles are computed using the mean of $N/M$ values for each side. The magnified part shows that the error is much smaller for the inner part than near the boundary. }
		\label{subfig:T3num-errorA}
	\end{subfigure}
	\hfill
	\begin{subfigure}[b]{0.484\textwidth}
		\centering
		\includegraphics[width=\textwidth]{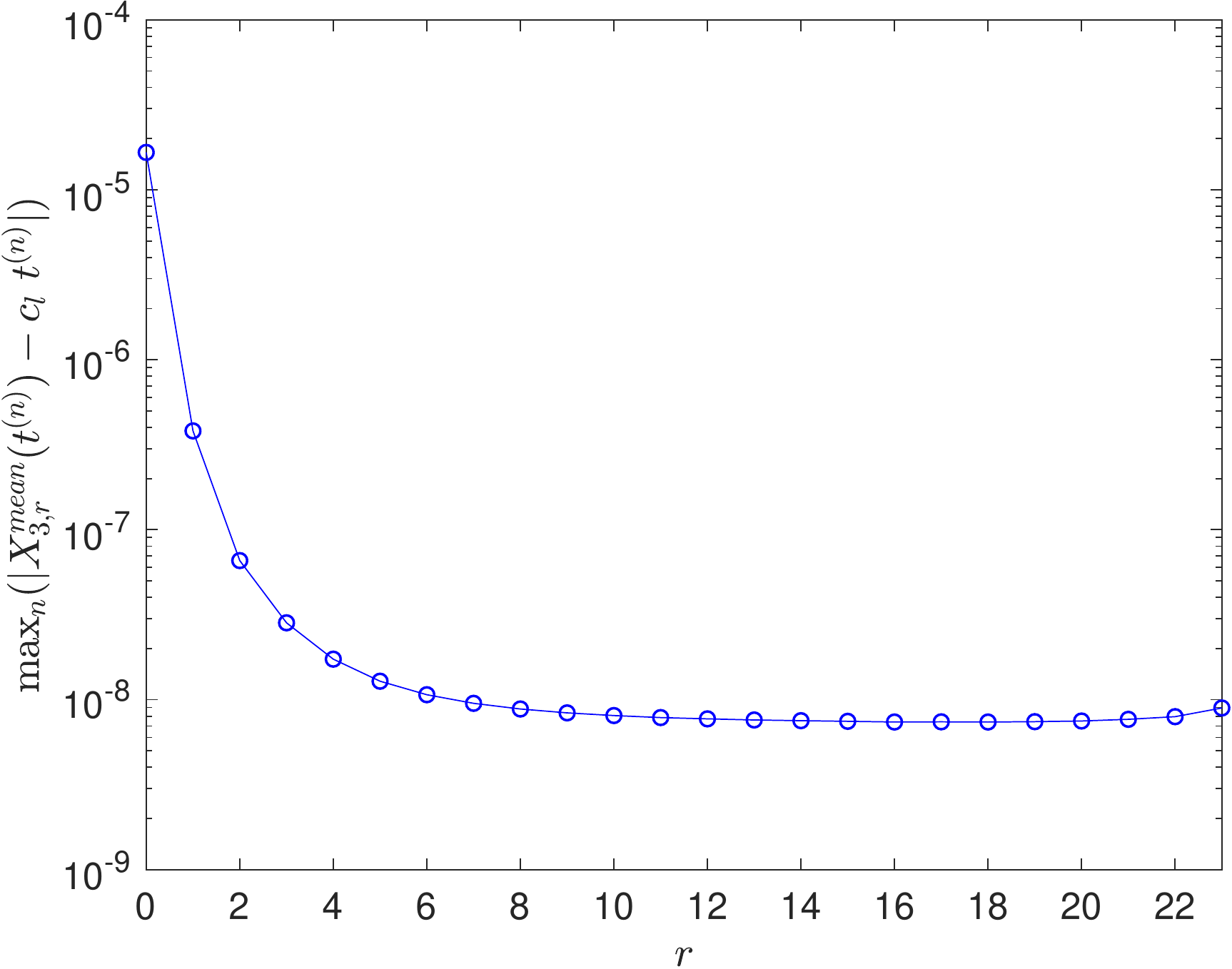}
		\caption{The error $\max_n(|X_{3,r}^{mean}(t^{(n)})-c_l t^{(n)}|)$, where the mean of $X_3(s_j,t^{(n)})$ has been computed using $4rN/M$ values, such that $s_j \in [-L/2+2rl,L/2-2rl]$, $r=0,1,2,\ldots, M/4-1$, $c_l =  1.000416458444891$. Clearly, the error reduces significantly when only the inner points are included in the computations.}
		\label{subfig:T3num-errorB}
	\end{subfigure}		
	\caption{The errors for $M=96$, $l=0.1$, $N/M=2^9$.}
	\label{fig:T3num-error}
\end{figure}

On the other hand, the values of $\Tnum$ (hence, those of $\Xnum$) corresponding to the inner grid points are found to be far more accurate than the ones close to the boundary. This is due to the exponential growth of the tangent vector and to the fact that we are approximating piecewise continuous functions using a finite difference scheme. For instance, Subfigure \ref{subfig:T3num-errorA}, shows the error $|T_{3,num}(s,t)-T_{3,alg}(s,t)|$, considering the third component of the algebraic and numerical solutions of $\T$, for $M=96$, $l=0.1$, $N/M=2^9$, $s\in[-L/2,L/2]$, at $t=T_f$. The error is of $\mathcal{O}(10^{-5})$ in the magnified part, whereas it is of $\mathcal{O}(10^{-2})$ near the boundary. This indicates that the inner part of the polygon is more accurate than the one close to the endpoints. Therefore, for a good approximation of $\X^{mean}$, we choose to work with the inner points of the discretized domain, where, in order to make a reasonable choice of the ``inner points'', we define 
\begin{equation}
\label{eq:htr-def}
X_{3,r}^{mean}(t) = \frac{1}{N_r} \sum_{j=2rN/M}^{N-2rN/M-1} X_3(s_j,t), \ r=0,1,\ldots,M/4-1,
\end{equation}
for $N_r = N-(4rN/M)$, i.e., the mean of $X_3(s_j,t)$, for $s_j \in [-L/2+2rl , L/2- 2rl ]$. Then, for each $r$, we compute the error $\max_n(|X_{3,r}^{mean}(t^{(n)})-c_l t^{(n)}|)$, i.e., the maximum difference between $X_{3,r}^{mean}(t)$ and its exact linear approximation $c_l t$. Figure \ref{subfig:T3num-errorB} shows that the error is smaller when the nodes closer to the boundary are avoided. It also shows that, after a certain value of $r$, the error does not vary much; consequently, without loss of generality, we choose $r=M/8$, i.e., $s_j \in [-L/4,L/4]$, $j=N/4+1,N/4+2,\ldots,3N/4$. Note that, although using the symmetries, $X_{3,0}^{mean}(t^{(n)})$ can be computed by using only $N/M$ values, we prefer to work with $N/2$ elements, due to the unevenness of errors discussed above.

In order to further strengthen our claim to \eqref{eq:ht-def}, we compute the error $\max_n(|X_{3,r}^{mean}(t^{(n)})-c_lt^{(n)}|)$, for different values of $l$ and a fixed $r$. Since a regular planar $l$-polygon is characterized by the parameter $l$, the speed of the center of mass depends only on it. However, we are approximating the infinitely long $l$-polygon with the parameter $M$ and, as a consequence, better results are obtained for larger values of $M$. In our simulations, we work with  moderately large values of $M$, and different values of $l$ and $N/M$. Table \ref{table:fixM96-diffl-NM} displays the corresponding errors, and it is evident that, whenever the number of grid points is doubled, the error decreases by a factor slightly lower than two, hence, suggesting a convergence of the order of $\mathcal{O}((N/M)^{-1})$. For small values of $l$, $c_l$ is very close to 1, so we provide the value of $c_l - 1$, for each $l$. Note that $c_l$ converges to $1$, as $l$ goes to zero, i.e., $\X(s,0)$ tends to a hyperbola. 

\begin{table}[h]
	\centering
	\begin{tabular}{|l| l| l | l| l| l|}
		\hline 
		$l$ &  $N/M=2^{6}$& $N/M=2^{7}$ & $N/M=2^{8}$ &  $N/M=2^9$ &$(c_l-1)$\\
		\hline 
		$0.15$ & $2.0578\cdot10^{-7}$  &$1.1334\cdot10^{-7}$&$6.3398\cdot10^{-8}$ & $3.7051\cdot10^{-8}$ & $9.3645\cdot10^{-4}$\\
		$0.12$ & $ 8.4311\cdot10^{-8}$ &$4.6133\cdot10^{-8}$&$2.6124\cdot10^{-8}$ & $1.5543\cdot10^{-8}$ & $5.9957\cdot10^{-4}$ \\
		$0.1$  & $4.0669\cdot10^{-8}$ &$ 2.2237\cdot10^{-8}$&$1.2801\cdot 10^{-8}$ & $7.7110\cdot 10^{-9}$ & $4.1646\cdot10^{-4}$ \\
		$0.05$ &$2.5460\cdot10^{-9}$ &$ 1.5103\cdot10^{-9}$& $9.2425\cdot10^{-10}$ &$6.1407\cdot10^{-10}$& $1.0415\cdot10^{-4}$ \\
		$0.025$ & $1.6008\cdot10^{-10}$&$1.2481\cdot10^{-10}$&$8.8786\cdot10^{-11}$ & $7.1406\cdot10^{-11}$&$2.6040 \cdot10^{-5}$ \\
		\hline
	\end{tabular}
	\caption{The error $\max_n(|X_{3,r}^{mean}(t^{(n)})-c_lt^{(n)}|)$, for $M=96$ and different $N/M$, $l$, where $X_{3,r}^{mean}(t^{(n)})$ is computed using \eqref{eq:htr-def}, for $r=M/8$. After doubling $N/M$, the error reduces by a factor close to two, showing first-order convergence.}
	\label{table:fixM96-diffl-NM}	
\end{table}
\begin{figure}[h ]
	\centering
	\begin{subfigure}[b]{0.24\textwidth}
		\centering
		\includegraphics[width=\textwidth, clip=true, valign =t]{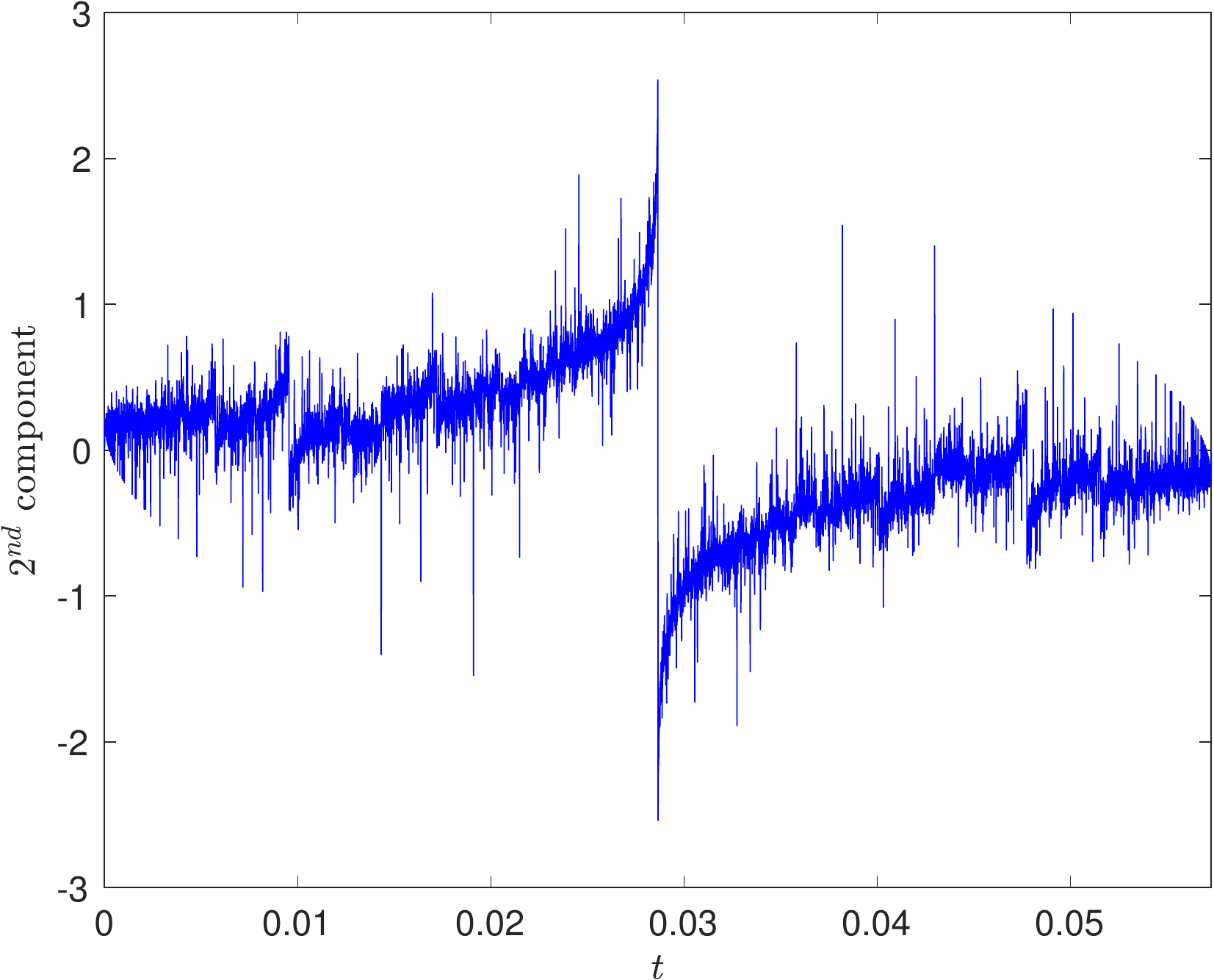}
		\caption{}
		\label{subfig:X23mean-algA}   
	\end{subfigure}
	\hfill
	\begin{subfigure}[b]{0.24\textwidth}	
		\includegraphics[width=\textwidth, clip=true, valign =t]{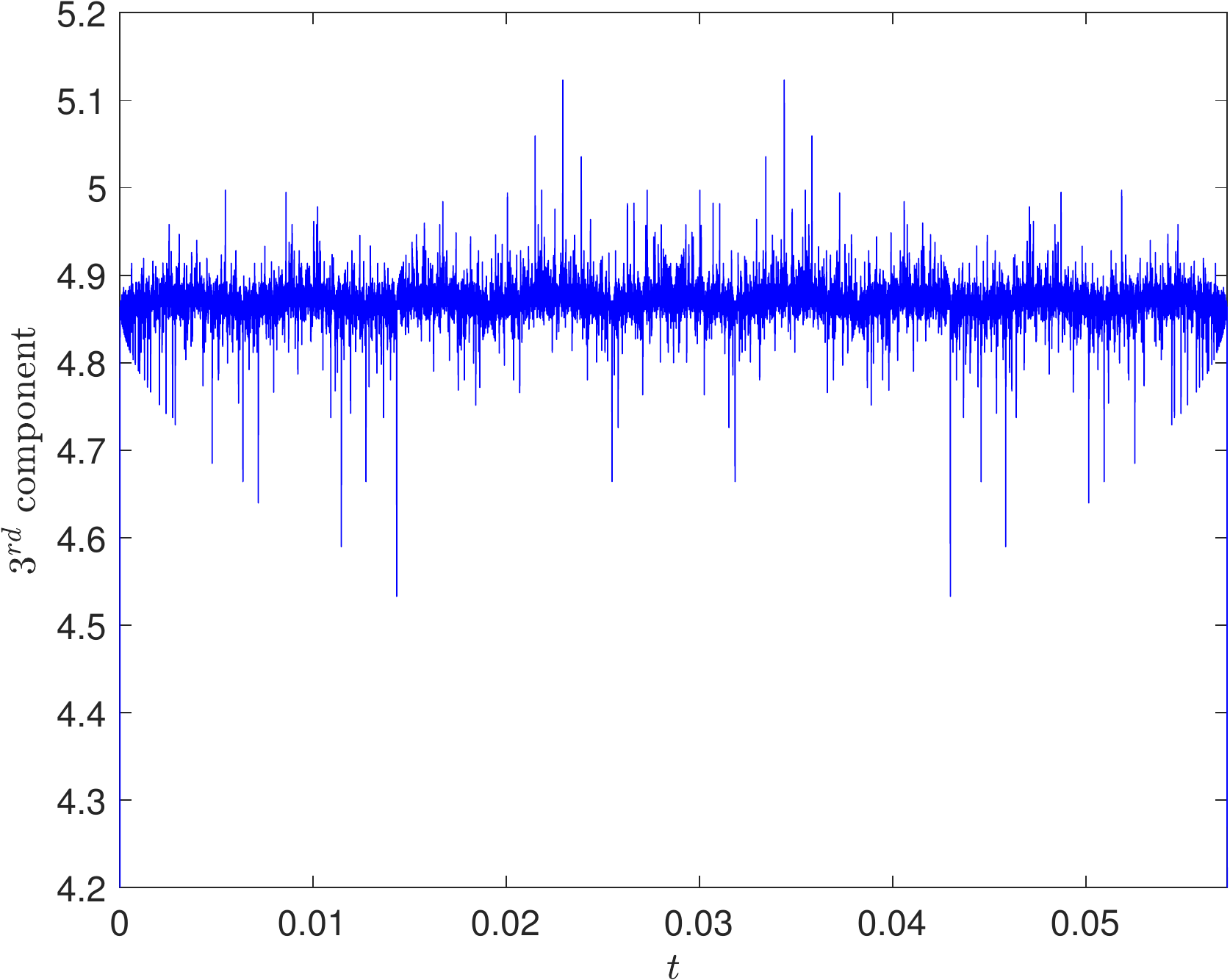}
		\caption{}
		\label{subfig:X23mean-algB}   	
	\end{subfigure}
	\hfill
	\begin{subfigure}[b]{0.24\textwidth}
		\includegraphics[width=\textwidth, clip=true, valign =t]{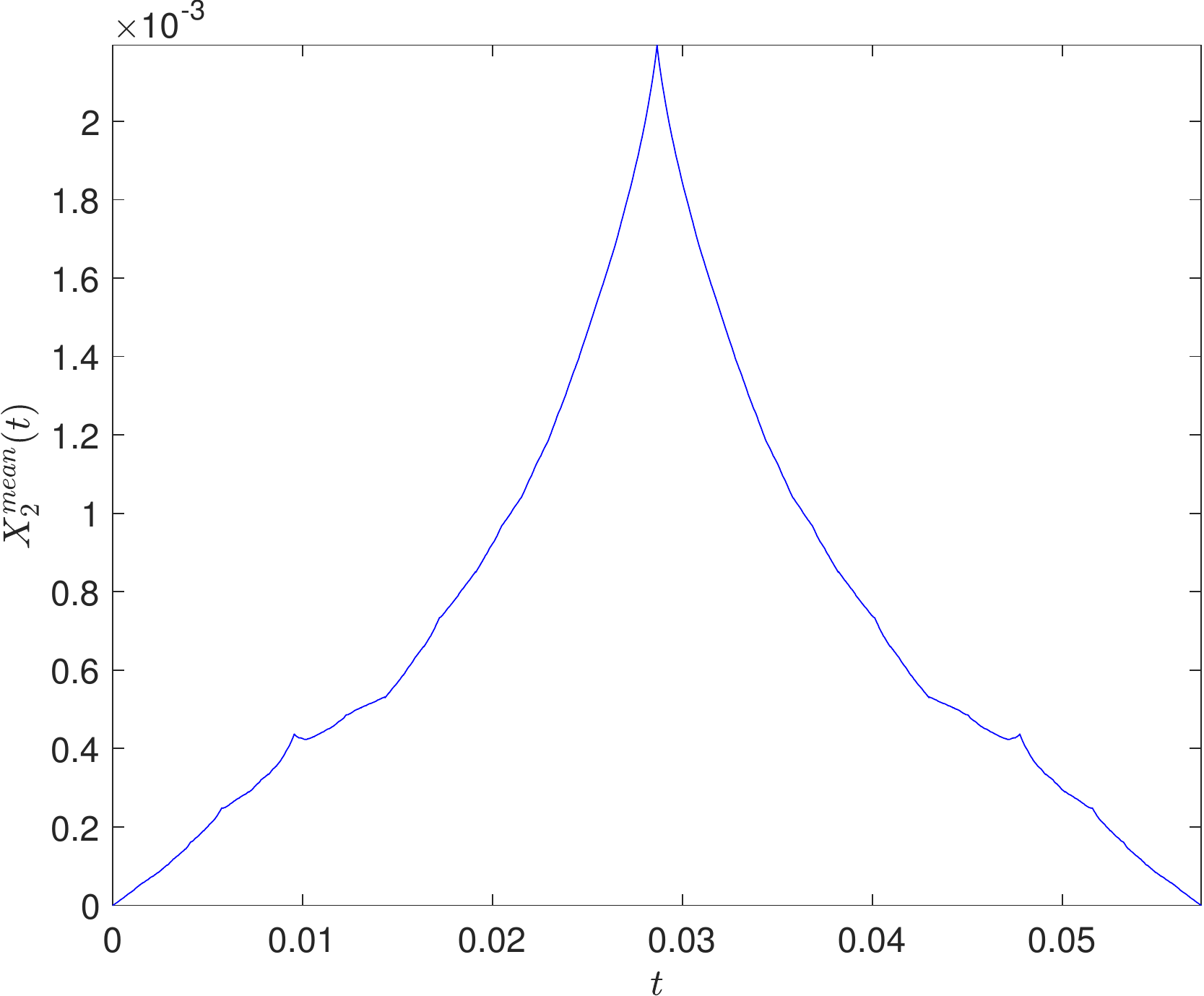}
		\caption{}
		\label{subfig:X23mean-algC}   	
	\end{subfigure}
	\hfill
	\begin{subfigure}[b]{0.24\textwidth}
		\includegraphics[width=\textwidth, clip=true, valign =t]{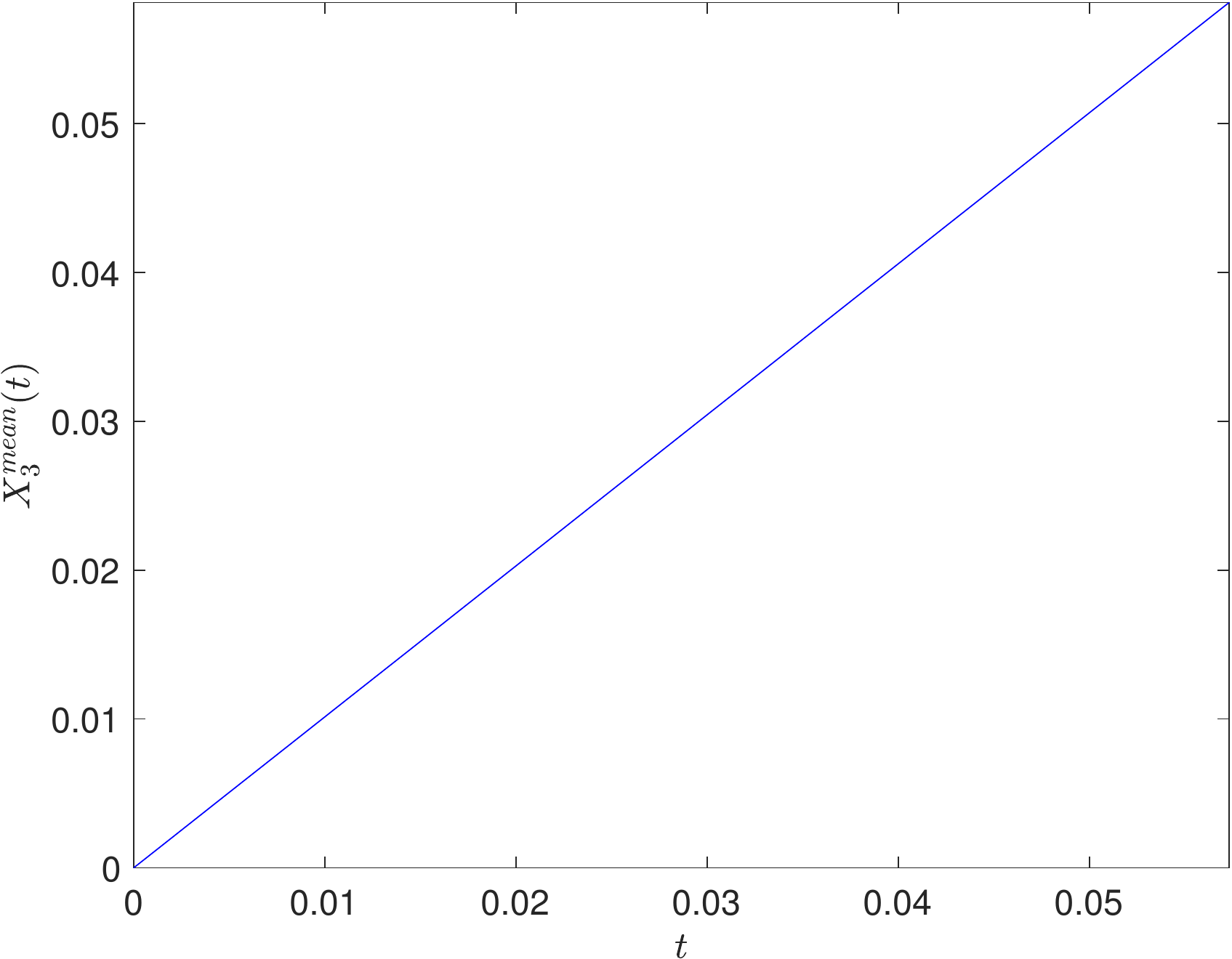}
		\caption{}	
		\label{subfig:X23mean-algD}   	
	\end{subfigure}
	\hfill	
	\caption{Second and third components of \eqref{eq:ht-CoM} and \eqref{eq:ht-CoM2}, computed for $M=8$, $l=0.6$, $q=7560$. After integrating with respect to time, the oscillations completely disappear, and we obtain a periodic curve and a straight line, respectively.} 
	\label{fig:X23mean-alg}   
\end{figure}
Furthermore, it is possible to approximate $\X^{mean}$ also from the algebraic solution, hence minimizing the numerical errors. Using the approach in \cite[Section 4]{HozVega2018}, we compute it as
\begin{align}\label{eq:ht-CoM}
\X^{mean}(t) = \int_{0}^{t} \mean (\Xt)(t^\prime) dt^\prime = \int_{0}^{t} \left[ \frac{1}{L} \int_{-L/2}^{L/2} \Xt(s,t^\prime) ds \right] dt^\prime.
\end{align}
For any rational time $\tpq$, the first integral is given precisely by
\begin{equation}
\label{eq:ht-CoM2}
\int_{-L/2}^{L/2} \Xt (s,\tpq) ds = \frac{l_q}{\sinh (l_q)} \sum_{k=0}^{Mq-1} \T_{alg,k} \wedge_- \T_{alg,k+1}, 	
\end{equation}
where $\T_{alg,k} = \Talg(s_k^+,\cdot)$, and, taking a large $q$, the integral with respect to time in \eqref{eq:ht-CoM} can be approximated with third-order accuracy. Taking $M=8$, $l=0.6$, $q=7560$, the interval $[0,T_f]$ has been divided into $q$ equally spaced segments, and we have plotted the integral in \eqref{eq:ht-CoM2}, whose first component is zero, and the other two seem to have a very oscillatory behavior, as shown in Subfigures \ref{subfig:X23mean-algA} and \ref{subfig:X23mean-algB}. However, after integrating in time, the oscillations disappear, and we obtain the components of $\X^{mean}$, where the second component is periodic and the third component is a straight line whose slope converges to $c_l$ with $q$, as shown in Subfigures \ref{subfig:X23mean-algC}, \ref{subfig:X23mean-algD}.  

\subsubsection{Comparison between the numerical and algebraic solutions}

\begin{figure}[htbp!]
	\centering
	\begin{subfigure}[t]{0.468\textwidth}
		\centering
		\includegraphics[width=\textwidth, clip=true, valign=t]{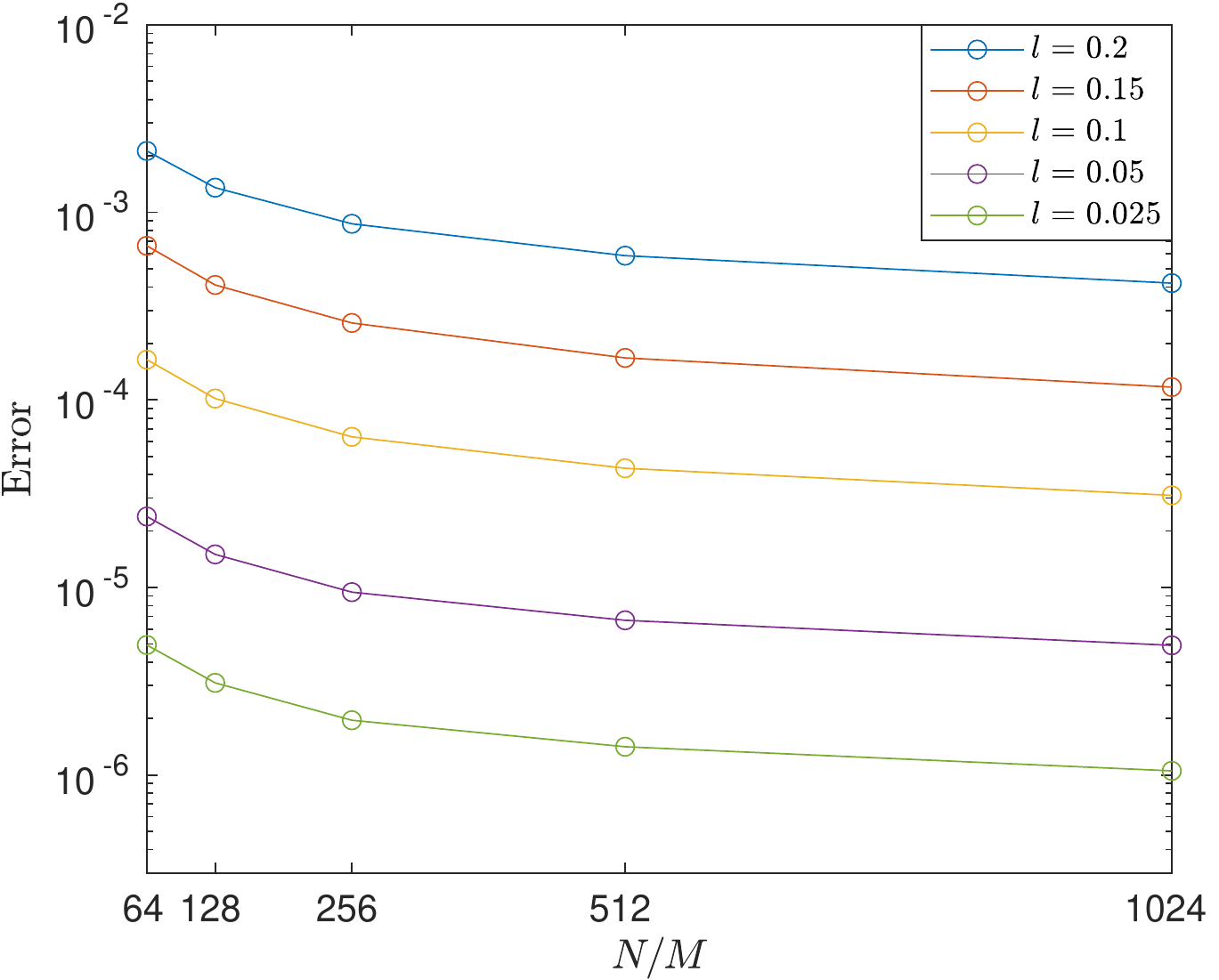}
		\caption{The error $\max_{j,n}(\| \Xnum(s_j,t^{(n)})-(0,X_{2,r}^{mean}(t^{(n)}),c_l \ t^{(n)}) - \Xalg(s_j,t^{(n)})\|)$, for different values of $N/M$, $l$, $M=48$. The plot shows that, even if slowly, convergence indeed occurs as $N/M$ increases.}
		\label{subfig:XnumA}
	\end{subfigure}
	\hfill
	\begin{subfigure}[t]{0.492\textwidth}
		\centering
		\includegraphics[width=\textwidth, clip=true, valign=t]{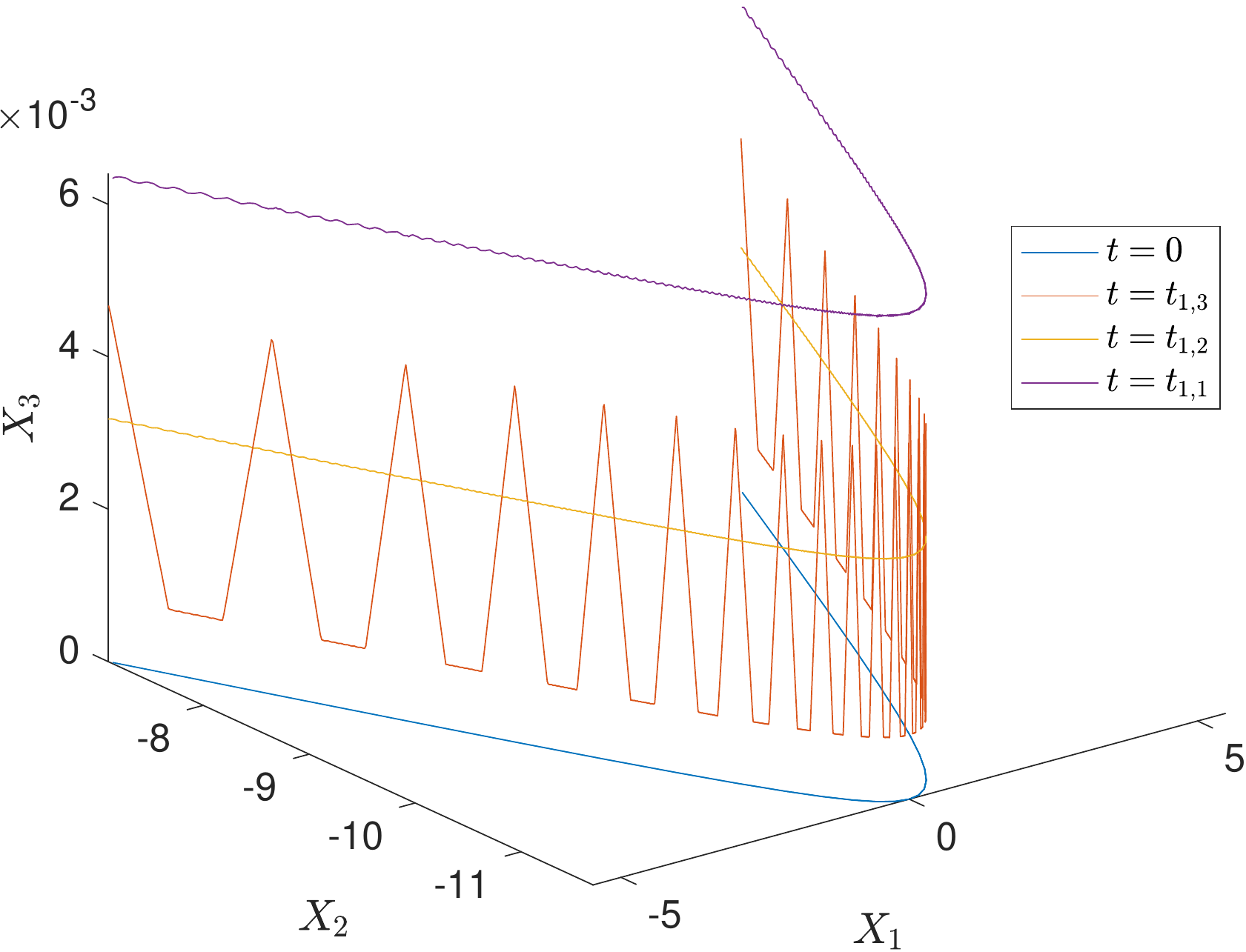}
		\caption{$\Xnum(s,t)$, for $l=0.2$, $M=48$, $N/M=2^{11}$. Besides the constant vertical movement, at half the time period and at the end of it, the planar polygon reappears, and at one-third of the time period, three times as many sides are formed in the non-planar polygon.}
		\label{subfig:XnumB}
	\end{subfigure}	
	\caption{}
	\label{fig:Xnum}   
\end{figure}

For given $M$ and $l$, we subtract the movement of the center of mass from $\Xnum$, and compare it with $\Xalg$. Recall that the algebraic solution $\Xalg$ corresponds to the vertices of the polygonal curve, and the non-vertex values can be computed using linear interpolation. We calculate the error $\max_{j,n}(\| \Xnum(s_j,t^{(n)})-(0,X_{2,r}^{mean}(t^{(n)}),c_l \ t^{(n)}) - \Xalg(s_j,t^{(n)})\|)$, where $\| \cdot \|$ is the Euclidean norm. On the other hand, given the size of the discretization, it is computationally very difficult to compare the solutions at all the $N_t+1$ time instants; therefore, we do it for a fairly large amount, e.g., $N_t=1260$. Continuing as previously, in Table \ref{table:fixM48-diffl-N2}, we show the errors for $M=48$, $r=M/8$ and different values of $l$, $N/M$; their plots in logarithmic scale in Subfigure \ref{subfig:XnumA} confirm that the errors decrease by a factor close to 1.6, when halving the space step size. Although the convergence is slow, bearing in mind that $\max\| \Xalg\|\gg1$, the results are satisfactory and show that, as $N$ grows larger, the numerical solution converges to the algebraic one. It also gives strong evidence that, up to the vertical height, the evolution of $\X$ is $T_f$-periodic in time. Subfigure \ref{subfig:XnumB} shows the inner $N/2$ points of $\Xnum$ for $M=48$, $l=0.2$, $N/M=2^{11}$, and it can be clearly observed that, besides the planar curve at the initial, middle and final times of the time period, three times as many sides appears at one-third of the time period.
\begin{table}
	\centering
	\begin{tabular}{c c c c c c  }
		\hline
		$l$ & $N/M=2^{6}$& $N/M=2^{7}$ & $N/M=2^{8}$ &  $N/M=2^{9}$ &  $N/M=2^{10}$ \\
		\hline
		$0.2$ & $2.1238\cdot10^{-3}$& $1.3533\cdot10^{-3}$ &$8.6836\cdot10^{-4}$ & $5.8715\cdot10^{-4}$ & $4.1921 \cdot10^{-4}$ \\		
		$0.15$ &$6.6274\cdot10^{-4}$ &$4.0982\cdot10^{-4}$& $2.5702\cdot10^{-4}$& $1.6740\cdot10^{-4}$ & $1.1695\cdot10^{-4}$ \\		
%		$0.12$ &$1.5620\cdot10^{-2}$ &$1.2359\cdot10^{-2}$ &$9.8698\cdot10^{-3}$ & $8.4188\cdot10^{-3}$ & $7.2017\cdot10^{-3}$\\
		$0.1$ & $1.6388\cdot10^{-4}$ & $ 1.0164\cdot10^{-4}$ & $6.3520\cdot10^{-5}$& $4.3188\cdot10^{-5}$ & $3.0973\cdot 10^{-5}$\\		
		$0.05$ &$2.3886\cdot10^{-5}$ & $1.5006 \cdot 10^{-5}$& $9.4398\cdot 
		10^{-6}$& $6.6809\cdot10^{-6}$ & $4.9127\cdot10^{-6}$ \\
		$0.025$ &$4.9348\cdot10^{-6}$ &$3.0992\cdot10^{-6}$ &$1.9590 \cdot10^{-6}$&  $1.4149\cdot10^{-6}$& $1.0525 \cdot 10^{-6}$\\
		\hline
	\end{tabular}
	\caption{The error $\max_{j,n}(\| \Xnum(s_j,t^{(n)})-(0,X_{2,r}^{mean}(t^{(n)}),c_l t^{(n)})- \Xalg(s_j,t^{(n)})\|)$, for $r=M/8$, $j=N/4+1,\ldots,3N/4+1$, $n=0,1,\ldots,1260$, $M=48$. }
	\label{table:fixM48-diffl-N2}
\end{table}

%-------------------------------------------------------
\begin{figure}[htbp!] \centering
	\begin{subfigure}[t]{0.48\textwidth}
		\centering
		\includegraphics[width=\textwidth, clip=true]{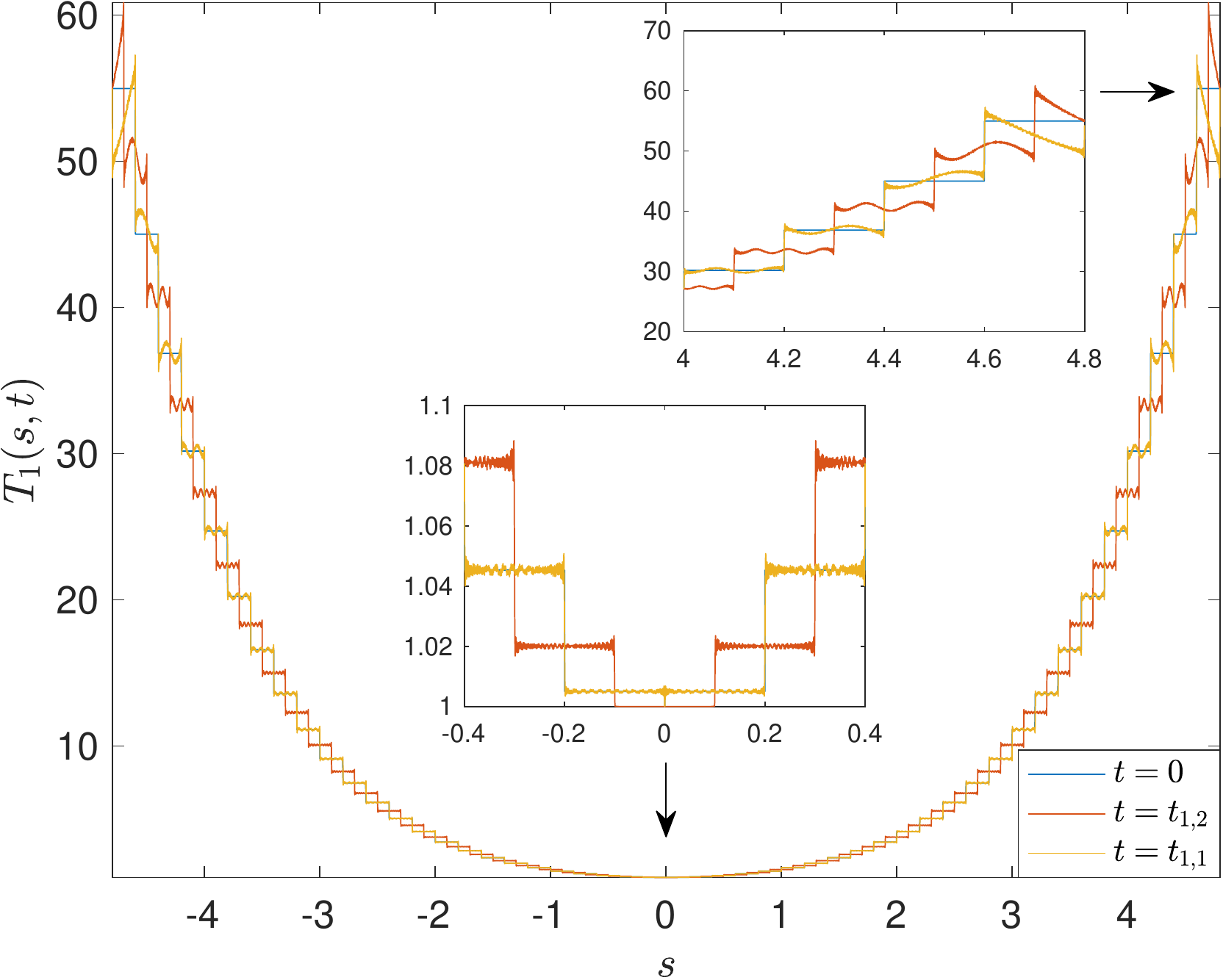}
		\caption{For $M=48$, $l=0.2$, $N/M=2^{11}$, the magnified part around $s=0$ confirms that, at $t=t_{1,2}$, the tangent vector is continuous at $s=0$, so there is no vertex located at $s=0$, as indicated in \eqref{eq:psi_tpqfinal} for $q=2$. At $t=t_{1,1}$, it matches the one at $t=0$, thus, showing time periodicity. The oscillations are more prominent near the boundary, shown in the magnified part on the right-hand side.}
		\label{subfig:T1num-M48l02A}				
	\end{subfigure}
	\hfill
	\begin{subfigure}[t]{0.48\textwidth}
		\centering
		\includegraphics[width=\textwidth, clip=true]{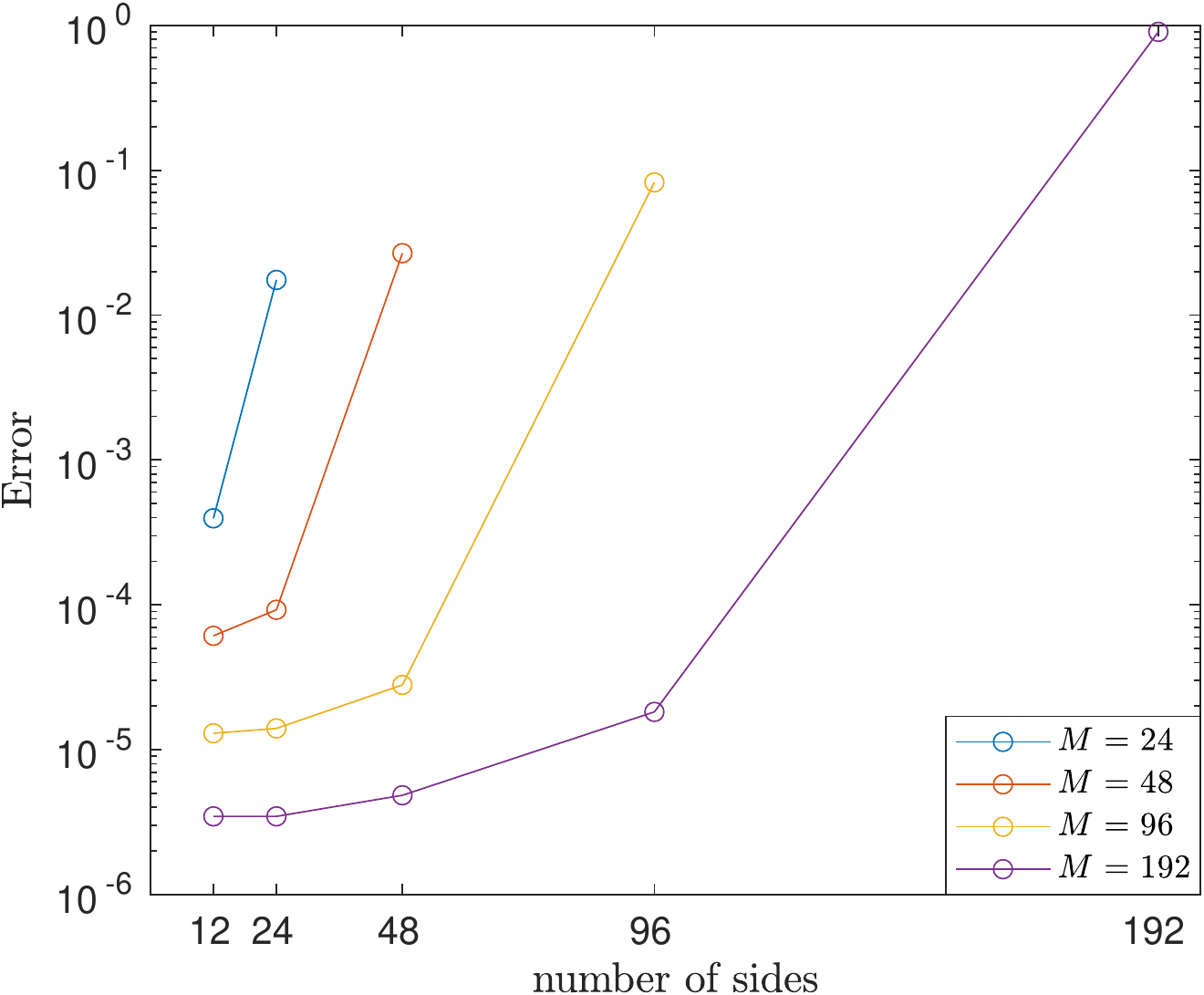}
		\caption{The error $\max_j\|\Talg(s_j,T_f)-\Tnum(s_j,T_f)\|$, for $j=M/2-M/2^n+1,\ldots,M/2+M/2^n$, $n=1,2,\ldots,\log_2(M/6)$, $l=0.025$. The error for an amount of sides equal to $24$ is much smaller if they are part of an $l$-polygon with $M = 48$ sides (in red), than if they correspond to an $l$-polygon with $M = 24$ sides (in blue), and so on. Thus, as $M$ increases, so does the amount of sides for which very accurate results are achieved.}
		\label{subfig:T1num-M48l02B}
	\end{subfigure}		
	\caption{}
	\label{fig:T1num-M48l02}   
\end{figure}
	
	On the other hand, Subfigure \ref{subfig:T1num-M48l02A} shows $T_{1,num}$, the first component of the tangent vector, at different rational times. From the magnified part, it is clear that, at half the time period, the tangent vector is continuous at $s=0$, i.e., there is no corner at that time, which is also consistent with \eqref{eq:psi_tpqfinal}, for $q=2$; moreover, the oscillations causing the errors are more prominent toward the boundary. At the end of one time-period, up to the numerical errors, the solution matches the one at the initial time, thus showing the time periodicity of $\Tnum$ (in yellow). 
	
It is worth emphasizing that, despite the challenges posed by the nature of the problem, e.g., lack of regularity in the initial data, exponential growth of the tangent vector near the boundary, etc., we are able to capture the evolution of the central part of the polygonal curve very accurately with the numerical simulations. Moreover, for a given $l>0$, this can be further improved by increasing $M$. For instance, in Subfigure \ref{subfig:T1num-M48l02B}, taking $l=0.025$, we have considered the corresponding $l$-polygons with $M$ sides, for $M = 24, 48, 96, 192$, and have compared the corresponding $\Talg$ and $\Tnum$ at $t=T_f$, for different amounts of inner sides. More precisely, we have calculated the error $\max_j\|\Talg(s_j,T_f)-\Tnum(s_j,T_f)\|$, for $j=M/2-M/2^n+1,\ldots,M/2+M/2^n$, $n=1,2,\ldots,\log_2(M/6)$, where $\Tnum$ is computed by taking the mean of the central part of each side. The plot shows that, if we fix the amount of sides on the $x$-axis, then the error decreases when these sides are considered as being part of an $l$-polygon with a larger number of sides. For example, when the amount of sides is equal to 24, the error is much smaller if they are part of an $l$-polygon with $M=48$ sides (in red), than if they correspond to an $l$-polygon with $M=24$ sides (in blue). This process can be continued further, and we observe the same behavior; hence, there is strong numerical evidence that the convergence indeed occurs when $M$ tends to infinity.

\subsubsection{Trajectory $\X(0,t)$}
\label{sec:X0t}
\begin{figure}[htbp!] 
	\centering
	\begin{subfigure}[t]{0.433\textwidth}
		\centering
		\includegraphics[width=\textwidth, clip=true]{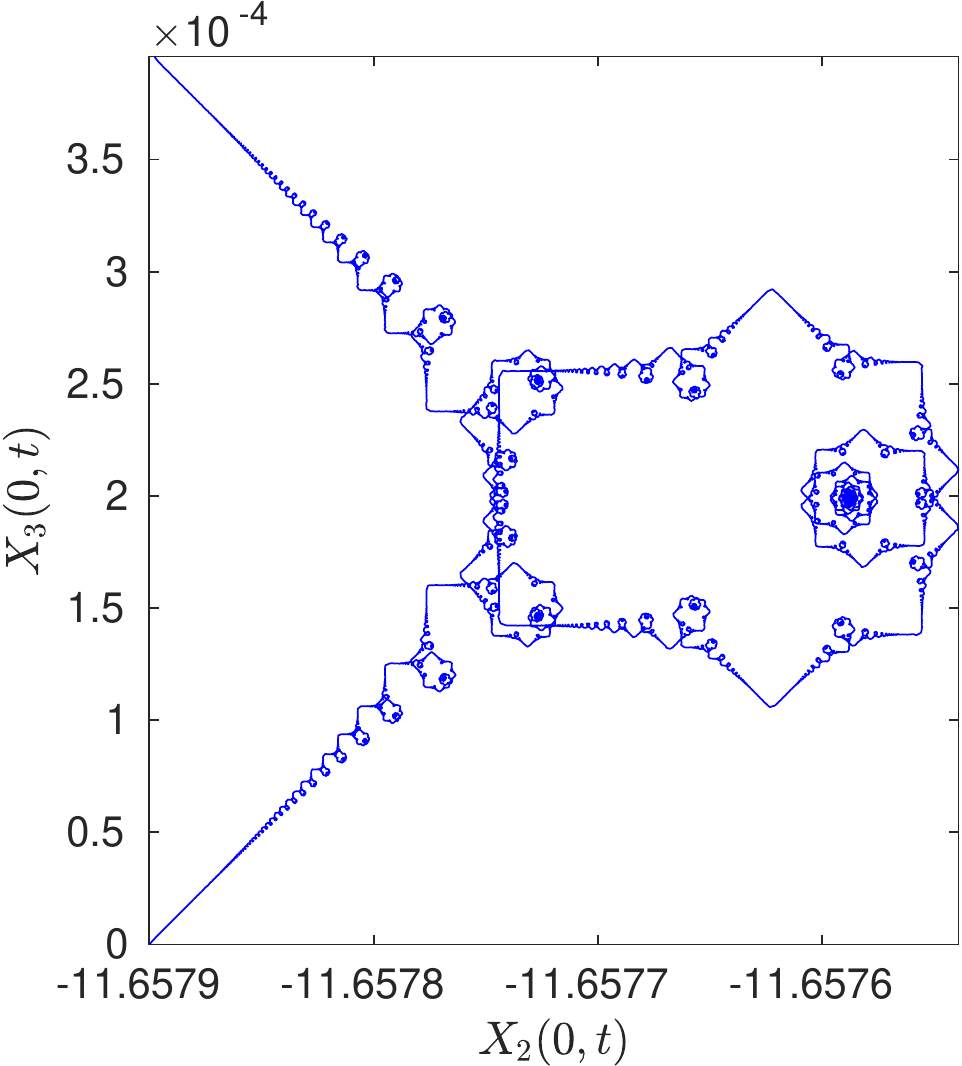}
		\caption{$z(t)$ as in \eqref{eq:zt}, for $M=192$, $l=0.05$.}
		\label{subfig:X0tM192l0025A}				
	\end{subfigure}
	\hfill	
	\begin{subfigure}[t]{0.527\textwidth}
		\centering
		\includegraphics[width=\textwidth, clip=true]{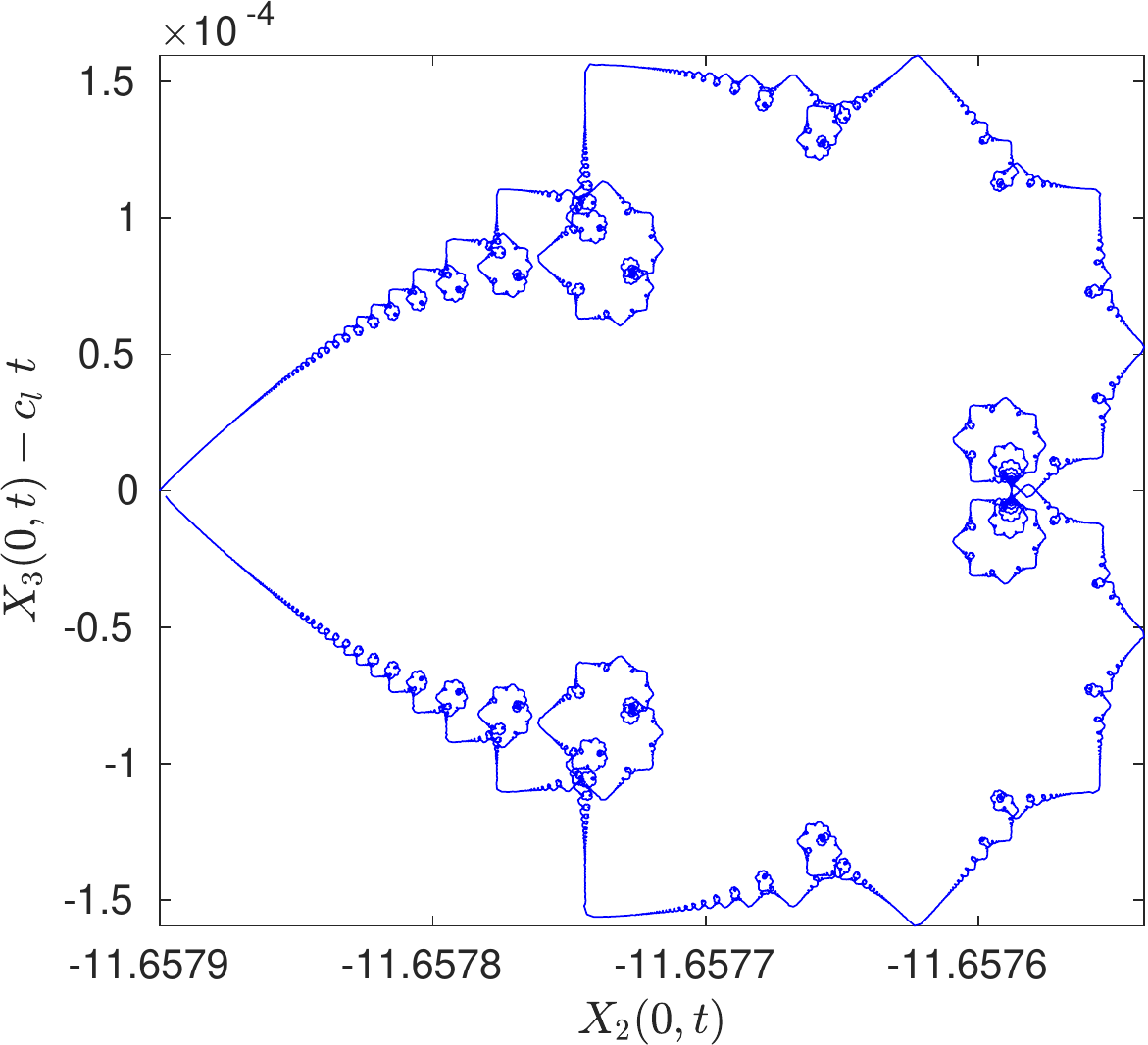}
		\caption{$z_l(t)$ as in \eqref{eq:z_lt}, for $M=192$, $l=0.05$, $t\in[0,T_f]$.}
		\label{subfig:X0tM192l0025B}				
	\end{subfigure}
	\caption{}
	\label{fig:X0tM192l0025}
\end{figure}

The choice of the initial data (i.e., an even number of sides for $\X$) allows us to capture the time evolution of $\X(0,t)$. Due to the mirror symmetries of $\X$ given in Section \ref{sec:Spatial-sym}, during the time evolution, the $z$-axis and $\X(-L/2+kl/2,t)$, for $k=0,1,\dots,2M$, always lie in the same plane, for all $t\geq0$, where an even value of $k$ corresponds to the vertices, and an odd value, to the middle point of the sides. For instance, the numerical simulations show that $\X(0,t)$ lies in the YZ-plane, whereas $\X(-L/2+lk,t)$, for $k=0,1,2,\ldots,M$, lies in the plane obtained after rotating counterclockwise the YZ-plane  by a time-like angle $L/2-lk$ about the space-like $z$-axis. Thus, without loss of generality, we choose to observe the trajectory of $\X(0,t)$, and, after projecting it onto $\mathbb{C}$, we define
\begin{equation}
\label{eq:zt}
z(t) = X_2(0,t) + i X_3(0,t).
\end{equation}
As $\X$ is periodic in time up to a constant vertical movement, we introduce, for a given $l$,
\begin{equation}
\label{eq:z_lt}
z_l(t) = z(t) - i c_l \ t, \ t\in[0,T_f],
\end{equation}
which is $T_f$-periodic. Subfigures \ref{subfig:X0tM192l0025A} and \ref{subfig:X0tM192l0025B} show respectively $z(t)$ and $z_l(t)$, for $M=192$, $l=0.05$, $N/M=2^{11}$, which remind us of the multifractal structures obtained in the case of regular planar polygons in the Euclidean space \cite{HozVega2014}. This motivates us to compare $z_l(t)$ with the graph of 
\begin{equation}
\label{eq:phi-org-def}
\phi(t) = \sum_{k=1}^{\infty}\frac{e^{i \pi k^2 t}}{i\pi k^2}, \ t \in [0,2].
\end{equation}
Let us mention that $\phi(t)$ appeared in \cite{Du}, where its real part $f(t) = \sum_{k=1}^{\infty} \frac{\sin(\pi k^2t)}{\pi k^2}$, also called Riemann's non-differentiable function, was considered. Its geometrical and regularity properties have been studied recently in \cite{E2020}.

Recall that, in the numerical simulations, for large values of $M$, the value of $l$ needs to be chosen very small; hence, we can have $z_l(t)$ only for certain values of $l$. However, the computation of $\X(0, t)$ through the algebraic solution does not depend on $M$, and it is free from numerical errors. As a result, we can work with any value of $l$ (bearing in mind that, due to the exponential growth of $\T$, $l$ cannot be very large), and compute $z(t)$ algebraically. Bearing this in mind, let us define
\begin{equation}
\label{eq:z_lalg}
\begin{cases}
z_{l,alg}(t) = - (X_{2,alg}(0,t)+X_{2,alg}^{mean}(t)) + i X_{3,alg}(0,t),
\\[0.5em]
z_{alg}(t) = z_{l,alg}(t) + i c_l t, \ t\in[0,T_f],
\end{cases}
\end{equation}
where $X_{2,alg}^{mean}(t)$ is the second component of \eqref{eq:ht-CoM}. On the other hand, we will work with  
\begin{equation}
\label{eq:phi}
\phi(t) = -\sum_{k=1}^{\infty} \frac{e^{2\pi i k^2 t}}{k^2}, \ t\in[0,1],
\end{equation}
rather than with \eqref{eq:phi-org-def}.
\begin{figure}[htbp!] \centering
	\centering
	\begin{subfigure}[b]{0.301\textwidth}
		\centering
		\includegraphics[width=\textwidth, clip=true]{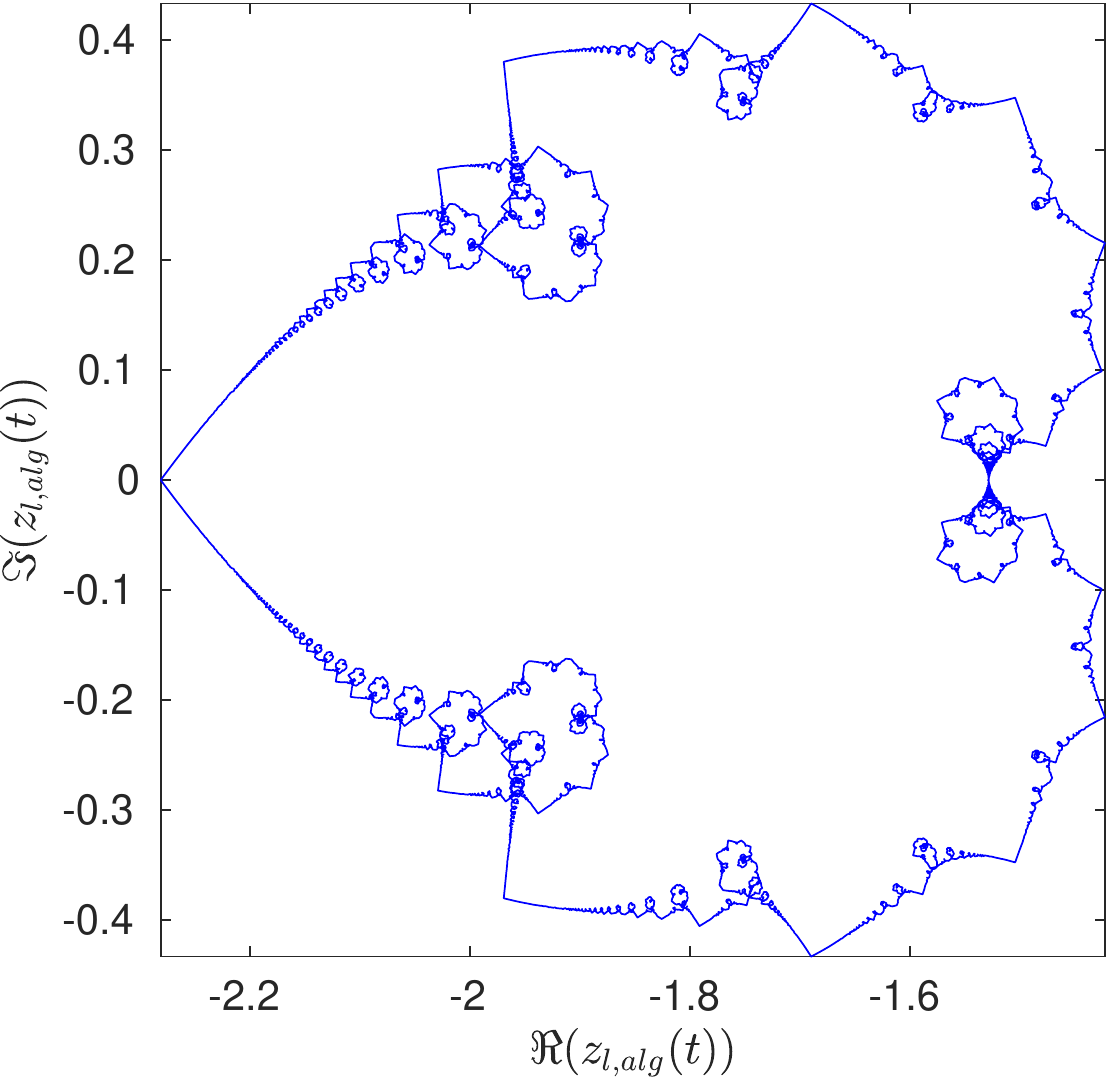}
		\caption{$l=2.6339\ldots$, i.e., $M_e=3$.}
		\label{subfig:zMtA}   
	\end{subfigure}
	\hfill
	\begin{subfigure}[b]{0.326\textwidth}
		\centering
	\includegraphics[width=\textwidth, clip=true]{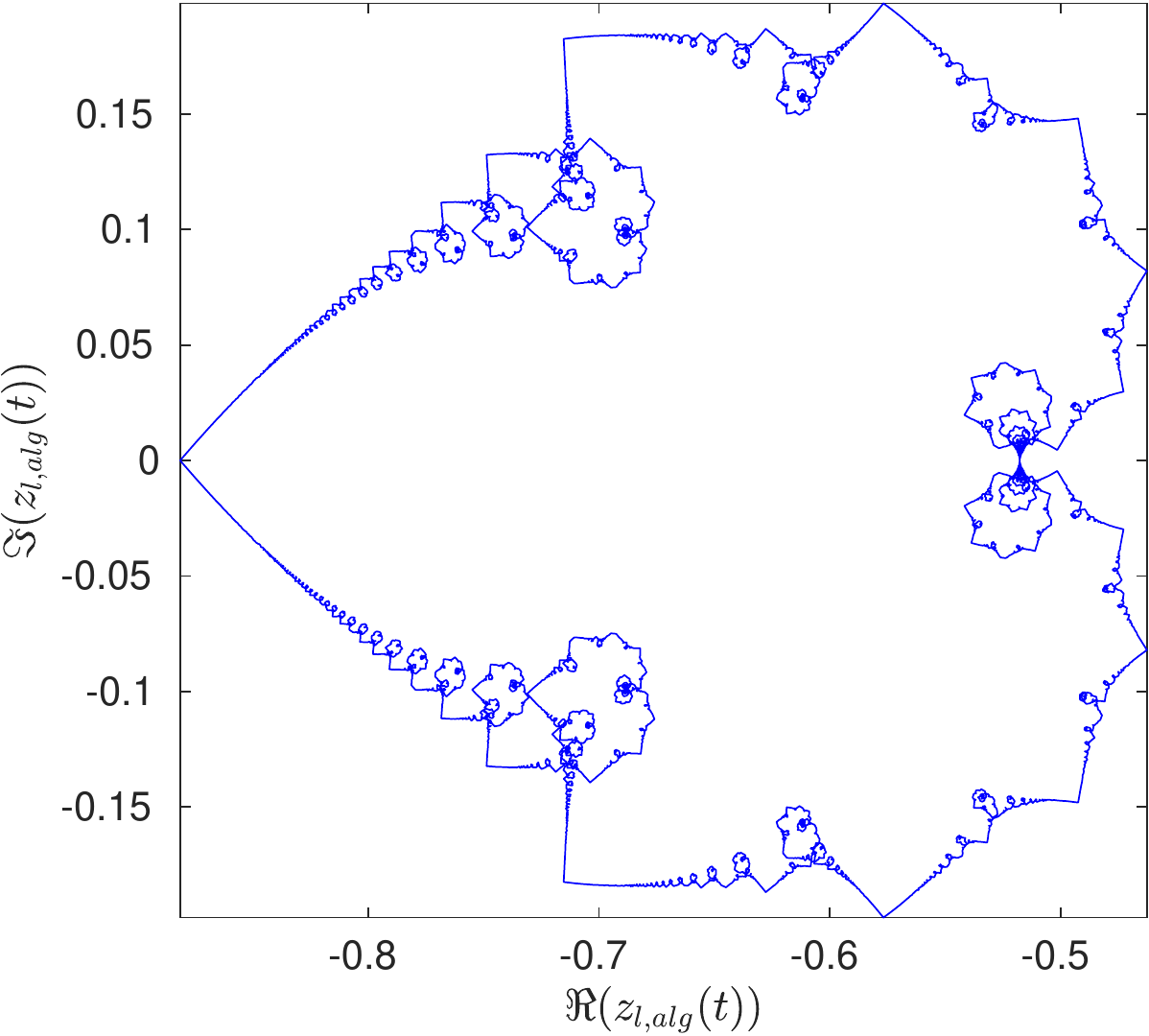}
		\caption{$l=1.7627\ldots$, i.e., $M_e=4$.}
		\label{subfig:zMtB}   
	\end{subfigure}
	\hfill
	\begin{subfigure}[b]{0.333\textwidth}
		\centering
	\includegraphics[width=\textwidth, clip = true]{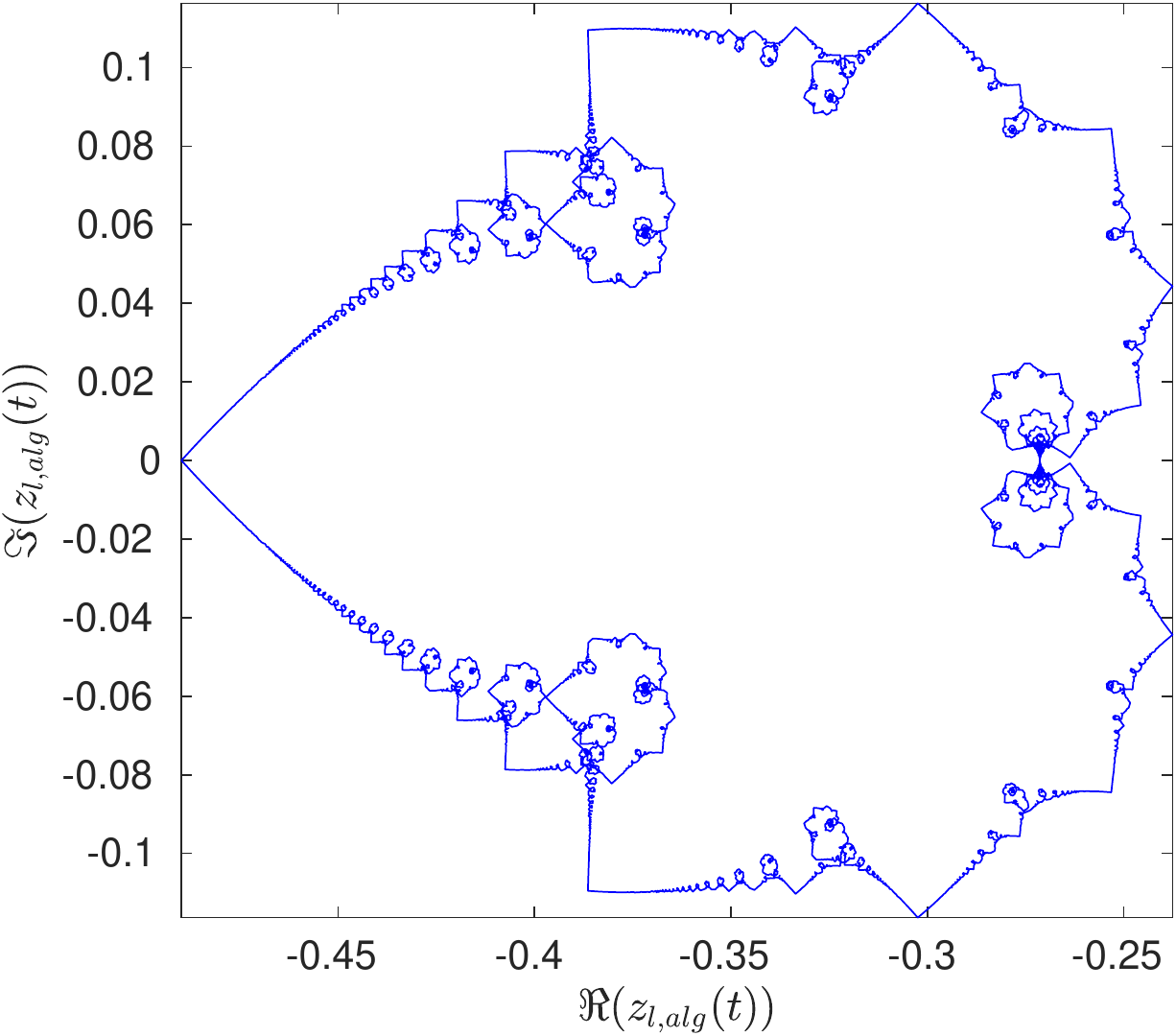}     
		\caption{$l=1.3486\ldots$, i.e., $M_e=5$.}
		\label{subfig:zMtC}   
	\end{subfigure}	
	\begin{subfigure}[b]{0.3217\textwidth}
		\centering
		\includegraphics[width=\textwidth, clip=true]{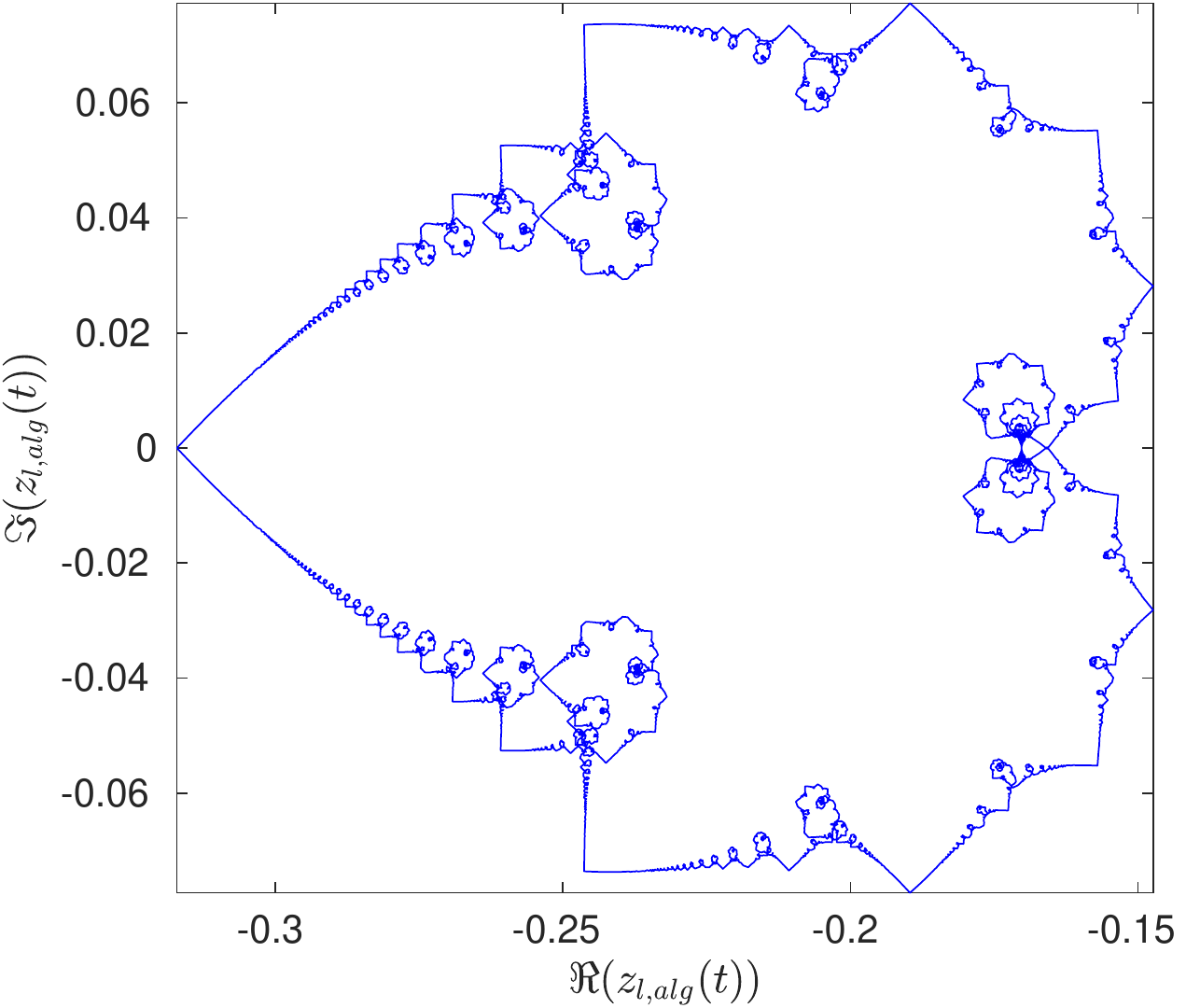}
		\caption{$l=1.0986\ldots$, i.e., $M_e=6$.}
		\label{subfig:zMtD}   
	\end{subfigure}	
	\hfill
	\begin{subfigure}[b]{0.3227\textwidth}
		\centering
	\includegraphics[width=\textwidth, clip=true]{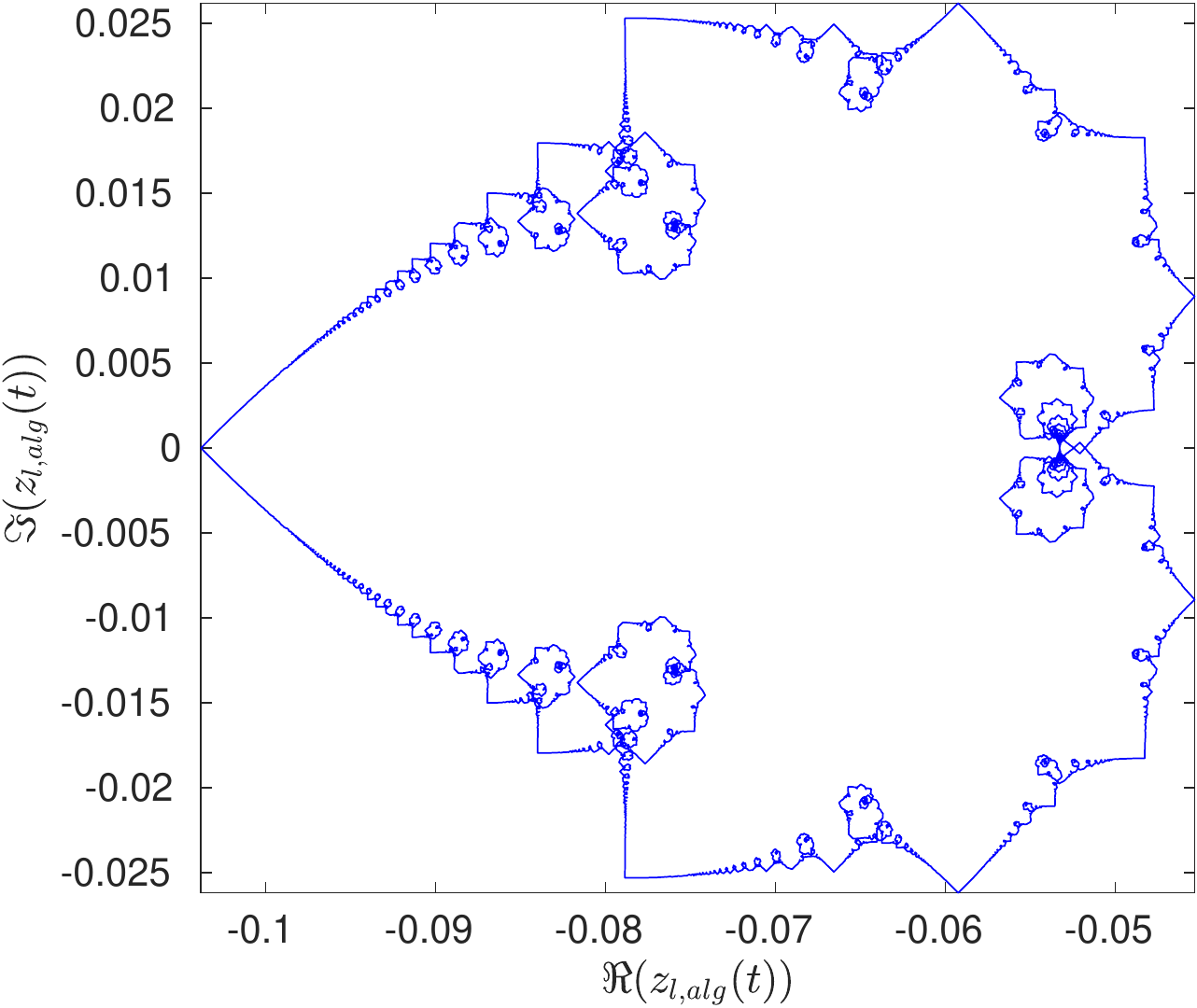}
		\caption{$l=0.6389\ldots$, i.e., $M_e=10$.}
		\label{subfig:zMtE}   
	\end{subfigure}	
	\hfill
	\begin{subfigure}[b]{0.3156\textwidth}
		\centering
	\includegraphics[width=\textwidth, clip = true]{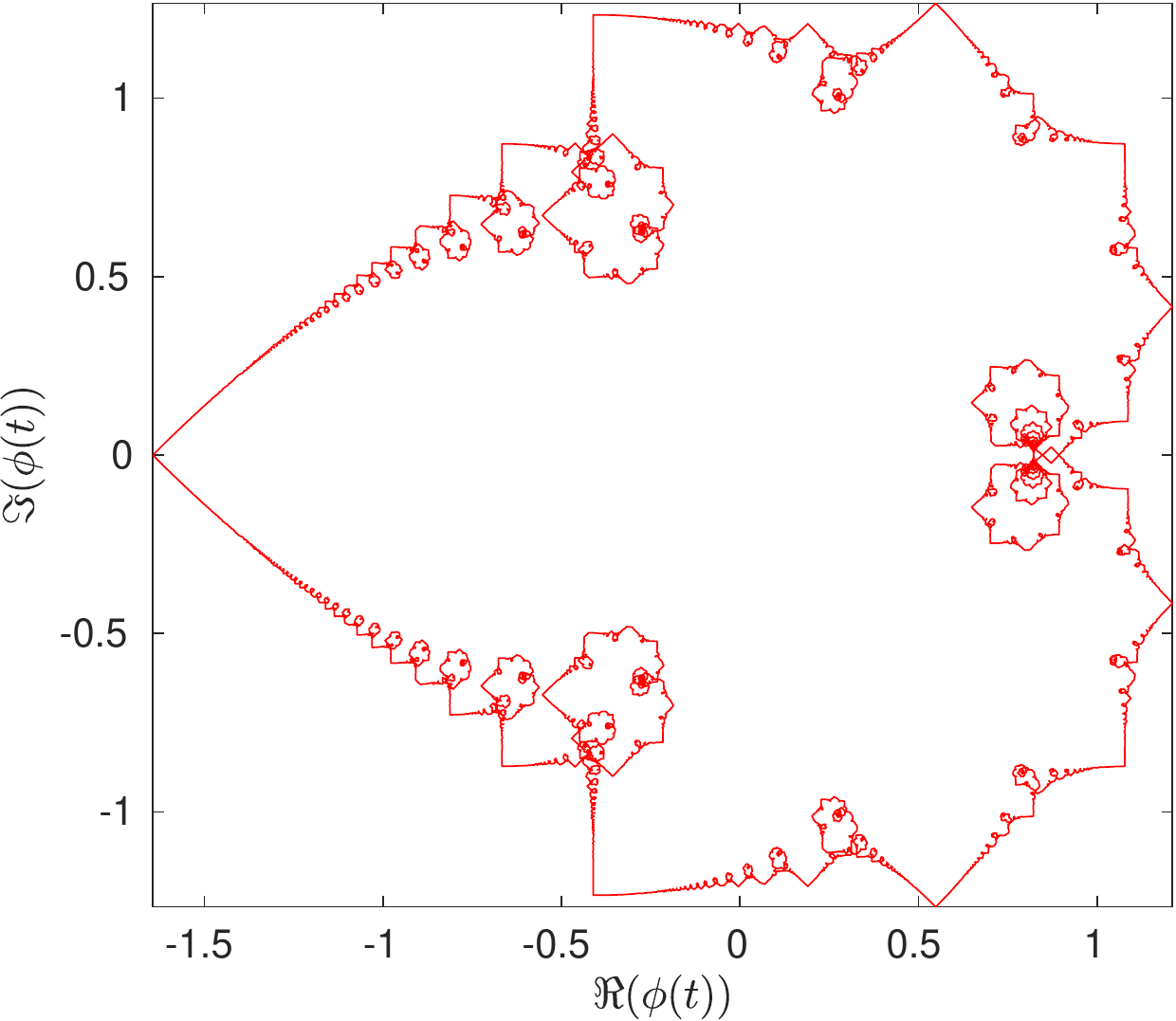}
		\caption{$\phi(t)$.}
		\label{subfig:zMtF}   
	\end{subfigure}
	\caption{$z_{l,alg}(t)$ (blue) as in \eqref{eq:z_lalg}, and $\phi(t)$ (red) as in \eqref{eq:phi}. $z_{l,alg}(t)$ has been generated for different values of $l$, computed using \eqref{eq:choice-of-l}, for $M_e=3,4,5,6$ and $10$. The respective intervals of the parameter $t$ are divided into $7561$ points.}
	\label{fig:zMt}
\end{figure}
From now on, the purpose of working with $\X(0,t)$ will be twofold. First, we would like to see its dependence on the parameter $l$, and second, we would want to compare its structure with the one for the $M$-sided regular planar polygon in the Euclidean space, with $M\geq3$ \cite{HozVega2014} (in the rest of Section \eqref{sec:X0t}, in order to avoid any confusion, we use $M_e$ rather than $M$ to refer to the number of sides of such polygons). In order to address the latter issue, we choose the value of $l$ such that the corresponding parameter $c_0$ is kept the same in both problems. More precisely, from \eqref{eq:c0_exp_hyp}, \cite[(4)]{BanicaVega2018},
\begin{align}
\label{eq:choice-of-l}
\left[\frac{2}{\pi}\ln\left(\cosh\left(\frac{l}{2}\right)\right)\right]^{1/2}= \left[-\frac{2}{\pi}\ln\left(\cos\left(\frac{\pi}{M_e}\right)\right)\right]^{1/2} \Longleftrightarrow l = 2\arccosh\left(\sec\left(\frac{\pi}{M_e}\right)\right).
\end{align} 
Figure \ref{fig:zMt} shows $z_{l,alg}(t)$ (blue) for the values of $l$ corresponding to $M_e=3,4,5,6$ and $10$, and $\phi(t)$ (red); the respective intervals for $t$ have been divided into $7561$ points. Observe that, although similar, the shape of $z_{l,alg}$, for $M_e=3$, is different from its Euclidean counterpart (see \cite[Figure 3]{HozVega2014}); this is discussed further in Section \ref{sec:frthr-rmrks}. On the other hand, except for a scaling, the $z_{l,alg}(t)$ corresponding to the value $M_e=10$ looks very close to $\phi(t)$. In order to further compare the two for different values of $l$, we compute $\phi-\lambda_l z_{l,alg} - \mu_l$, where $\lambda_l \in \mathbb{R}$  and $\mu\in\mathbb{C}$ are given by the least squares fitting method; more precisely,
\begin{equation}
\label{eq:lambda-mu}
\begin{cases}
\lambda_l = \Re\left(\dfrac{\mean[(z_{l,alg}(t)-\mean(z_{l,alg}(t))) (\bar{\phi}(t)-\mean(\bar{\phi}(t) )) ]} {\mean(|\phi(t)-\lambda_l\mean(z_l(t))|^2)} \right), \\ \mu_l = \mean(\phi(t)) - \lambda_l(z_{l,alg}(t)).
\end{cases}
\end{equation}
Thus, for $l$ in \eqref{eq:choice-of-l} corresponding to $M_e=3,4,\ldots,20$, Subfigure \ref{subfig:phi-zl-errorA} shows a log-log plot of the scaling factor $\lambda_l$ in \eqref{eq:lambda-mu}, which behaves linearly with respect to $l$. More precisely, it can be approximated with a straight line $-1.989\,l + 2.994$, shown in red in the same plot. This allows us to claim that the scale of $z_{l,alg}$, when compared to $\phi$, decreases as $l^{-2}$, as $l$ tends to zero. Moreover, in Subfigure \ref{subfig:phi-zl-errorB}, we plot the absolute error $\max_t|(\phi(t) - \lambda_l  z_{l,alg}(t) - \mu_l)|$ (in circled points) and the relative error $\max_t|(\phi(t) - \lambda_l  z_{l,alg}(t) - \mu_l)/ \phi(t)|$ (in starred points), where the maximum is taken over $7560$ values. Clearly, as $l$ gets smaller, both errors decrease, thus, showing that $z_{l,alg}$ converges to $\phi$. Finally, to illustrate the comparison visually, in Subfigures \ref{subfig:zMphiA}, \ref{subfig:zMphiB}, we have plotted scaled $z_{l,alg}$ (blue) superimposed on $\phi$ (red), for $M_e=3, 10$, respectively.	
\begin{figure}[htbp!] 
	\centering
	\begin{subfigure}[t]{0.493\textwidth}
		\centering
		\includegraphics[width=\textwidth]{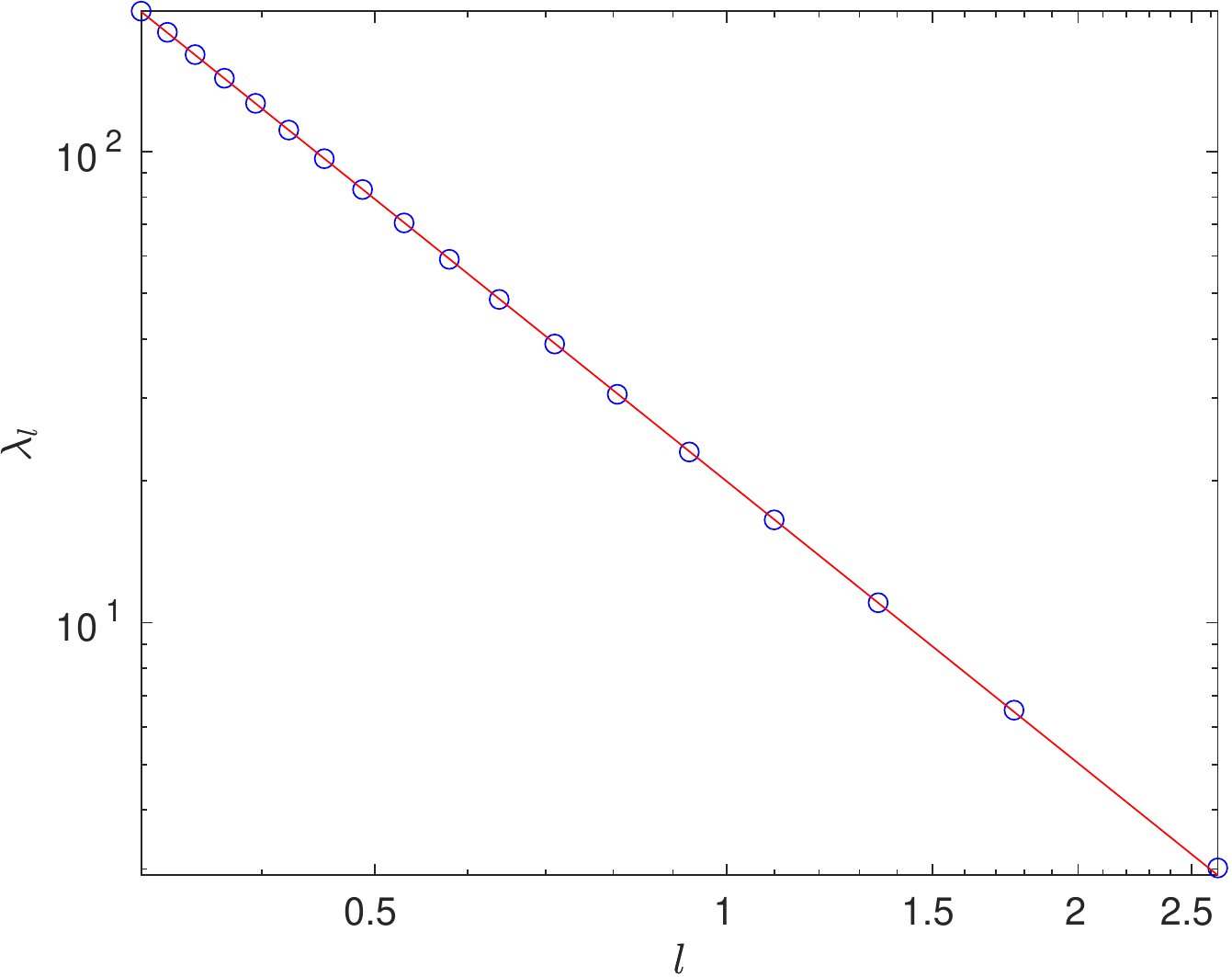}
		\caption{A log-log plot of $\lambda_l$ in \eqref{eq:lambda-mu}, for different values of $l$ as in \eqref{eq:choice-of-l} , which can be approximated with a straight line $-1.989 \,l +2.994$ (in red). This implies that, when compared with $\phi(t)$, the size of $z_{l,alg}$ becomes smaller with a rate tending to $l^{-2}$. } 
		\label{subfig:phi-zl-errorA}   
	\end{subfigure}	
\hfill
	\begin{subfigure}[t]{0.467\textwidth}
		\centering
	\includegraphics[width=\textwidth]{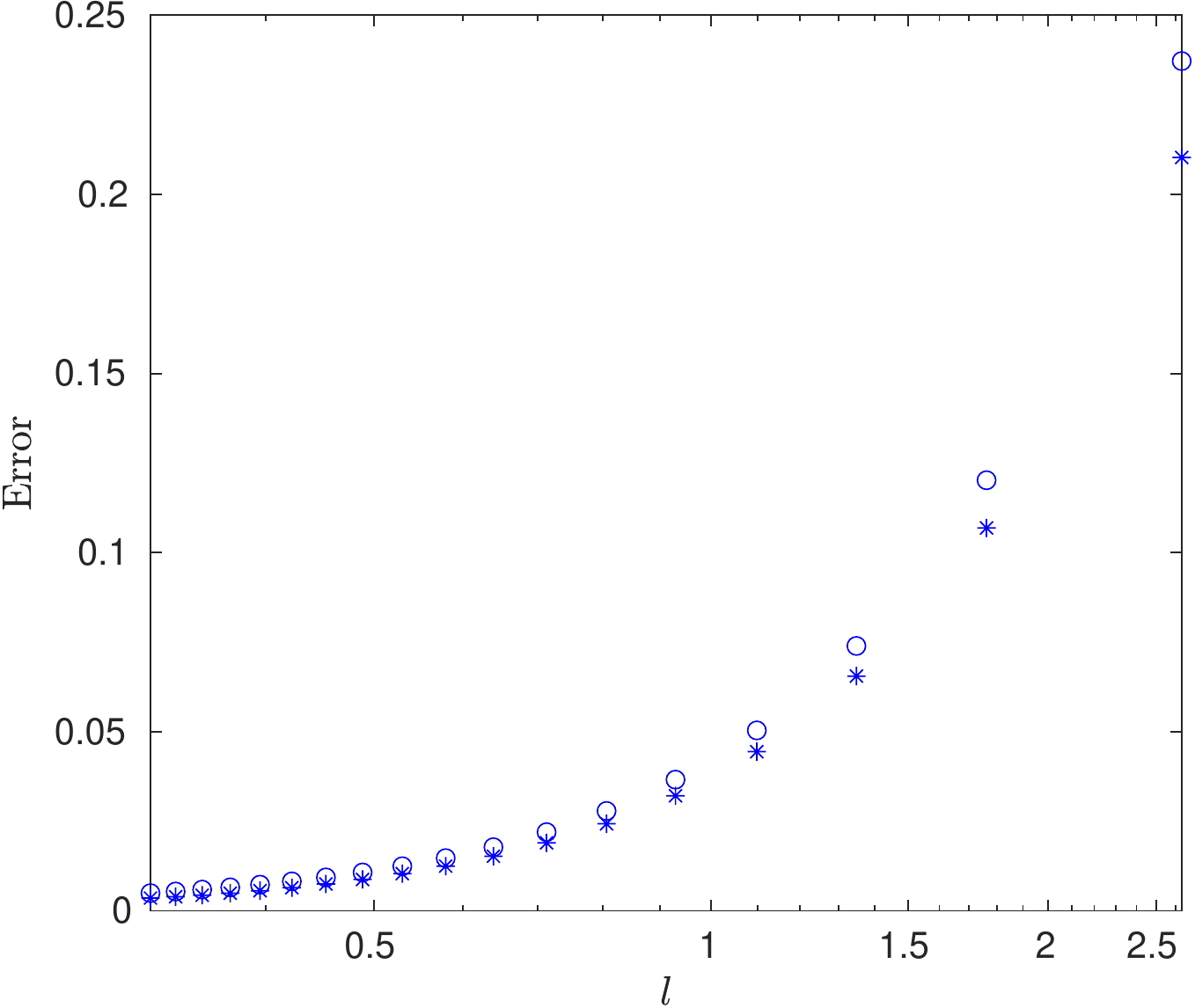}
		\caption{The maximum relative error (starred) and the absolute error (circled) between $\phi(t)$ in \eqref{eq:phi} and $z_{l,alg}(t)$, for different values of $l$. The sum for $\phi(t)$ is taken over $k=1,2,\ldots,2048$, and both $\phi(t)$ and $z_{l,alg}(t)$ have been evaluated at $7561$ points.}
		\label{subfig:phi-zl-errorB}   
	\end{subfigure}		
	\caption{}
		\label{fig:phi-zl-error}
\end{figure}
\begin{figure}[htbp!] 
	\centering
	\begin{subfigure}[b]{0.475\textwidth}
		\centering
		\includegraphics[width=\textwidth]{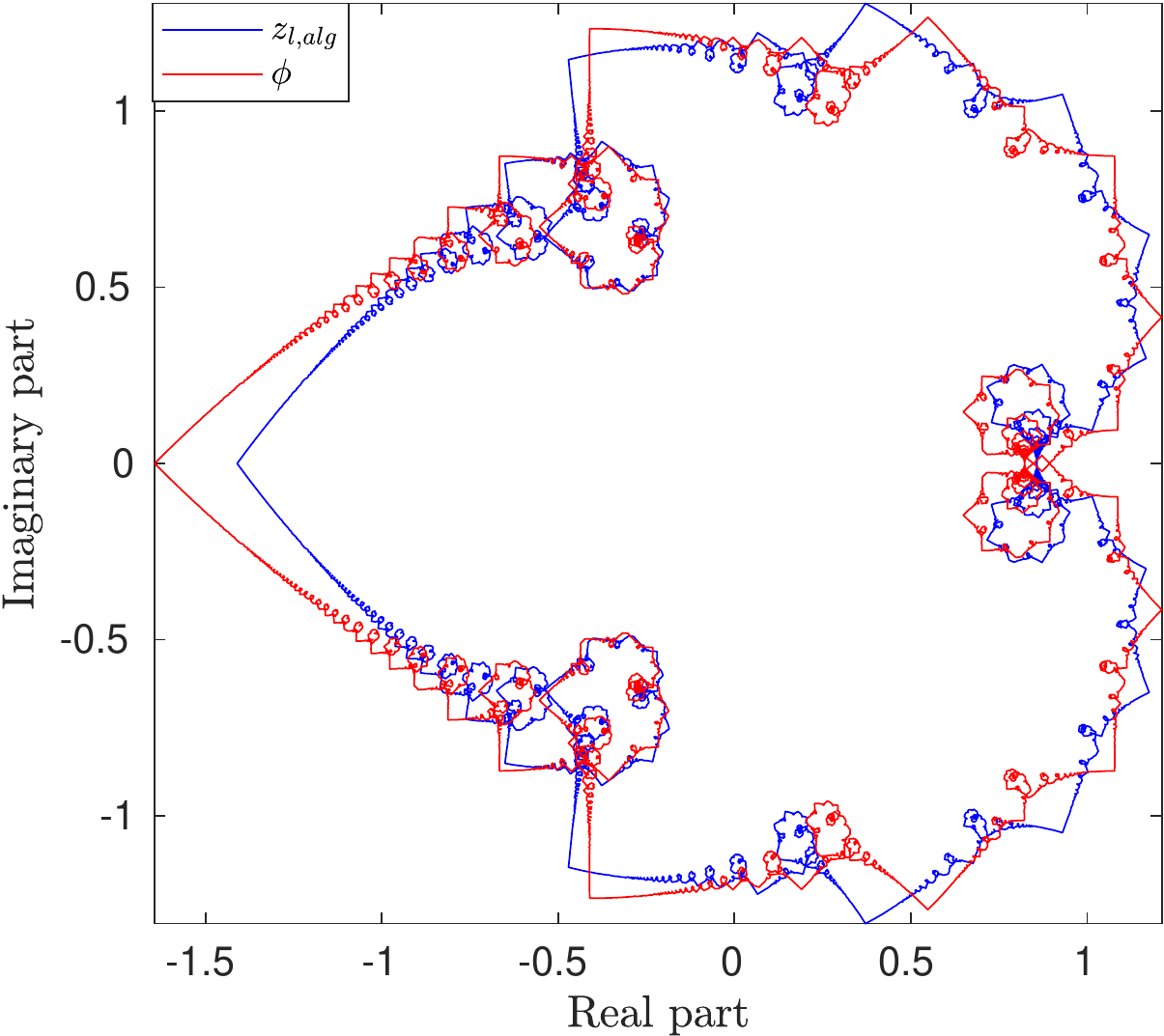}
		\caption{$l=2.6339\ldots$, i.e., $M_e=3$.}
		\label{subfig:zMphiA}   
	\end{subfigure}
\hfill
	\begin{subfigure}[b]{0.485\textwidth}
		\centering
		\includegraphics[width=\textwidth]{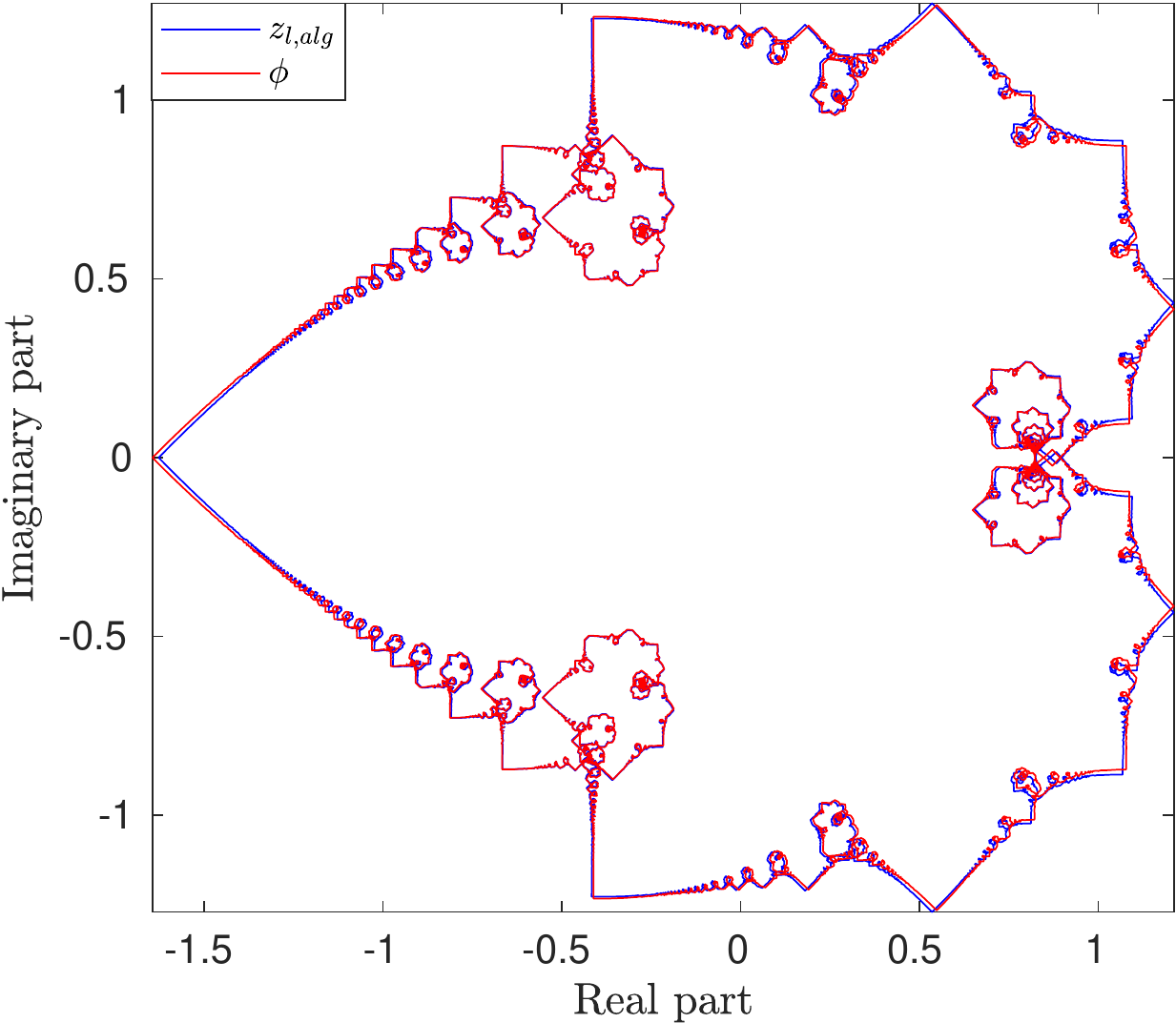}
		\caption{$l=0.6389\ldots$, i.e., $M_e=10$.}
		\label{subfig:zMphiB}   
	\end{subfigure}	
	\caption{A comparison of $\phi(t)$ (red) and the scaled $z_{l,alg}(t)$ (blue) for two different $l$ values. Clearly, as $l$ decreases, the scaled $z_{l,alg}(t)$ tends to $\phi(t)$.}
	\label{fig:zMphi}
\end{figure}
% --------------------------------------------------------------------------------
\subsection{$\T(s,\tpq), \ q\gg1$}
\label{sec:T-irr}
\begin{figure}[htbp!] 
	\centering
	\begin{subfigure}[b]{\textwidth}
		\centering
	\includegraphics[width=\textwidth]{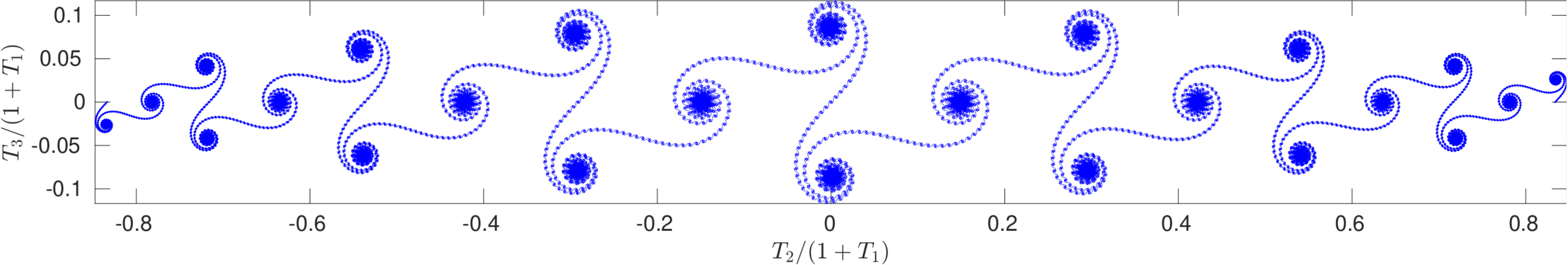}
		\caption{$p=8002$, $q=23997$.}
		\label{subfig:Talg-IrrA}   
	\end{subfigure}
	\begin{subfigure}[b]{\textwidth}
		\centering
	\includegraphics[width=\textwidth]{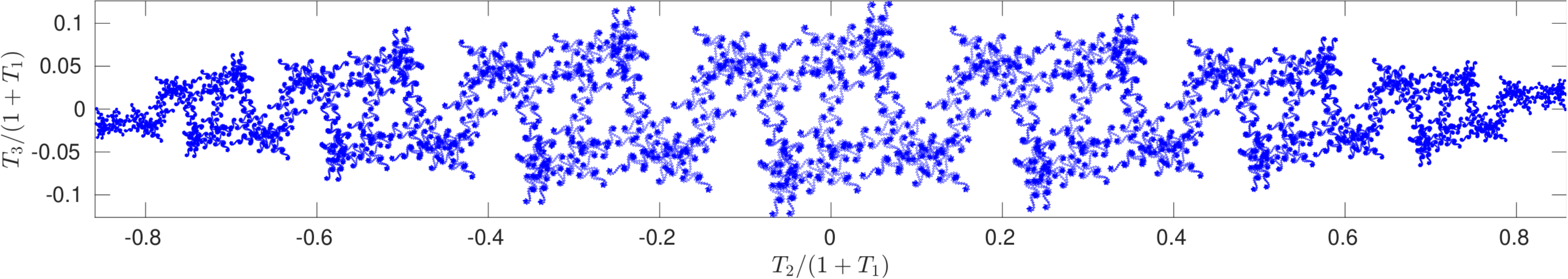}
		\caption{$p=10327$, $q=27993$.	}
		\label{subfig:Talg-IrrB}   
	\end{subfigure}	
	\caption{The stereographic projection of $\Talg(s,\tpq)$ onto $\mathbb{C}$, for $M=8$, $l=0.6$.}
	\label{fig:Talg-Irr}   
\end{figure}

% --------------------------------------------------------------------------------
Having observed the evolution of regular $M$-polygons in the Euclidean case at rational times $\tpq$, with $q\gg1$, we are curious about the behavior of a planar $l$-polygon \cite{HozVega2014,HozKumarVega2019}. In this respect, as in \cite{HozVega2014}, we have examined two cases; first, we consider $\tpq$ with a \textit{small} $q$ and compute the evolution at $t=\tpq+\epsilon$, $|\epsilon|\ll1$.  More precisely, we take $\epsilon = \tfrac{T_f}{q^\prime}$, such that $q^\prime\gg1$, $\gcd(q,q^\prime)=1$, and $\tfrac{p}{q}+\tfrac{1}{q^\prime} = \tfrac{pq^\prime+1}{qq^\prime}$. Therefore, at $\tpq+\epsilon$, there will be $qq^\prime$ or $qq^\prime/2$ times as many sides. We consider the stereographic projection of $\Talg$, projecting it from $(-1,0,0)$ onto the complex plane $\mathbb{C}$; Subfigure \ref{subfig:Talg-IrrA} shows it for $M=8$, $l=0.6$, $p=1$, $q=3$, $q^\prime=7999$. Note that $8\times 23997$ values of the tangent vector form spiral-like structures whose center corresponds to the values of $\T$ at $t=T_f/3$. These spirals can be compared with the Cornu spiral which also appeared in \cite{delahoz2007,HozVega2014}. Next, we look at the rational times $\tpq$, with a \textit{large} $q$, such that there is no pair $\tilde p, \tilde q$, with both $\tilde q$ and $|\tfrac{p}{q}-\tfrac{\tilde p}{\tilde q} |$ being small. In particular, for the same parameters as before, we have taken $t=\left( \tfrac{1}{3} + \tfrac{1}{31}+\tfrac{1}{301}\right) T_f = \tfrac{10327}{27993}T_f$. The stereographic projection of $\Talg$ is shown in Sufigure \ref{subfig:Talg-IrrB}, where the spiral structures at a smaller scale can be observed, thus exhibiting a fractal-like phenomenon.

% --------------------------------------------------------------------------------
\section{Relationship between the $l$-polygon and one-corner problems}
\label{sec:rel-l-1-corner}

Following the approach in the Euclidean case, we conjecture that, at infinitesimal times, the $l$-polygon problem can be seen as a superposition of several one-corner problems \cite{HozVega2018}. In order to compare them, we solve the	one-corner problem for $t=t_{1,q}$, $q\gg1$, and rotate it in such a way that it is oriented with respect to the $l$-polygon problem. We denote the resulting solution by $\X_{rot}$ and $\T_{rot}$, where $\X_{rot}=\mathbf{K}\cdot\X_{c_0}$, $\T_{rot}=\mathbf{K}\cdot\T_{c_0}$, with $\X_{c_0}$, $\T_{c_0}$ being the solution of the one-corner problem, for some rotation matrix $\mathbf K$. Recall that 
\begin{equation*}
\lim\limits_{s\rightarrow -\infty} \T_{c_0} = \A^-= (A_1,-A_2,-A_3)^T, \  \lim\limits_{s\rightarrow \infty} \T_{c_0} = \A^+= (A_1,A_2,A_3)^T,  \ 
\end{equation*}
where $-A_1^2+A_2^2+A_3^2=-1,$ and they are given by \eqref{eq:A1}, \eqref{eq:A2}, \eqref{eq:A3}. Thus, the matrix $\mathbf{K}$ can be computed by enforcing that $\lim_{s\rightarrow \pm \infty} \T_{rot}(s)$ corresponds to the tangent vector of the $l$-polygon at $s=0^\pm$, $t=0$:
\begin{align*}
\lim\limits_{s\rightarrow -\infty} \T_{rot}(s) =  \left( \cosh \left(l/2 \right) , -\sinh \left( l/2\right), 0\right)^T,  \  \lim\limits_{s\rightarrow  \infty} \T_{rot}(s) &=  \left( \cosh \left(l/2 \right) , \sinh \left(l/2\right), 0\right)^T.
\end{align*}
Furthermore, 
\begin{equation}
\label{eq:Xrot-Trot-def}
\begin{aligned}
\X_{rot} &\equiv
(X_{rot,1}, X_{rot,2}, X_{rot,3} )^T
=\mathbf{K} \cdot
(X_{c_0,1}, X_{c_0,2}, X_{c_0,3} )^T
+\frac{l/2}{\sinh(l/2)} \X(0,0), \\ 
\T_{rot} &\equiv
(T_{rot,1},T_{rot,2}, T_{rot,3})^T
=\mathbf{K} \cdot
(T_{c_0,1}, T_{c_0,2},T_{c_0,3})^T,
\end{aligned}
\end{equation}
where $\X(0,0)$ corresponds to the location of the corner of the planar $l$-polygon in \eqref{eq:X-ini-hyp}. 
\subsection{Numerical experiments}
\label{sec:num-exp-l-1-rel}
To solve the two problems numerically, depending on whether $q$ is even or odd, we have different discretizations; however, we restrict ourselves to the case when $q/2$ is odd, as the other two cases can be addressed in a similar way \cite{HozVega2018}. For the $l$-polygon problem, given a value of $q$, we compute the algebraic solution $\Talg(s,t_{1,q})$ at those $s=s_k \in [-l/2,l/2]$ which belong to the middle points of the sides of the corresponding hyperbolic polygon. Thus, for $q/2$ odd, $s_k = 2l k/q$, $k=-(q-2)/4, \ldots, (q-2)/4$, $\Delta s= 2l/q$. On the other hand, after discretizing the interval $[-l/2,l/2]$ with a step size $\Delta s = l/2^4q$, we solve the one-corner problem numerically. 
In this way, $\T_{c_0}(s,t_{1,q})$ can be computed for the same $s=s_k$ as in the $l$-polygon problem, and then, from \eqref{eq:Xrot-Trot-def}, we obtain $\T_{rot}(s,t_{1,q})$. 
\begin{figure}[htbp!] 
	\centering
	\begin{subfigure}[b]{0.32\textwidth}
		\centering
		\includegraphics[width=\textwidth,valign=t]{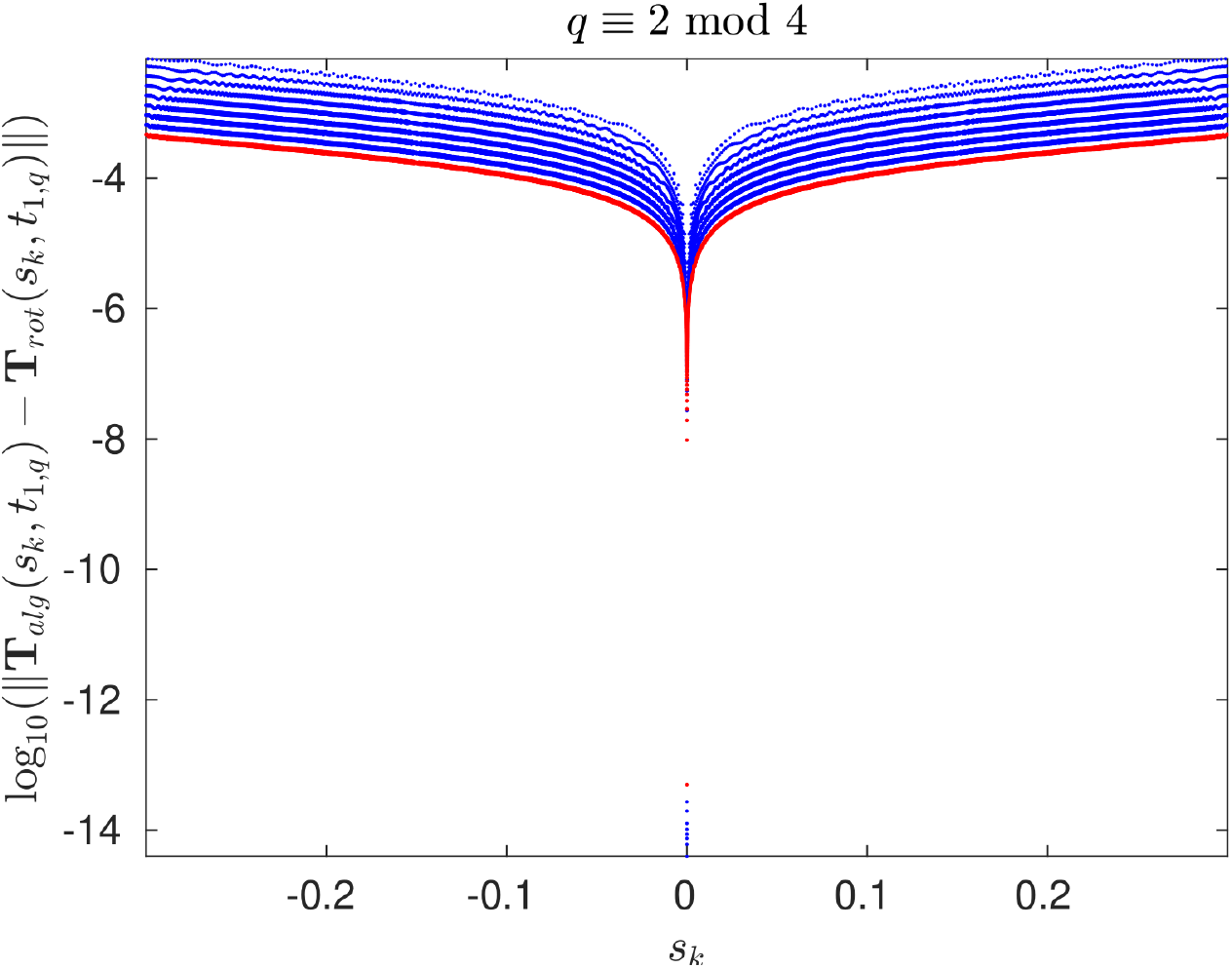}
		\caption{Error between $\Talg$ and $\T_{rot}$ for $q=502, 1002, \ldots,64002$ (in blue), $q=128002$ (in red), $l=0.6$.}
		\label{subfig:Error-mult_one_cornerA}   
	\end{subfigure}			
	\hfill		
	\begin{subfigure}[b]{0.64\textwidth}
		\centering
		\includegraphics[width=\textwidth]{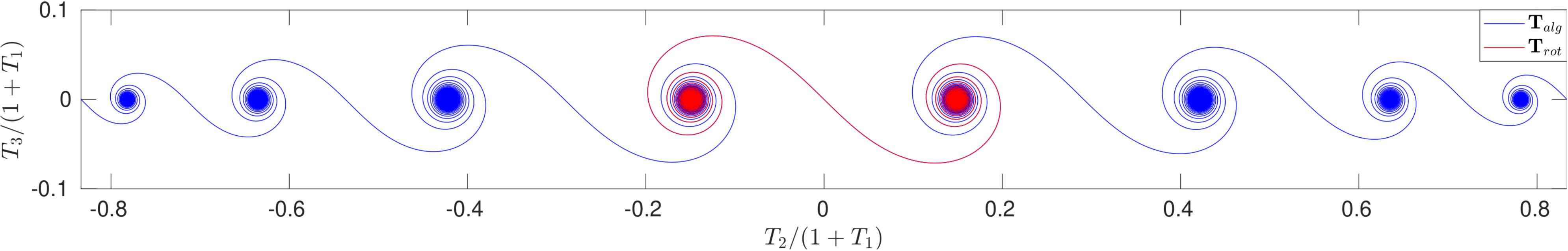}
		\caption{The stereographic projection of $\Talg$ and $\T_{rot}$ onto $\mathbb{C}$ at $t=t_{1,q}$, $q=64002$, $M=8$, $l=0.6$. The two curves are visually indistinguishable from each other. }
		\label{subfig:Talg-tpqB}   
	\end{subfigure}		
%	\caption{}
	\label{fig:Error-mult_one_corner}   
\end{figure}

Subfigure \ref{subfig:Error-mult_one_cornerA} shows the error $\log_{10}(\left \|\Talg(s_k,t_{1,q})-\T_{rot}(s_k,t_{1,q}) \right \|)$ against $s_k$, for $M=8, l=0.6$, and $q=502,1002,2002,\ldots,128002$, where the Euclidean distance $\| \cdot \|$ is computed for each $s_k$. Note that, for a given $q$, the minimum error is attained at $s=0$; and, in general, the best results (in red)  are gotten when the largest value of $q$ is taken. It is possible to check that the maximum of the errors taken over all the values of $s_k$ decreases as $\mathcal{O}(1/\sqrt{q}) = \mathcal{O}(t_{1,q})$. Subfigure \ref{subfig:Talg-tpqB} shows simultaneously the stereographic projection of $\Talg$ and $\T_{rot}$ onto $\mathbb{C}$, at $t=t_{1,q}$, $q=64002$; remark that the red curve is visually indistinguishable from the blue one. 

We can also recover the coefficient $c_0$ in \eqref{eq:cur-ini-pol-hyp}. In order to do it, from the one-corner problem, we write the curvature at $s=0$ and $t>0$ as $c_0(t)=\sqrt{t} | \Ts(0,t)|_0$ \cite{delahoz2007,DelahozGarciaCerveraVega09}. Then, as in \cite{HozVega2018}, at $t=t_{1,q}$, we approximate the derivative with respect to $s$ using a finite difference. Without loss of generality, after taking $q\equiv2\bmod4$, we write
\begin{equation}\label{eq:c0-approximation}
c_0= \lim\limits_{\stackrel{q\rightarrow \infty}{q\equiv2\bmod4}} \sqrt{t_{1,q}} \frac{\left |\Talg(2l/q, t_{1,q}) - \Talg(-2l/q, t_{1,q}) \right |_0 }{4l/q},
\end{equation}
where $\Talg(s,t_{1,q})$ is continuous at $s=0, 2l/q, -2l/q$. Next, using \eqref{mat:Hm}, we obtain
\begin{equation*}
\begin{cases}
\Talg(2l/q,t_{1,q}) = \left(\cosh (l_q), \cos (\theta_1) \sinh (l_q), \sin (\theta_1) \sinh (l_q) \right)^T, \\
\Talg(-2l/q,t_{1,q}) = \left(\cosh (l_q), -\cos (\theta_{q-1}) \sinh (l_q), -\sin (\theta_{q-1}) \sinh (l_q) \right)^T,
\end{cases}
\end{equation*}
and, by substituting them in \eqref{eq:c0-approximation} and computing the limit, we get $c_0$ as in \eqref{eq:c0_exp_hyp} (see \cite[Section 2]{HozVega2018} for the intermediate steps). In Table \ref{table:c0approximation}, we display the error between $c_0$ and its approximated value, computed using \eqref{eq:c0-approximation}, for $l=0.6$ and different values of $q$. Clearly, the error reduces as $\mathcal{O}(1/q)$, thus, showing a complete agreement between the two.

\begin{table}
	\centering
	\begin{tabular}{l c c c r c}
		\hline
		$q$ & Error & $q$ & Error &$q$ & Error \\ 
		\hline
		$502$ & $4.4527 \cdot 10^{-5}$ & $4002$ & $5.5847\cdot 10^{-6}$ & $32002$& $6.9837\cdot 10^{-7}$ \\		
		$1002$ & $2.2306\cdot 10^{-5}$ & $8002$ & $2.7930\cdot 10^{-6}$ & 		$64002$& $3.4920\cdot 10^{-7}$\\		
		$2002$ & $1.1164\cdot 10^{-5}$ & $16002$& $1.3967\cdot 10^{-6}$ & $128002$ & $1.7461\cdot 10^{-7}$\\
		\hline
	\end{tabular}
	\caption{The error $|c_0-\sqrt{t_{1,q}} \left |\Talg(\Delta s,t_{1,q})-\Talg(-\Delta s,t_{1,q}) \right |_0 / (2\Delta s)|$, where $c_0=0.1680\ldots$, $l=0.6$, $\Delta s = 2l/q$.}
	\label{table:c0approximation}
\end{table}

\begin{figure}[htbp!] 
	\centering
	\begin{subfigure}[t]{0.475\textwidth}
		\centering
		\includegraphics[width=\textwidth,clip=true, valign =t]{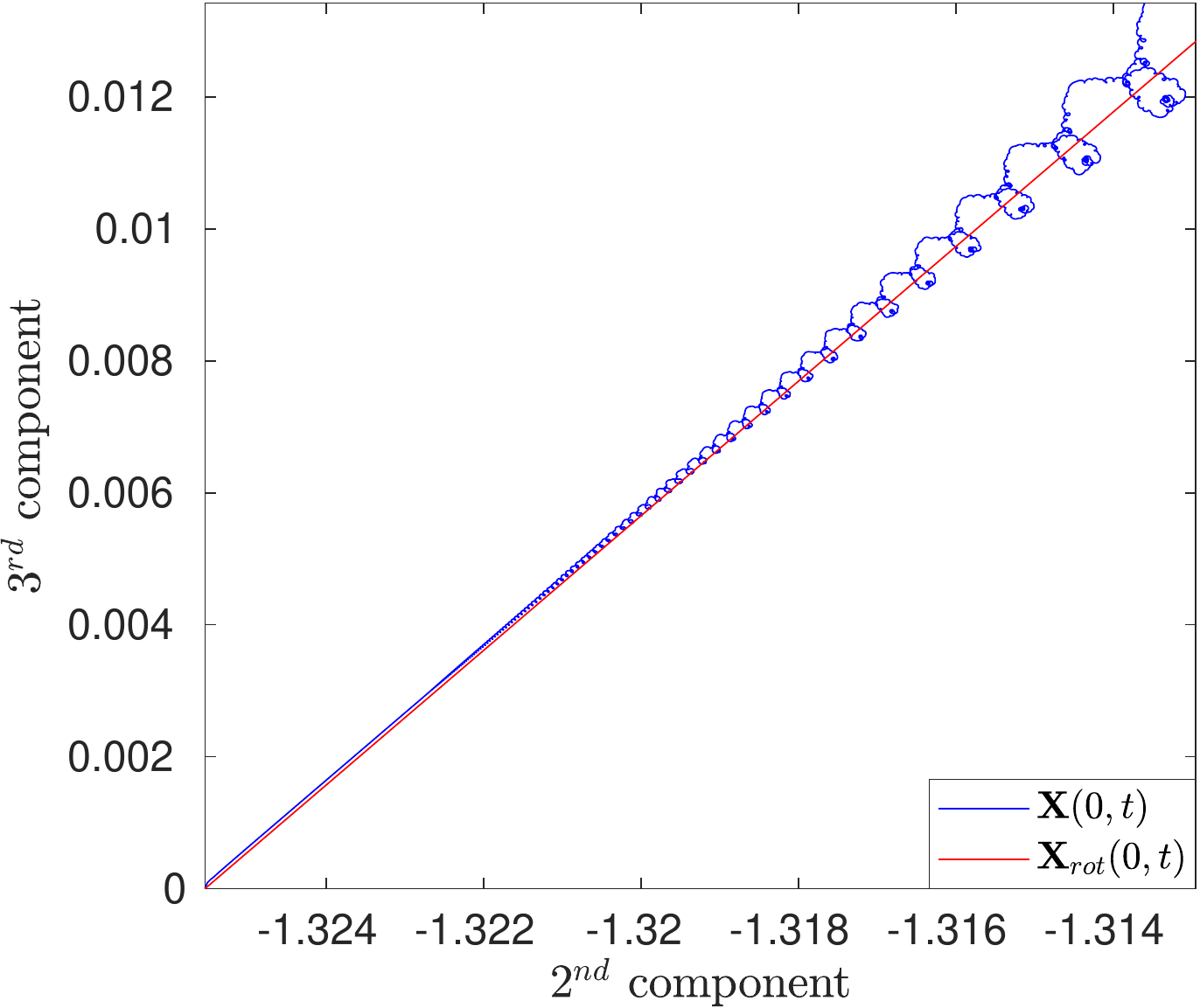}
		\caption{The time evolution for short times can be well approximated with a straight line with slope $A_2/A_3$. Here, we have taken $M = 8$, $l = 0.6$, $t\in [0,t_{1,20}]$. }
		\label{subfig:M-1-XtrajA}   
	\end{subfigure}	
	\hfill	
	\begin{subfigure}[t]{0.485\textwidth}
		\centering
		\includegraphics[width=\textwidth,clip=true, valign =t]{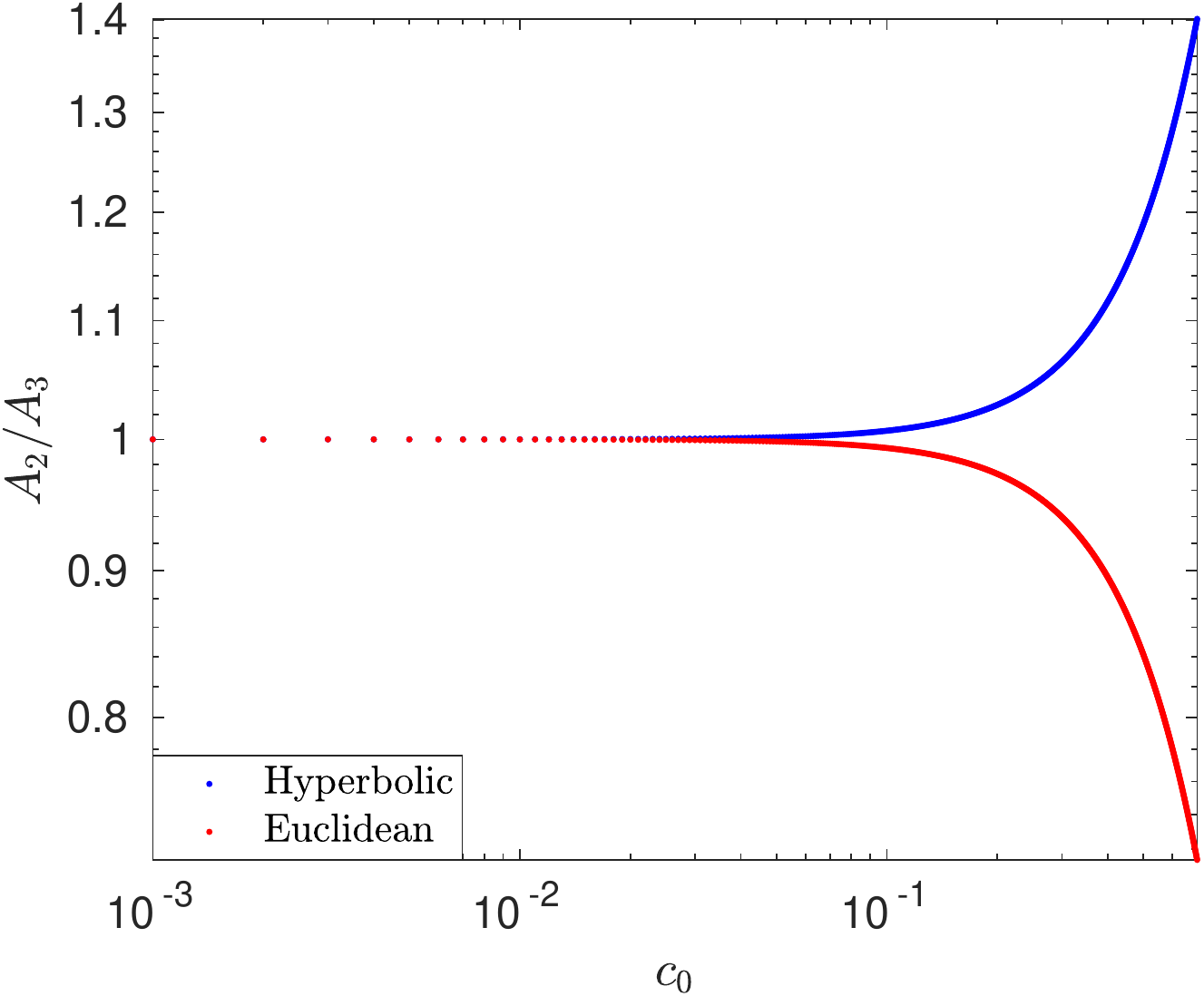}
		\caption{A log-log plot of $A_2/A_3$ as a function of $c_0$. Clearly, $A_2/A_3$ tends to $1$, as $c_0$ tends to $0$.}
		\label{subfig:A2A3}   
	\end{subfigure}				
	\caption{}
		
	\label{fig:M-1-Xtraj}
\end{figure}

Next, we compare the time evolution of a point in both problems and compute $\X(0,t)$ and $\X_{rot}(0,t)$, for $t\in [0,t_{1,20}]$. More precisely, using \eqref{eq:ini-cond-1-crnr},
\begin{equation}
\begin{aligned}
\label{eq:Xrot-def}
\X_{rot}(0,t) &\equiv
\mathbf{K} \cdot
(X_{c_0,1}(0,0), X_{c_0,2}(0,0), X_{c_0,3}(0,0))^T 
+\frac{l/2}{\sinh(l/2)} \X(0,0) \\
&=2c_0\sqrt{\frac{t}{A_2^2+A_3^2}} (0,A_3,A_2)^T +\frac{l/2}{\sinh(l/2)} \X(0,0).
\end{aligned}
\end{equation}
Both $\X(0,t)$ and $\X_{rot}(0,t)$ lie in the YZ-plane, and after plotting simultaneously their projection onto $\mathbb{C}$, we note that, for small times, $\X(0,t)$ (in blue) can be very well approximated by a straight line (in red) with slope $A_2/A_3$ (Subfigure \ref{subfig:M-1-XtrajA}). 

In our numerical simulations, we have taken $M=8$, $l=0.6$, but the results hold true for any $M\geq2$, where $M$ has been chosen to be even, in order to take advantage of the symmetries of the hyperbolic polygon. Thus, there is strong numerical evidence that, at small times, the $l$-polygon problem can be seen a superposition of several one-corner problems. 

\subsection{Further remarks}
\label{sec:frthr-rmrks}
Note that the quantity $A_2/A_3$ also determines the angle $\varphi$ that the curve $\X_{rot}(0,t)$ makes with the plane containing $\X_{rot}(s,0)$. Interestingly, $\varphi$ is the angle corresponding to the corner of $z_l(t)$ located at $t=0$, and this holds true for the Euclidean case as well. To compare the two cases simultaneously, we have computed $A_2/A_3$ for several different values of $c_0$, using \eqref{eq:A2}--\eqref{eq:A3} in the hyperbolic case, and \cite[(57)]{GutierrezRivasVega2003} in the Euclidean case. The values thus obtained have been plotted in Subfigure \ref{subfig:A2A3}, where it can be observed that $A_2/A_3$ is greater (respectively, smaller) than the one in the hyperbolic (respectively, Euclidean) case, and tends to one, as $c_0$ approaches zero; in fact, from \eqref{eq:A2}--\eqref{eq:A3}, $A_2(0)/A_3(0)=1$. On the other hand, in the hyperbolic case, we have
\begin{equation}
\begin{aligned}
\varphi = \arctan\left(\frac{A_2}{A_3}\right)=\arctan\left(\frac{\Re\{\Upsilon\}}{\Im\{\Upsilon\}}\right)=\arg(i\bar\Upsilon), 
\end{aligned}
\end{equation}
with $\Upsilon=e^{i\pi /4} \Gamma(1-ic_0^2/4) \Gamma(1/2+ic_0^2/4)$. Thus, for a given $c_0$, $\varphi$ is larger (respectively, smaller) than $\pi/2$ and, in the limit, it converges to $\pi/2$, as in the case of Riemann's function. 
%
%\begin{figure}[htbp!] 
%	\centering
%	\includegraphics[width=0.5\textwidth,clip=true, valign =t]{A2A3b}	
%	\caption{Semilogarithmic plot of $A_2/A_3$ as a function of $c_0$. Clearly, as $c_0$ tends to $0$, $A_2/A_3$ tends to $1$ (dashed dotted line).}
%	\label{fig:A2A3}
%\end{figure}

The relationship between the $l$-polygon problem and the one-corner problem has several deep implications, and it is a leap forward in understanding the evolution of corner-shaped initial data. Thanks to this relationship, using similar arguments as in \cite[Section 3.1]{HozVega2018}, the speed of the center of mass $c_l$ can be obtained by computing the integral of $\X_{rot,3}(s/\sqrt{t},1)=\X_{rot,3}(s,t)$. Consequently, we have the following result.
\begin{thm}
	\label{thm:c_l_expression}
	\begin{equation}
	\label{eq:X3rot-integral}
	\int_{-\infty}^{\infty} \X_{rot,3}(s) \ ds =  \frac{2\pi c_0^2}{\sqrt{1-e^{-\pi c_0^2}}}.
	\end{equation}
	Therefore, we can express $c_l$ in terms of $c_0$ (and in terms of $l$):
	\begin{equation}
	\label{eq:c_l-expression}
	c_l = \frac{2\pi c_0^2}{l\sqrt{1-e^{-\pi c_0^2}}} =\frac{4 \ln \cosh (l/2) }{l\sqrt{1-\sech^2(l/2)}}  
	= -\frac{ \ln (1-\tanh^2(l/2)) }{l/2 \tanh(l/2)}.
	\end{equation}
\end{thm}
\noindent We omit the proof as it follows similar steps to those in \cite[Theorem 3.1]{HozVega2018}. Furthermore, by solving \eqref{eq:ini-cond-1-crnr}--\eqref{eq:X-ODE} numerically, as in \cite[Section 3.3]{HozVega2018}, we also have a numerical proof for \eqref{eq:X3rot-integral}--\eqref{eq:c_l-expression}.
\section{Conclusions}
\label{sec:conclusion}
In this paper, we have studied the evolution of \eqref{eq:SMP-hyp}--\eqref{eq:VFE-hyp} for a regular planar $l$-polygon. The motivation to work with such kind of initial data comes in fact from the one-corner problem in the hyperbolic case \cite{delahoz2007}, and recent work on the regular polygons in the Euclidean case \cite{HozVega2014}. In \cite{delahoz2007}, it was observed that, due to the exponential growth of the Euclidean length of the tangent vector, the numerical treatment of the one-corner problem in the hyperbolic case poses restrictions on the value of the parameter $c_0$, and the same is observed in the planar $l$-polygon case as well.
After trying several different numerical schemes, we have concluded that a finite difference scheme with fixed boundary conditions on $\T$ gives the best results, which are also in agreement with their algebraic counterparts. The evolution is periodic in time with a period $l^2/2\pi$, and, at intermediate rational times $\tpq=(l^2/2\pi)(p/q)$, $\gcd(p,q)=1$, depending on the parity of $q$, the polygonal curve has $q$ or $q/2$ times as many sides. As in the Euclidean case, this intermittent behavior of formation/annihilation of the corners can be seen as a nonlinear Talbot effect \cite{HozVega2014,HozVega2018}. 

We have also analyzed the multifractal trajectory of a corner $\X(0,t)$, by comparing it with Riemann's non-differentiable function and its equivalent in the Euclidean case; and this has been supported with adequate numerical experiments. Furthermore, as in \cite{HozVega2018}, we have established a relationship between the one-corner problem and the $l$-polygon problem, and, as a consequence, a precise expression for the speed of the center of mass of an $l$-polygon has been obtained. Finally, we have obtained explicit expressions for the components of the tangent vector $\A^\pm$, whose knowledge has been essential in this work. 
\section*{Acknowledgments}
Sandeep Kumar would like to thank Carlos J. Garc\'{\i}a-Cervera for the discussions on Section \ref{sec:num-sol} which took place during his visit to the University of California, Santa Barbara (UCSB), USA. The authors would like to thank the anonymous referees for their valuable comments and suggestions that have largely improved the presentation of this paper.

This paper was partially supported by the ERCEA Advanced Grant 2014 669689 - HADE, by the MICINN Projects PGC2018-094522-B-I00 and SEV-2017-0718, by the Basque Government Grant IT1247-19, and by the Basque Government BERC Program 2018-2021.
\appendix
\section{Expression for $\A^\pm=(A_1,\pm A_2,\pm A_3)^T$}
\label{apdx}
For the one-corner problem in the Euclidean case, a precise expression for each of the components of the tangent vector $\A^+ = (A_1,A_2,A_3)^T$ was given in \cite{GutierrezRivasVega2003}; and, later in \cite{delahoz2007}, an expression for $A_1$, i.e., \eqref{eq:A1-c0}, was obtained for the hyperbolic case as well. In the following lines, using a completely different approach, we rederive $A_1$ by means of the Laplace transform, and, continuing the calculations in \cite{delahoz2007}, we calculate $A_2$ and $A_3$, whose knowledge has been extremely useful in Section \ref{sec:rel-l-1-corner}.

\subsection{Computation of $A_1$ using the Laplace transform}
\label{sec:A1-comp}
Recall that the proof of Theorem \ref{thm:c_l_expression} mainly involves working with the even solution of \eqref{eq:X-ODE-FT}, which is also analytic. However, in the following lines, we consider the odd solution, which behaves like $\delta^\prime$ near the origin, and can be expressed as
\begin{equation*}
\hat \chi (\xi) = b_0\delta^\prime + b_1 \sign(\xi) + b_2 \sign(\xi) \xi + \ldots,
\end{equation*}
where the coefficients $b_0 = 1, b_1 = - c_0^2,  b_2 = -c_0^4/(1- 4c_0^2), \ldots$, are obtained by introducing $\hat \chi(\xi)$ into \eqref{eq:X-ODE-FT}. If we write the first component of $\hat{\X}(\xi)$ as $\hat{X}_1(\xi) = -i A_1 \hat \chi (\xi)$, and define
\begin{align*}
\WW_1(\eta) &= \WW_1(\xi^2) = \xi^2 \hat{X}_1(\xi), \ \eta > 0,
\end{align*}
then, it solves
\begin{align}
\label{eq:W-hat-ODE}
\WW^{\prime\prime} + \WW\left(1+ \frac{c_0^2}{\eta}\right) = 0 \Leftrightarrow\eta \WW^{\prime\prime} + \eta \WW + c_0^2 \WW = 0,
\end{align}
with
\begin{equation}
\label{eq:W-hat-0-Lap}
\WW_1(0)=0, \ \WW_1^\prime(0) = \lim\limits_{\eta \to 0} \frac{\WW_1(\eta)}{\eta} =  iA_1 c_0^2. 
\end{equation}
\noindent On the other hand, the Laplace transform of $\WW_1(\eta)$, 
\begin{equation}
\label{eq:Lap-t-def}
\mathcal{L}(t) = \mathcal{L}\{\WW_1(\eta)\} = \int_{0}^{\infty} \WW_1(\eta) e^{-t\eta} d\eta, \ t>0,
\end{equation}
satisfies
\begin{equation}
\label{eq:Lap-t-eqn}
t^2 \mathcal{L}^\prime(t) + 2t\mathcal{L}(t) + \mathcal{L}^\prime(t) - c_0^2 \mathcal{L}(t) = 0.
\end{equation}
Furthermore, 
\begin{equation}
\label{eq:Lap0}
\mathcal{L}(0) = \int_{0}^{\infty} \WW_1(\eta) d\eta 
= 2 \int_{0}^{\infty} \xi^3 \hat{X}_1(\xi) d\xi =  \int_{-\infty}^{\infty} \xi^3 \hat{X}_1(\xi) d\xi= i X_1^{\prime\prime\prime}(0)=ic_0^2,
\end{equation}
where we have used the fact that $\hat{X}_1$ is odd. 
Rewriting \eqref{eq:Lap-t-def} as
\begin{equation*}
\mathcal{L}(t) = \int_{0}^{\infty} \WW_1(\eta) e^{-t\eta}  d\eta = \frac 1t \int_{0}^{\infty} \WW_1^\prime(\eta) e^{-t\eta} d\eta 
= \frac{\WW_1^\prime(0)}{t^2} + \frac{1}{t^2} \int_{0}^{\infty} \WW_1^{\prime\prime}(\eta) e^{-t\eta} d\eta,
\end{equation*}
we have
\begin{equation*}
t^2 \mathcal{L}(t)  = \WW_1^\prime(0) + \int_{0}^{\infty} \WW_1^{\prime\prime} e^{-t\eta} d\eta, 
\end{equation*}
which, as $t\to \infty$, becomes
\begin{equation}
\label{eq:Lap-t-inf}
\lim\limits_{t\to\infty} t^2 \mathcal{L}(t) = \WW^\prime(0). 
\end{equation}
Hence, from \eqref{eq:Lap-t-eqn}-\eqref{eq:Lap0}, we have an initial value problem whose solution $\mathcal{L}(t)$ satisfies
\begin{align*}
\lim\limits_{t\to \infty} t^2 \mathcal{L}(t) = \lim\limits_{t\to \infty} t^2\frac{\mathcal{L}(0)}{1+t^2} e^{c_0^2 \arctan(t)} = ic_0^2 e^{c_0^2 \pi /2}.
\end{align*}
Combining this with \eqref{eq:W-hat-0-Lap} and \eqref{eq:Lap-t-inf},  we conclude that
\begin{equation} \label{eq:A1}
ic_0^2 A_1 = \WW^\prime(0) = ic_0^2 e^{c_0^2 \pi /2} 
\implies A_1 = e^{ c_0^2 \pi /2}. 
\end{equation}
The above approach works the same for the Euclidean case as well; hence, the corresponding expression for $A_1$ can be obtained.

\subsection{Computation of $A_2$ and $A_3$}
\label{sec:A2A3-comp}
By continuing the computations of \cite[Theorem 1]{delahoz2007}, we can also obtain the expressions for $A_2$, $A_3$. In this regard, writing them componentwise, the solutions of the Frenet-Serret formulas with $\kappa=c_0$, $\tau = s/2$, i.e., $\T\equiv(T_j)$, $\nn\equiv (n_{j})$, $\bb \equiv (b_{j})$, satisfy
$$
|n_{j}|^2 + |b_{j}|^2 - |T_j|^2 = 
\begin{cases}
-1, & \text{if} \ j = 1, \\ 
1, & \text{if} \ j = 2,3,
\end{cases}
$$
where
$$
\T(0) = (1,0,0)^T, \ \nn(0) = (0,1,0)^T,  \ \bb(0) = (0,0,1)^T.
$$
Recall that, from \cite[Theorem 1]{delahoz2007}, $A_j = \lim\limits_{s\to\infty} T_j(s), \ j =1,2,3,$
with
\begin{equation}
\label{eq:Apdx-T-def}
T_j(s) = i (1+\theta_j \bar{\vartheta_j})(s), \ j=2, 3,
\end{equation}
where $\theta_j$ and $\vartheta_j$ satisfy \cite[(49)]{delahoz2007} and $\theta_j^\prime \bar\vartheta_j^\prime - (c_0^2/4) \theta_j \bar\vartheta_j = E_j$, and can be represented as  
\begin{equation}
\label{eq:Apdx-th-verth-def}
\begin{cases}
\theta_j(s) = a_{1,j} \beta_1(s) + a_{2,j} \beta_2(s), \\
\vartheta_j(s) = b_{1,j} \beta_1(s) + b_{2,j} \beta_2(s),
\end{cases}
\end{equation}
where $E_j$ is chosen later, and $\beta_1(s)$, $\beta_2(s)$ are as in \cite[ (55)]{delahoz2007}. Hence, our first goal is to compute $a_{1,j}, a_{2,j}, b_{1,j}$, $b_{2,j}$, for $j=2,3$.

Differentiating \eqref{eq:Apdx-th-verth-def} gives
\begin{equation}
\label{eq:Apdx-th-verth-1drv-def}
\begin{cases}
\theta_j^\prime(s) = a_{1,j} \beta^\prime_1(s) + a_{2,j} \beta^\prime_2(s), \\
\vartheta_j^\prime(s) = b_{1,j} \beta^\prime_1(s) + b_{2,j} \beta^\prime_2(s),
\end{cases}
\end{equation}
and, from \cite[p. 77]{delahoz2007}, for $j=2,3$, the asymptotics of $\theta_j(s)$ and $\vartheta_j(s)$ are given by 
\begin{equation}
\label{eq:Apdx-asympts-th-varth}
\begin{cases}
\theta_j(s) = (a_{1,j} \gamma_1 + a_{2,j} \gamma_2) e^{-i\frac{c_0^2}{2}\log s} + \mathcal{O}(1/s), \ s \to \infty, \\
\vartheta_j(s) = (b_{1,j} \gamma_1 + b_{2,j} \gamma_2) e^{-i\frac{c_0^2}{2}\log s} + \mathcal{O}(1/s), \ s \to \infty,
\end{cases}
\end{equation}
where 
$$
\gamma_1 = 2e^{-\pi c_0^2/4} \Gamma(1+ic_0^2/2),  \ 
\gamma_2 = -2e^{\pi c_0^2/4} \Gamma(1+ic_0^2/2).
$$
By taking $E_2 = c_0^2/2$, $T_2(0)=0$ in \eqref{eq:Apdx-T-def}, we get
\begin{equation}
\label{eq:Apdx-th-varth-1drv-at0}
\theta_2^\prime(0) \bar\vartheta_2(0) = c_0^2/4,
\end{equation}
and from \cite[(53)]{delahoz2007}
\begin{equation}
n_{j} - i b_{j} = (2i/c_0) \theta_j \bar\vartheta_j^\prime,
\end{equation}
so, if $\theta_2(0)=1$, then, by using \eqref{eq:Apdx-T-def}, \eqref{eq:Apdx-th-varth-1drv-at0}, we obtain 
\begin{equation}
\label{eq:Apdx-th-varth-and-1drv-at0}
\theta_2^\prime(0) = ic_0/2, \ \bar\vartheta_2(0)=-1, \bar\vartheta_2^\prime(0) = -ic_0/2. 
\end{equation}
Thus, by evaluating \eqref{eq:Apdx-th-verth-def}, \eqref{eq:Apdx-th-verth-1drv-def} at $s=0$, and using $\beta_1(0) = -\beta_2(0)$, $\beta_1^\prime(0) = \beta_2^\prime(0)$, and \eqref{eq:Apdx-th-varth-and-1drv-at0},
\begin{equation}
\begin{cases}
(a_{1,2} - a_{2,2}) \beta_1 = 1,  \ (a_{1,2} + a_{2,2}) \beta_1^\prime = ic_0/2, \\
(b_{1,2} - b_{2,2}) \beta_1 = -1,  \ (b_{1,2} + b_{2,2}) \beta_1^\prime = ic_0/2, 
\end{cases}
\end{equation}
where 
\begin{align*}
\beta_1\equiv\beta_1(0) &=  2e^{-\pi c_0^2/8} \Gamma(1+ic_0^2/4), \quad
\beta_1^\prime\equiv\beta_1^\prime(0)  = -(c_0^2/2) e^{i\pi /4} e^{-\pi c_0^2/8} \Gamma(1/2+ic_0^2/4).
\end{align*}
As a result,
\begin{equation}
\label{eq:Apdx-A2-coef}
a_{1,2} = \frac{ic_0\beta_1 + 2\beta_1^\prime}{4\beta\beta_1^\prime}, \ a_{2,2} = \frac{ic_0\beta_1 - 2\beta_1^\prime}{4\beta\beta_1^\prime}, \ 
b_{1,2} = \frac{ic_0\beta_1 - 2\beta_1^\prime}{4\beta\beta_1^\prime}, \ b_{2,2} = \frac{ic_0\beta_1 + 2\beta_1^\prime}{4\beta\beta_1^\prime}.
\end{equation}
\noindent Similarly, for $j=3$, taking $E_3 = c_0^2/2$, $\theta_3(0) = i$ and continuing in the same way yields
\begin{equation}
\label{eq:Apdx-A3-coef}
a_{1,3} = i\frac{-c_0\beta_1 + 2\beta_1^\prime}{4\beta\beta_1^\prime}, \ a_{2,3} = -i\frac{c_0\beta_1 + 2\beta_1^\prime}{4\beta\beta_1^\prime}, \ 
b_{1,3} = -i\frac{c_0\beta_1 + 2\beta_1^\prime}{4\beta\beta_1^\prime}, \ b_{2,3} = i\frac{-c_0\beta_1 + 2\beta_1^\prime}{4\beta\beta_1^\prime}.
\end{equation}
Remember that our aim is to compute $A_j$, which, from \eqref{eq:Apdx-T-def} and \eqref{eq:Apdx-asympts-th-varth}, implies computing 
\begin{equation*}
\lim\limits_{s\to\infty} \theta_j(s) \bar\vartheta_j(s) = \left(a_{1,j}(s) \gamma_1 + a_{2,j}(s) \gamma_2\right)
\left(\overline{b_{1,j}(s) \gamma_1 + b_{2,j}(s) \gamma_2}\right),  \ j=2,3.
\end{equation*}
Therefore, from \eqref{eq:Apdx-A2-coef},  
\begin{equation}
\label{eq:Apdx-th2-varth2-exp1}
\lim\limits_{s\to\infty} \theta_2(s) \bar\vartheta_2(s)= \left( \frac{ic_0}{4\beta_1^\prime}(\gamma_1+\gamma_2) + \frac{1}{2\beta_1} (\gamma_1-\gamma_2)\right)  \overline{\left( \frac{ic_0}{4\beta_1^\prime}(\gamma_1+\gamma_2) - \frac{1}{2\beta_1} (\gamma_1-\gamma_2) \right)} , 
\end{equation}
and by using the following identities for $y\in\RR$:
\begin{equation}
\label{eq:Apdx-Gamma-prop}
\begin{aligned}
|\Gamma (1+iy)|^2 &= y^2 | \Gamma(iy)|^2, \quad |\Gamma(iy)|^2 = \frac{\pi}{y\sinh \pi y}, \\
|\Gamma(1/2+iy)|^2 &= \frac{\pi}{\cosh \pi y}, \quad \overline{\Gamma(1+iy)} = \Gamma(1-iy),
\end{aligned}
\end{equation}
we obtain 
\begin{align*}
|\gamma_1 - \gamma_2 |^2 = 4\pi c_0^2\frac{ (1+e^{-\pi c_0^2/2})^2}{1-e^{-\pi c_0^2}}, \ |\gamma_1 + \gamma_2 |^2 = \pi c_0^2 \frac{e^{-\pi c_0^2/4}}{\sinh(\pi c_0^2/4)},  \ (\gamma_1 - \gamma_2) \overline{(\gamma_1 + \gamma_2)} = -4c_0^2 \pi.
\end{align*}
Consequently, we can write 
\begin{align}
\label{eq:A2}
A_2 &= \lim\limits_{s\to\infty} T_2(s) = \lim\limits_{s\to\infty} T_2(s) (1+\theta_2(s) \bar \vartheta_2(s)) \nonumber \\
&=  \frac{2}{\pi c_0} e^{\pi c_0^2/4} \sinh(\pi c_0^2/2) \Re \{ e^{i\pi /4} \Gamma(1-ic_0^2/4) \Gamma(1/2+ic_0^2/4) \} 
=\frac{2}{\pi c_0} e^{\pi c_0^2/4} \sinh(\pi c_0^2/2) \Re\left\lbrace\Upsilon\right\rbrace,
\end{align}
with $\Upsilon=e^{i\pi /4} \Gamma(1-ic_0^2/4) \Gamma(1/2+ic_0^2/4)= \frac{\sqrt{\pi}}{2}\ e^{i\pi/4}B(1-ic_0^2/4,1/2+ic_0^2/4)$, 
where $B(\cdot,\cdot)$ is the beta function. Finally, using \eqref{eq:Apdx-A3-coef} and following the same steps as above, we can conclude that
\begin{equation}
\label{eq:A3}
A_3 =  \frac{2}{\pi c_0} e^{\pi c_0^2/4} \sinh(\pi c_0^2/2) \Im \{ \Upsilon\}.
\end{equation}
\bibliography{VFEver2ArXiv}
\bibliographystyle{ieeetr}

\end{document}